\documentclass[onefignum,onetabnum]{siamonline171218}

\usepackage{lipsum}
\usepackage{amsfonts}
\usepackage{graphicx}
\usepackage{epstopdf}
\usepackage{algorithmic}
\usepackage{amssymb}
\usepackage{mathpazo}
\usepackage[final]{showkeys}

\ifpdf
  \DeclareGraphicsExtensions{.eps,.pdf,.png,.jpg}
\else
  \DeclareGraphicsExtensions{.eps}
\fi
\newcommand{\dd}{\mathrm{d}}
\newcommand{\srm}{\mathrm{s}}
\newcommand{\hrm}{\mathrm{h}}
\newcommand{\cc}{\mathrm{c}}
\newcommand{\Brm}{\mathrm{B}}
\newcommand{\ee}{\mathrm{e}}
\newcommand{\ii}{\mathrm{i}}
\newcommand{\Tc}{\ell}
\renewcommand{\varepsilon}{\epsilon}
\newcommand{\Ai}{\mathrm{Ai}}

\newcommand{\im}[1]{\mathrm{Im} \left( #1 \right)}
\newcommand{\R}{\mathbb{R}}

\newcommand{\Quad}[1][1]{\hspace*{#1em}\ignorespaces}
\newcommand{\bigO}[1]{O \left( #1 \right)}

\usepackage{enumitem}
\setlist[enumerate]{leftmargin=.5in}
\setlist[itemize]{leftmargin=.5in}

\newsiamremark{remark}{Remark}
\newsiamremark{hypothesis}{Hypothesis}
\crefname{hypothesis}{Hypothesis}{Hypotheses}
\newsiamthm{claim}{Claim}
\newsiamthm{example}{Example}

\headers{Benjamin-Ono Universality}{E. Blackstone, P. D. Miller, and M. Mitchell}

\title{Universality in the Small-Dispersion Limit of the Benjamin-Ono Equation}

\author{Elliot Blackstone\thanks{Department of Mathematics, University of Michigan, Ann Arbor, MI 
  (\email{eblackst@umich.edu}).}
\and Peter D. Miller\thanks{Department of Mathematics, University of Michigan, Ann Arbor, MI 
  (\email{millerpd@umich.edu}).}
\and Matthew D. Mitchell\thanks{Department of Physics, University of Michigan, Ann Arbor, MI (\email{mattdmit@umich.edu}).}
}

\usepackage{amsopn}

\makeatletter
\newcommand*{\addFileDependency}[1]{
  \typeout{(#1)}
  \@addtofilelist{#1}
  \IfFileExists{#1}{}{\typeout{No file #1.}}
}
\makeatother

\newcommand*{\myexternaldocument}[1]{%
    \externaldocument{#1}%
    \addFileDependency{#1.tex}%
    \addFileDependency{#1.aux}%
}

\ifpdf
\hypersetup{
  pdftitle={Zero-Dispersion Universality for Benjamin-Ono},
  pdfauthor={Elliot Blackstone and Peter D. Miller and Matthew Mitchell}
}
\fi

\myexternaldocument{ex_supplement}

\usepackage{comment}
\usepackage{tikz}
\usetikzlibrary{arrows}
\usetikzlibrary{decorations.pathmorphing}
\usetikzlibrary{decorations.markings}
\usetikzlibrary{patterns}
\usetikzlibrary{automata}
\usetikzlibrary{positioning}
\usetikzlibrary{fadings}
\usetikzlibrary{calc}
\usepackage{tikz-cd}

\tikzset{->-/.style={decoration={markings,mark=at position #1 with {\arrow[thick]{>}}},postaction={decorate}}}
	
\tikzset{-<-/.style={decoration={markings,mark=at position #1 with{\arrow[thick]{<}}},postaction={decorate}}}

\tikzset{
  frac arrow/.style={postaction={decorate,decoration={
        markings,
        mark=at position #1 with {\arrow[]{Stealth[round]}}
      }}},
}

\tikzset{
  mid arrow/.style={postaction={decorate,decoration={
        markings,
        mark=at position .5 with {\arrow[#1]{Stealth[round]}}
      }}},
}

\tikzset{
  midp arrow/.style={postaction={decorate,decoration={
        markings,
        mark=at position .6 with {\arrow[#1]{Stealth[round]}}
      }}},
}

\tikzset{
  midpp arrow/.style={postaction={decorate,decoration={
        markings,
        mark=at position .8 with {\arrow[#1]{Stealth[round]}}
      }}},
}
\tikzset{
  midm arrow/.style={postaction={decorate,decoration={
        markings,
        mark=at position .4 with {\arrow[#1]{Stealth[round]}}
      }}},
}

\tikzset{
  midmm arrow/.style={postaction={decorate,decoration={
        markings,
        mark=at position .3 with {\arrow[#1]{Stealth[round]}}
      }}},
}

\newtheorem{conjecture}[theorem]{Conjecture}

\begin{document}

\maketitle

\begin{abstract}
  We examine the solution of the Benjamin-Ono Cauchy problem for rational initial data in three types of double-scaling limits in which the dispersion tends to zero while simultaneously the independent variables either approach a point on one of the two branches of the caustic curve of the inviscid Burgers equation, or approach the critical point where the branches meet.  The results reveal universal limiting profiles in each case that are independent of details of the initial data.  We compare the results obtained with corresponding results for the Korteweg-de Vries equation found by Claeys-Grava in three papers \cite{ClaeysG09,ClaeysG10a,ClaeysG10b}.   Our method is to analyze contour integrals appearing in an explicit representation of the solution of the Cauchy problem, in various limits involving coalescing saddle points.  
\end{abstract}



\section{Introduction}
We consider in this paper the Benjamin-Ono equation
\begin{equation}
\partial_t u + \partial_x(u^2) = \epsilon \partial_x|D_x|u,\quad x\in\mathbb{R},
\quad t\ge 0,
\quad \epsilon>0,
\label{eq:BO}
\end{equation}
with dispersive term defined either by the Fourier multiplier $\widehat{|D_x|u}(\xi)=|\xi|\widehat{u}(\xi)$  or equivalently by the singular integral
\begin{equation}|D_x|f(x):=-\frac{1}{\pi}{\mathrm{P.V.}}\int_\mathbb{R}{(x-y)^2}.
\frac{f'(y)\,\dd y}{y-x},
\label{eq:AbsDx}
\end{equation}
for a real-valued function $u=u(t,x;\epsilon)$ subject to a real-valued rational initial condition of the form
\begin{equation}
u(0,x;\epsilon)=u_0(x):=\sum_{n=1}^N\left(\frac{c_n}{x-p_n}+\frac{{c_n^*}}{x-{p_n^*}}\right),
\quad \operatorname{Im}(p_n)>0,\quad \text{$p_m\neq p_n$ if $m\neq n$}.
\label{eq:rationalIC}
\end{equation}
This initial-value problem arises in the study of internal waves in two-layer stratified fluids as an asymptotic model for long small-amplitude waves propagating in one direction when the lower fluid layer is idealized to be infinitely deep.   When $\epsilon>0$ is small, it appears feasible to neglect the dispersive term entirely, however the resulting initial-value problem for the inviscid Burgers equation generally does not have a global solution $u^\mathrm{B}(t,x)$.  Typically, there is a positive least upper bound $t_\mathrm{c}>0$  for the time of existence of $u^\mathrm{B}(t,x)$, and we assume in this paper that there is a unique point $x=x_\cc$ at which the solution breaks down via a generic gradient catastrophe.
\begin{definition}[generic gradient catastrophe] 
$(t_\cc,x_\cc)\in\mathbb{R}^2$ is a generic gradient catastrophe point if the solution $u^\mathrm{B}(t,x)$ is classical for $t<t_\cc$ and for $t=t_\cc$ but $x\neq x_\cc$, and if for $(t,x)=(t_\mathrm{c},x_\mathrm{c})$ the equation $x=2u_0(y)t+y$ for intercepts $y$ of characteristic lines has a root of multiplicity exactly $3$.  \label{def:gc}
\end{definition}
The terminology means that the Burgers solution $u^\mathrm{B}(t,x)$ is finite with value $u_\cc:=u^\mathrm{B}(t_\cc,x_\cc)$ at the catastrophe point, while $u_x^\mathrm{B}(t_\mathrm{c},x_\mathrm{c})=-\infty$ and when $t=t_\cc$ the inverse function $u^\mathrm{B}\mapsto x$ has a root of exactly third order at $u_\cc$.

For times $t>t_\mathrm{c}$, the formally small dispersive term $\epsilon \partial_x|D_x|u$ in \eqref{eq:BO} is no longer negligible for $x$ in a certain interval that expands in length as $t$ increases.  Within this interval, the solution $u(t,x;\epsilon)$ displays oscillations of wavelength proportional to $\epsilon$ and amplitude that does not decay to zero, forming a so-called \emph{dispersive shock wave}.  If it is assumed that a dispersive shock wave forms with microstructure modeled by a periodic traveling wave solution of the Benjamin-Ono equation, then a version of Whitham's modulation theory applies to yield modulation equations describing the space-time evolution of the slowly modulated wave parameters \cite{DobrokhotovK91}.  It turns out that these equations admit Riemann invariants even in the case of multiphase modulation, as is also known to be true for the Korteweg-de Vries (KdV) equation \cite{FlaschkaFM80}.  However unlike the KdV case in which the characteristic velocities are complicated nonlinear functions of the Riemann invariants, for Benjamin-Ono the modulation equations in Riemann-invariant form take the form of an odd number of \emph{uncoupled} copies of the inviscid Burgers equation $u^\mathrm{B}_t + 2u^\mathrm{B}u^\mathrm{B}_x=0$.  Thus one prediction of this formal theory is that the dispersive shock wave formed in the Benjamin-Ono equation occupies exactly the same space-time domain on which the implicit solution of the approximating inviscid Burgers equation for the initial condition $u_0(x)$ is multi-valued.  Working within this formalism, one can also couple the dispersive shock wave with the slowly-varying background at its left and right boundaries to describe global dynamics beyond the gradient catastrophe point \cite{JorgeMS98,Matsuno98b,Matsuno98a}.

Several results have been obtained going beyond the Whitham modulation theory formalism. With motivation from the Lax-Levermore method for KdV \cite{LL}, which is based on the idea of a \emph{soliton ensemble} that approximates the given initial data by a reflectionless potential consisting of $O(\epsilon^{-1})$ solitons, a weak convergence result for Benjamin-Ono soliton ensembles was proved in \cite{MillerX10}.  With the help of a new formula for the solution of Benjamin-Ono, the same limiting formula was recently shown to be valid for the exact Cauchy problem for \eqref{eq:BO} (i.e., without the soliton ensemble approximation) with $L^2(\mathbb{R})$ initial data \cite{Gerard23}, and it has also been proven that this formula also governs the weak limit in the periodic setting \cite{Gassot23}.  This weak limit formula is particularly easy to explain; letting $u_0^\mathrm{B}(t,x)<u_1^\mathrm{B}(t,x)<\cdots<u_{2J}^\mathrm{B}(t,x)$ denote the $2J+1$ branches of the generally multivalued solution of the inviscid Burgers equation obtained by the method of characteristics (so $J(0,x)=0$ for all $x$ and generally the integer $J=J(t,x)$ depends on $(t,x)$), one defines the alternating sum
\begin{equation}
    \bar{u}(t,x):=\sum_{n=0}^{2J(t,x)}(-1)^nu_n^{\mathrm{B}}(t,x).
\end{equation}
Then all of the aforementioned papers prove the convergence of $u(t,x;\epsilon)$ to $\bar{u}(t,x)$ in the weak $L^2$ topology as $\epsilon\to 0$, where $u(t,x;\epsilon)$ denotes the solution of the Cauchy problem \eqref{eq:BO} in either the periodic \cite{Gassot23} or whole line \cite{Gerard23} settings, or in the latter case its soliton ensemble approximation \cite{MillerX10}.  Some partial results on strengthening the convergence in the soliton ensemble setting were obtained in \cite{BGM23}.   In a recent breakthrough, the new formula used to prove weak convergence in \cite{Gerard23} was adapted to general rational initial data of the form \eqref{eq:rationalIC}, resulting in a representation for $u(t,x;\epsilon)$ as the real part of a ratio of $(N+1)\times (N+1)$ determinants with entries that are contour integrals explicitly depending on $(t,x;\epsilon)$ \cite{RationalBO}; this formula was then used in \cite{SmallDispersionBO-1} to prove a locally uniform $O(\epsilon)$ estimate on the difference between $u(t,x;\epsilon)$ and an explicit and generally highly oscillatory approximation $u^\mathrm{ZD}(t,x;\epsilon)$  that is modeled on the Dobrokhotov-Krichever multiphase wave formula \cite{DobrokhotovK91}.

\begin{figure}
    \centering
    \includegraphics[width=0.5\linewidth]{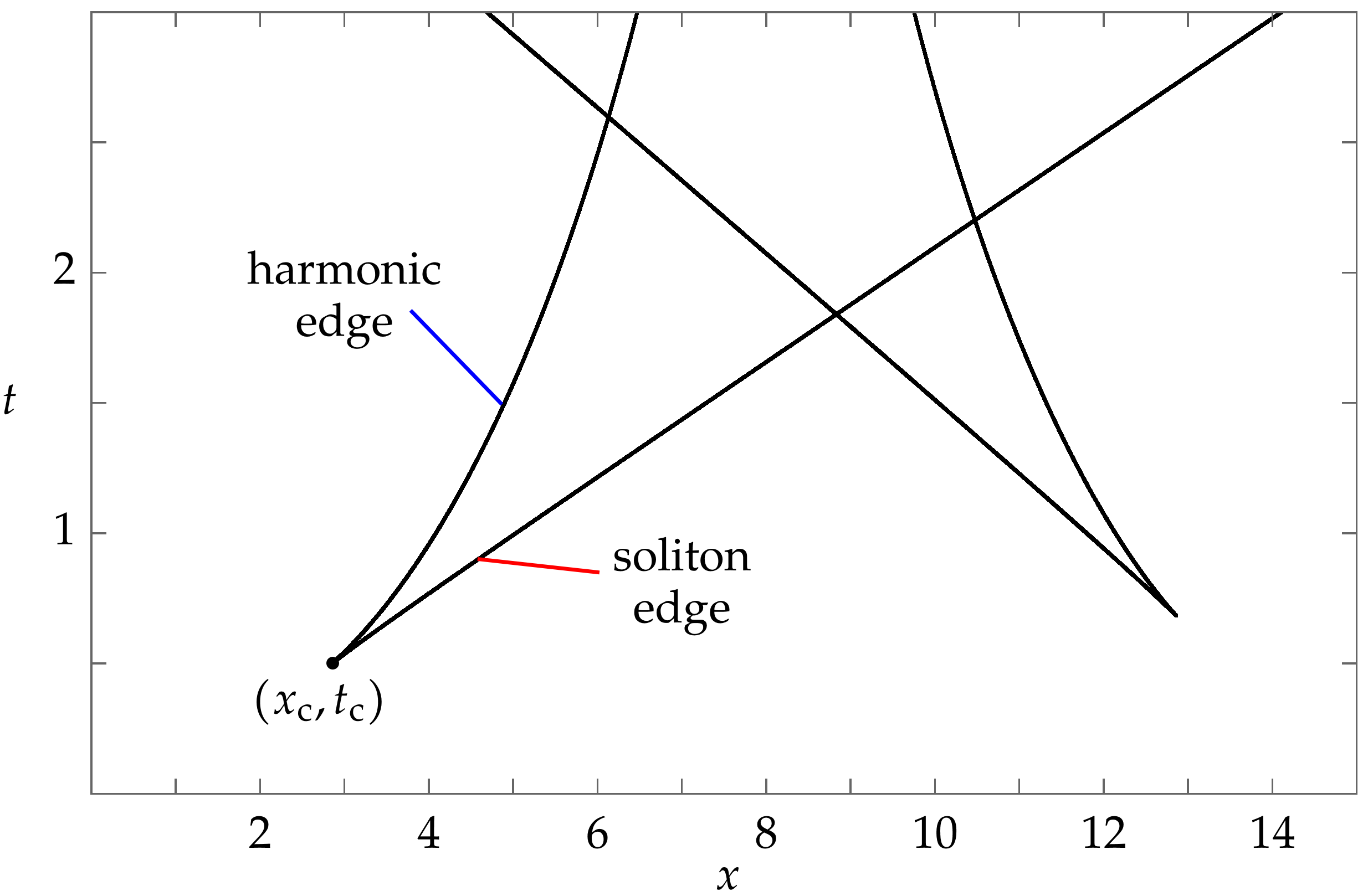}
    \caption{The caustic locus corresponding to the initial condition $u_0(x)$ of the form~\eqref{eq:rationalIC} for $N=2$ with $c_1=1-\ii$, $c_2=1+\ii/\sqrt{2}$, $p_1=\ii$, $p_2=16+\ii$. We can see the harmonic and soliton edges emerging from the critical point $(x_\cc,t_\cc)$.}
    \label{fig:edge}
\end{figure}

The locally uniform approximation by $u^\mathrm{ZD}(t,x;\epsilon)$ is valid on compact subsets of spacetime that are disjoint from the \emph{caustic locus} $C$ consisting of bifurcation points for the implicit solution of the inviscid Burgers equation:  $u^\mathrm{B}=u_0(x-2tu^\mathrm{B})$.  In a neighborhood of the critical point $(t_\cc,x_\cc)$, $C$ consists of two curves $x=x_\hrm(t)$ and $x=x_\srm(t)$ defined for $t>t_\cc$ with $x_\hrm(t)<x_\srm(t)$ and such that $x_\hrm(t_\cc)=x_\srm(t_\cc)=x_\cc$, see Figure~\ref{fig:edge}. We refer to the caustic curves represented by $x=x_\hrm(t)$ and $x=x_\srm(t)$ as the harmonic and soliton edges respectively, due to the demonstrated behavior of the solution $u(t, x; \epsilon)$ near these caustic curves. These results are of a similar nature to what is expected based on Whitham modulation theory for other dispersive regularizations of the inviscid Burgers equation, and what has been rigorously proven \cite{DeiftVZ97} in some integrable cases such as the Korteweg-de Vries equation.  It should be pointed out that in the Korteweg-de Vries case and other cases different from Benjamin-Ono, the caustic locus $C$ consists of different curves that have nothing to do with the multivalued solution of the invsicid approximation, although those curves still emerge from the same critical point $(t_\cc,x_\cc)$.  These curves generally come from matching conditions at the two edges of the dispersive shock wave and are described by a degeneration of the single-phase modulation equations provided by Whitham theory.

In this paper, we give new asymptotic formul\ae\ for the solution of the Benjamin-Ono equation \eqref{eq:BO} with arbitrary rational initial data $u_0$ independent of $\epsilon$ of the form~\eqref{eq:rationalIC} that are valid in the limit that $\epsilon\to 0$ while simultaneously the spacetime coordinates are localized at a suitable scale to a point on the caustic locus $C$.  In this situation, the limiting behavior exhibits universal features that are mainly independent of initial data.  These results are analogues of those obtained by Claeys and Grava for the Korteweg-de Vries equation in a series of three well-known papers \cite{ClaeysG09,ClaeysG10a,ClaeysG10b}.  We will describe our results and compare them with those for the Korteweg-de Vries equation now.

\subsection{Asymptotic behavior near the soliton edge}

\subsubsection{Formulation of main results}

Assume that $(t_\cc,x_\cc)$ is a generic gradient catastrophe point (see Definition~\ref{def:gc}). 
Then, when $x=x_\srm(t)$ for $t-t_\cc$ positive and sufficiently small there are two distinct real solutions $y=y_\srm(t)$ and $y=y_0(t)$ with $y_\srm(t)<y_0(t)$ of the characteristic equation $x=2u_0(y)t+y$, and $y_\srm(t)$ is a double root while $y_0(t)$ is a simple root.
\begin{theorem}[Benjamin-Ono universality at the soliton edge]
\label{thm:soliton-edge}
    Let $u_0$ be a rational initial condition of the form~\eqref{eq:rationalIC}, and assume that $(t_\cc,x_\cc)$ is a generic gradient catastrophe point and moreover that all complex roots $y$ of the characteristic equation $x_\cc=2u_0(y)t_\cc+y$ are simple. 
    Set $u_\srm^\Brm(t) := u_0(y_\srm(t))$, $u_0^\Brm(t) := u_0(y_0(t))$ and $k_\srm(t) := 2 t [-\frac{1}{2}u_0''(y_\srm(t))]^{1/3} > 0$. 
    Then, there exists $\delta>0$ sufficiently small that for $t\in K:=(t_\cc,t_\cc+\delta)$,
    \begin{equation}
        u(t, x; \epsilon) = u^\Brm_0(t) + L\left( \epsilon^{-1/3} k_\srm(t) \frac{\Ai(X_\srm)}{\Ai'(X_\srm)} +  2 t \lambda(t) ; \, c(t)\right) + O(\epsilon^{1/3}) 
        \label{eq:BOSolitonApproxFormula}
    \end{equation}
    as $\epsilon\to 0$, where
    \begin{equation}
        L(\xi; c) := \frac{2 c}{ c^2 \xi^2+1},\quad \xi\in\mathbb{R},\quad c>0
        \label{eq:BOSolitonProfile}
    \end{equation}
    is the Benjamin-Ono soliton Lorentzian profile,
    \begin{equation}
        c(t):=2\left(u_\srm^\Brm(t)-u_0^\Brm(t)\right),
        \label{eq:SolitonSpeed}
    \end{equation}
    and
    \begin{equation}
        \lambda(t) := \frac{1}{2\pi} \int_0^{+\infty} \left[ \frac{\partial_yg_\srm(y_\srm(t) + s, t)}{g_\srm(y_\srm(t) + s, t)} - \frac{\partial_yg_\srm(y_\srm(t) - s, t)}{g_\srm(y_\srm(t) - s, t)} \right] \frac{\dd s}{s},
        \label{eq:IntroSolGammaDef}
    \end{equation}
    where $g_\srm(y, t)$ is the positive rational function defined for $y \in \R$ and $t \in K$ by
    \begin{equation}
        g_\srm(y, t) := \frac{y + 2t u_0(y) - x_\srm(t)}{\big( y + 2t u_0\big( y_0(t) \big) - x_\srm(t) \big) \big(y + 2t u_0\big( y_\srm(t) \big) -x_\srm(t) \big)^2}.
        \label{eq:IntroSolgDef}
    \end{equation}
    The approximate formula \eqref{eq:BOSolitonApproxFormula} holds uniformly for $X_\srm$ bounded in $\mathbb{R}$ where
    \begin{equation}
        X_\srm := \frac{x-x_\srm(t)}{k_\srm(t) \epsilon^{2/3}}.
        \label{eq:Xs-coordinate}
    \end{equation}
    and $t$ in compact subsets of $K$.
\end{theorem}

\begin{remark}
    \label{rem:OtherSolitonEdges}
    For a particular initial condition $u_0$, instead of only considering the soliton edge curve emanating from the first gradient catastrophe point $(t_\cc, x_\cc)$, we can extend the definition of the soliton edge to the locus of points $(t, x)$ where the characteristic equation $x=2u_0(y)t+y$ has exactly one real and double root to the left of one real and simple root. This locus of points can always be decomposed into a finite union of curves $\{ (t, x^i_\srm(t)) \, | \, t \in I_\srm^i \}$ where $I_\srm^i$ is an open (possibly semi-infinite) interval of $\R_+$ and $x^i_\srm$ is a differentiable function on $I_\srm^i$. Each of these constituent curves can be viewed as a soliton edge curve. By careful consideration of the conditions required for the proof, we note that the asymptotic formula in Theorem \ref{thm:soliton-edge} above follows on soliton edges more generally with uniform error for $t \in K$, a compact subset of $I_\srm^i$ for a particular soliton edge $x = x^i_\srm(t)$ on which all of the non-real roots of the characteristic equation $x^i_\srm(t)=2u_0(y)t+y$ remain simple.  For example, one can consider the right-most curve emanating from the second cusp point on $C$ near $x\approx 13$ visible in Figure~\ref{fig:edge} as another soliton edge curve.
\end{remark}

\begin{corollary}
    Under the same conditions as Theorem~\ref{thm:soliton-edge}, 
\begin{equation}
    u(t,x;\epsilon)=u_0^\Brm(t)+\sum_{n=1}^\infty L\left(\epsilon^{-1/3}k_\srm(t)(X_\srm-a_n)+ 2t\lambda(t);c(t)\right) + O(\epsilon^{1/3}),
    \label{eq:u-soliton-edge-sum}
\end{equation}
holds, wherein  $a_1>a_2>a_3>\cdots$ denote the zeros of the Airy function $\Ai(\cdot)$.  Also, define for each $n=1,2,3,\dots$ an offset
\begin{equation}
    x_0^{(n)}(t;\epsilon):=x_\srm(t) + \epsilon^{2/3}k_\srm(t)a_n-2\epsilon t\lambda(t).
    \label{eq:soliton-offset}
\end{equation}
Then, given any $t_0\in K=(t_\cc,t_\cc+\delta)$, and any $n=1,2,3,\dots$,
\begin{equation}
    u(t,x;\epsilon)=u_0^\Brm(t_0) + L\left(\epsilon^{-1}\left(x-x_0^{(n)}(t_0;\epsilon)-2u_\srm^\Brm(t_0)(t-t_0)\right);c(t_0)\right) + O(\epsilon^{1/3})
    \label{eq:u-single-soliton}
\end{equation}
holds as $\epsilon\to 0$ uniformly for $x-x_0^{(n)}(t_0;\epsilon)=o(\epsilon^{2/3})$ and $t-t_0=o(\epsilon^{2/3})$.  
\label{cor:soliton-edge}
\end{corollary}
The exact soliton solutions of the Benjamin-Ono equation in the form \eqref{eq:BO} are
\begin{equation}
    u(t,x;\epsilon)=L(\epsilon^{-1}(x-x_0-c(t-t_0));c),
\end{equation}
and, by a Galilean boost, we obtain a soliton propagating on a constant background level $b$:
\begin{equation}
    u(t,x;\epsilon)=b+L(\epsilon^{-1}(x-x_0-(c+2b)(t-t_0));c).
    \label{eq:soliton-on-background}
\end{equation}
The explicit terms on the right-hand side of \eqref{eq:u-single-soliton} therefore constitute an exact soliton solution of the form \eqref{eq:soliton-on-background} with parameters
\begin{equation}
    b:=u_0^\Brm(t_0),\quad c:=c(t_0)=2(u_\srm^\Brm(t_0)-u_0^\Brm(t_0)),\quad x_0:=x_0^{(n)}(t_0;\epsilon).
\end{equation}
Note that when $t=t_0$, the soliton is located at $x=x_0^{(n)}(t_0;\epsilon)$, which is a point displaced from the curve $x=x_\srm(t_0)$ by a small constant $-2\epsilon t_0\lambda(t_0)$ and by the location of the $n^\mathrm{th}$ zero of the Airy function $a_n$, scaled by the factor $\epsilon^{2/3}k_\srm(t_0)$.

\begin{figure}[h]
\begin{center}
    \includegraphics[width=0.49\linewidth]{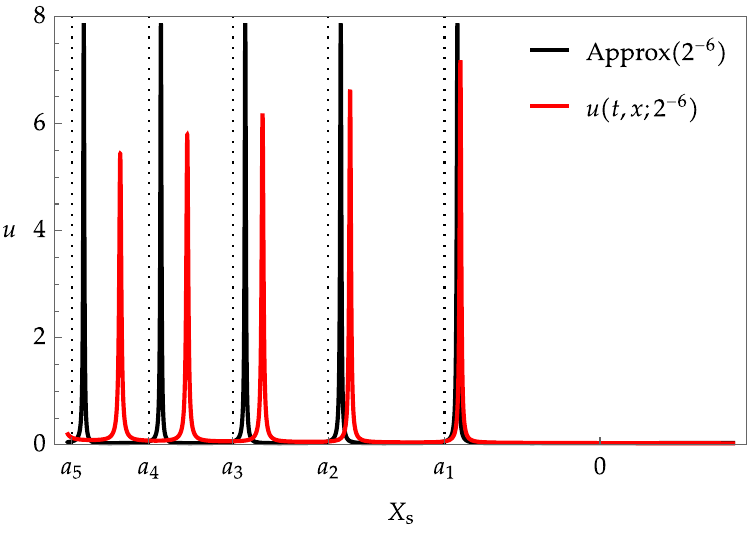}\hfill%
    \includegraphics[width=0.49\linewidth]{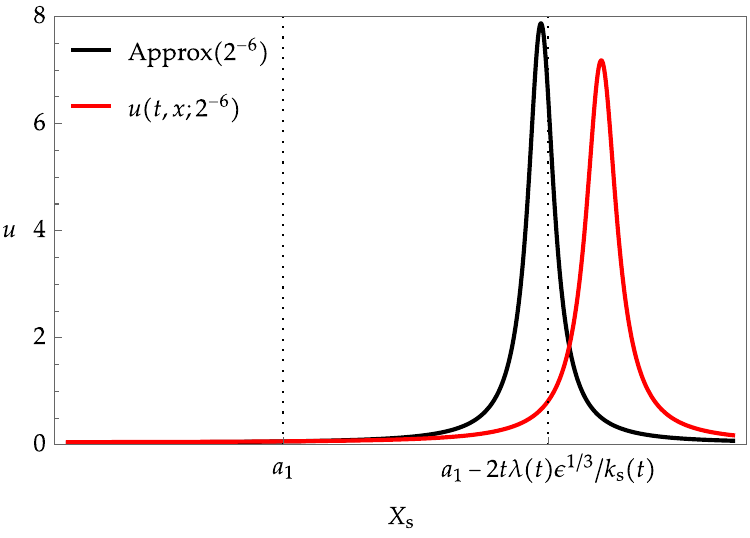}
\end{center}
\caption{Left:  comparison between the exact solution $u(t,x;\epsilon)$ with initial data $u_0(x)=2/(1+x^2)$ for $t=2>t_\cc$ (red; for a discussion of this computation, see \cite[Section 8]{SmallDispersionBO-1}) and the explicit terms on the right-hand side of \eqref{eq:BOSolitonApproxFormula} (black) for $\epsilon=2^{-6}$, plotted as a function of $X_\srm$.  Right:  the same zoomed in near the leading pulse in the train to the smaller interval $-2.5<X_\srm<-2$.  The Airy zeros $a_1,a_2,a_3,\dots$ are shown on the $X_\srm$-axis and with vertical dotted lines for reference.}
\label{fig:soliton-edge-6}
\end{figure}

\begin{figure}[h]
\begin{center}
    \includegraphics[width=0.49\linewidth]{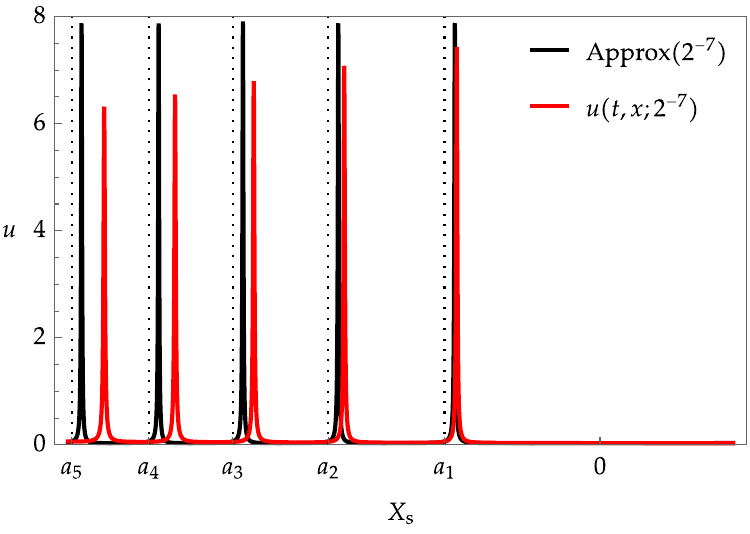}\hfill%
    \includegraphics[width=0.49\linewidth]{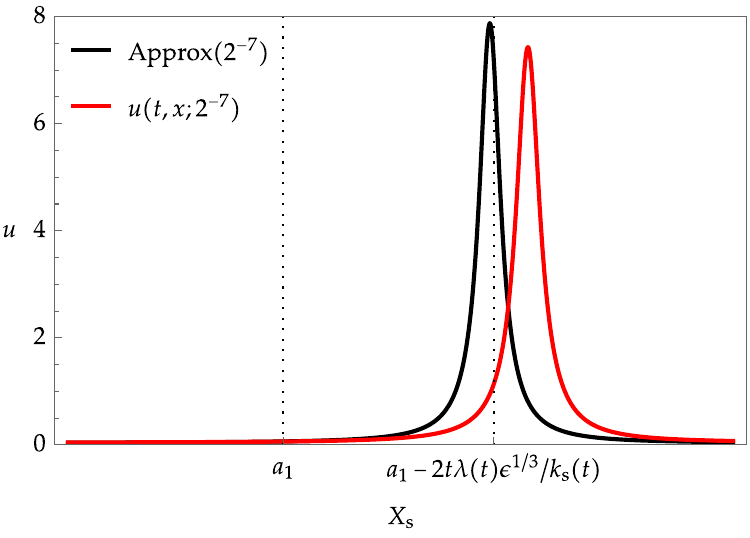}
\end{center}
\label{fig:soliton-edge-7}
\caption{The same as Figure~\ref{fig:soliton-edge-6} but for $\epsilon=2^{-7}$.}
\end{figure}

\begin{figure}[h]
\begin{center}
    \includegraphics[width=0.49\linewidth]{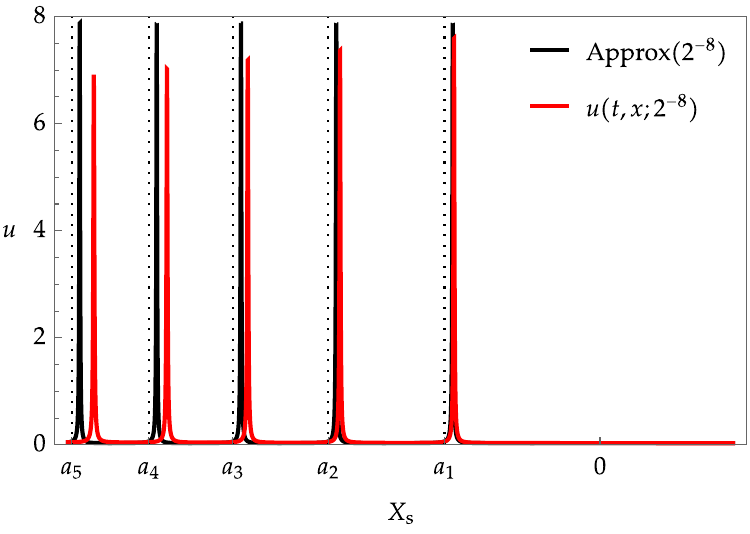}\hfill%
    \includegraphics[width=0.49\linewidth]{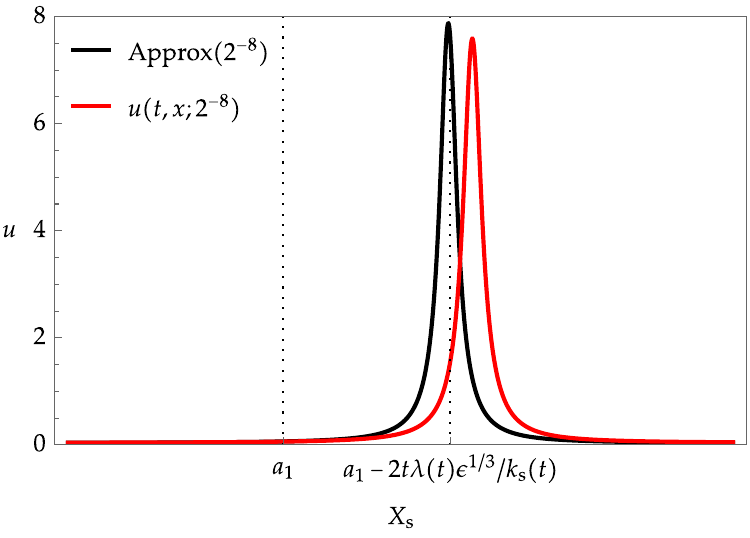}
\end{center}
\caption{The same as Figure~\ref{fig:soliton-edge-6} but for $\epsilon=2^{-8}$.}
\label{fig:soliton-edge-8}
\end{figure}

\begin{figure}[h]
\begin{center}
    \includegraphics[width=0.49\linewidth]{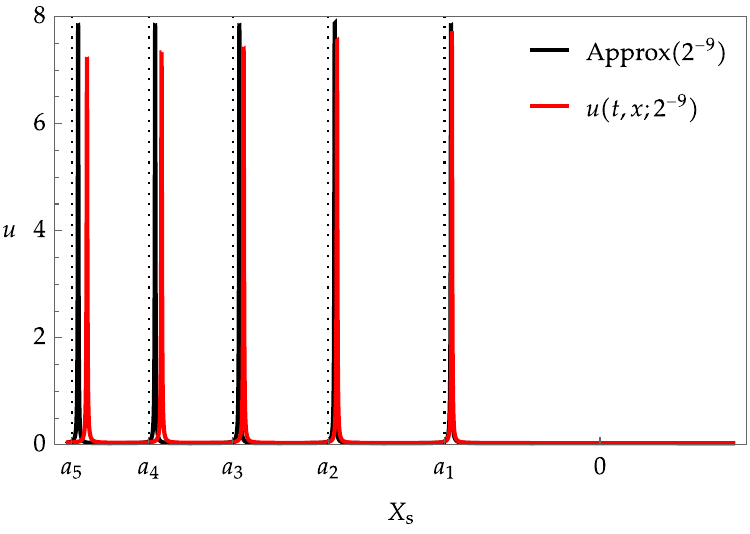}\hfill%
    \includegraphics[width=0.49\linewidth]{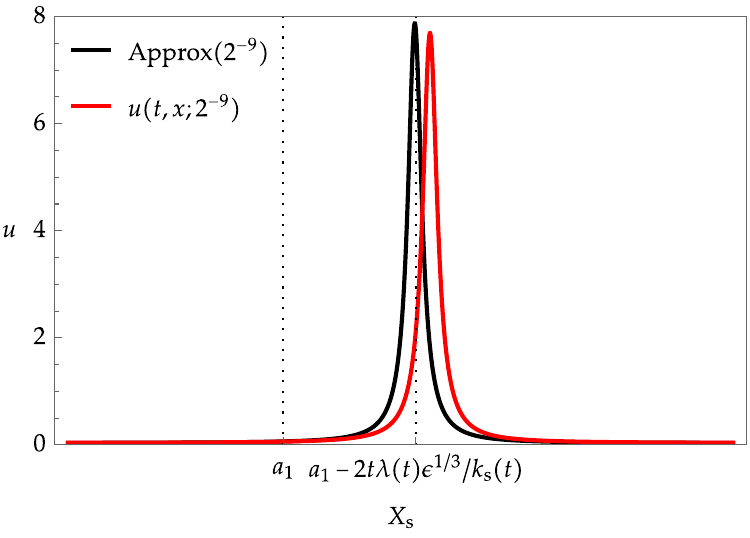}
\end{center}
\caption{The same as Figure~\ref{fig:soliton-edge-6} but for $\epsilon=2^{-9}$.}
\label{fig:soliton-edge-9}
\end{figure}

The accuracy of the approximate formula \eqref{eq:BOSolitonApproxFormula} is illustrated in Figures~\ref{fig:soliton-edge-6}--\ref{fig:soliton-edge-9}.  From these plots, one can see that when $\epsilon$ decreases by a factor of $8$, the maximum absolute error over the fixed interval of $X_\srm$ in the plots appears to decrease by a factor of about $2$, consistent with the $O(\epsilon^{1/3})$ error estimate in \eqref{eq:BOSolitonApproxFormula}.

\subsubsection{Comparison with KdV}
\label{sec:soliton-comparison}
In their paper \cite{ClaeysG10b}, Claeys and Grava proved an analogue of this result for the KdV equation, which we write here in the form\footnote{We write the KdV equation in the form \eqref{eq:KdV} so that it has the same inviscid Burgers formal limit as the Benjamin-Ono equation in the form \eqref{eq:BO}.  Since in the papers \cite{ClaeysG09,ClaeysG10a,ClaeysG10b} the KdV equation is written in the form $\partial_t u^\mathrm{CG}+6u^\mathrm{CG}\partial_xu^\mathrm{CG} +\epsilon^2\partial_x^3u^\mathrm{CG}=0$, we make the identification that $u=3u^\mathrm{CG}$ and keep the same independent variables $(t,x)$.  All of the formul\ae\ we quote from the papers \cite{ClaeysG09,ClaeysG10a,ClaeysG10b} will be modified accordingly.}
\begin{equation}
    \partial_tu +2u\partial_x u + \epsilon^2\partial_x^3u=0.
    \label{eq:KdV}
\end{equation}
We now describe the result of \cite{ClaeysG10b}, which was obtained
under the assumption of real-analytic initial data $u_0(x)$ with sufficient decay in a region containing the real line, and such that $u_0$ is a negative bump with a single critical point (minimizer), for which the inviscid Burgers equation generates a generic gradient catastrophe (see \cite[Assumptions 1.1]{ClaeysG10b} for details).

Although for the same initial data $u_0(x)$, both the Benjamin-Ono equation \eqref{eq:BO} and the KdV equation 
share the same gradient catastrophe point $(t_\cc,x_\cc)$ because both equations are dispersive perturbations of the same inviscid Burgers equation $u_t^\Brm+2u^\Brm u^\Brm_x=0$, as mentioned above, the caustic locus consists of different curves for KdV than for Benjamin-Ono.  For KdV, the analogue of the soliton-edge curve $x=x_\srm(t)$ is a curve $x=x^\mathrm{KdV}_\srm(t)$ determined via degeneration of the Whitham equations along with auxiliary functions $U_\srm(t)$ and $V_\srm(t)$ by the system of equations
\begin{equation}
    \begin{split}
        2tU_\srm(t)+f_\mathrm{L}(U_\srm(t))&=x^\mathrm{KdV}_\srm(t)\\
        2t + \theta_\srm(V_\srm(t),U_\srm(t))&=0\\
        \int_{U_\srm(t)}^{V_\srm(t)}\left(2t+\theta_\srm(\zeta,U_\srm(t))\right)\sqrt{\zeta-U_\srm(t)}\,\dd\zeta &=0
    \end{split}
\end{equation}
where $f_\mathrm{L}(U)$ denotes the inverse function of the decreasing part of the initial condition $u_0$ and where for $V>U$,
\begin{equation}
    \theta_\srm(V,U):=\frac{1}{2\sqrt{V-U}}\int_U^V\frac{f_\mathrm{L}'(\xi)\,\dd\xi}{\sqrt{V-\xi}}.
\end{equation}
See \cite[Eqns.\@ (1.11)--(1.14)]{ClaeysG10b}.  Like $u_0^\Brm(t)$, $U_\srm(t)$ is the continuation of the inviscid Burgers solution from the right, here evaluated on the KdV soliton-edge curve $x=x^\mathrm{KdV}_\srm(t)$.  The KdV soliton profile is the function
\begin{equation}
    S(\xi;c):=\frac{3c}{2}\mathrm{sech}^2\left(\frac{\sqrt{c}}{2}\xi\right),\quad \xi\in\mathbb{R},\quad c>0.
\end{equation}
The KdV solitons are the exact solutions of \eqref{eq:KdV} given by
\begin{equation}
    u(t,x;\epsilon)=S(\epsilon^{-1}(x-x_0-c(t-t_0));c),
\end{equation}
or with a Galilean boost,
\begin{equation}
    u(t,x;\epsilon)=b+S(\epsilon^{-1}(x-x_0-(c+2b)(t-t_0));c).
\end{equation}
Setting
\begin{equation}
    c^{\mathrm{KdV}}(t):=\frac{4}{3}\left(V_\srm(t)-U_\srm(t)\right),
\end{equation}
a rescaled and recentered spatial variable is defined by
\begin{equation}
    X_\srm^\mathrm{KdV}:=\sqrt{c^\mathrm{KdV}(t)}\frac{x-x^\mathrm{KdV}_\srm(t)}{\epsilon\ln(\epsilon^{-1})}.
\end{equation}
Defining the leading coefficients of the normalized Hermite polynomials by
\begin{equation}
    h_k:=\frac{2^{k/2}}{\pi^{1/4}\sqrt{k!}},\quad k=0,1,2,\dots,
\end{equation}
a shift (the KdV analogue of the Airy zeros) is defined by:
\begin{multline}
    a_n^\mathrm{KdV}(t;\epsilon):=-n-\frac{1}{2}\\
    {}-\frac{1}{\ln(\epsilon^{-1})}\left[\ln(2\pi h_{n-1}^2)+(2n-1)\ln\left(\frac{4}{3^{1/4}}(V_\srm(t)-U_\srm(t))^{5/4}\sqrt{-\partial_V\theta_\srm(V_\srm(t),U_\srm(t))}\right)\right].
\end{multline}
Then, the main result of \cite{ClaeysG10b} is that
\begin{equation}
    u(t,x;\epsilon)=U_\srm(t)+\sum_{n=1}^NS\left(\frac{\ln(\epsilon^{-1})}{\sqrt{c^\mathrm{KdV}(t)}}(X_\srm^\mathrm{KdV}-a^\mathrm{KdV}_n(t;\epsilon));c^\mathrm{KdV}(t)\right)+O(\epsilon\ln(\epsilon^{-1})^2)
    \label{eq:KdV-soliton-edge}
\end{equation}
holds as $\epsilon\to 0$, uniformly for $X_\srm^\mathrm{KdV}$ bounded, where the upper index $N$ is finite but sufficiently large given the bound on $X_\srm^\mathrm{KdV}$.

Like Corollary~\ref{cor:soliton-edge}, this result describes the wave field $u(t,x;\epsilon)$ near the soliton edge of the dispersive shock wave as a sum of identical solitons propagating on a background that is slowly-varying in time.  When re-expanded in suitable local coordinates anchored to a point on the soliton edge curve, each soliton in the sum is an exact solution of the relevant equation. 

In both cases, the solitons have a characteristic width in the $x$ coordinate proportional to $\epsilon$, and the spacing between the solitons is large by comparison while still small in the absolute sense.  However, there are important differences as well:
\begin{itemize}
\item
    for the Benjamin-Ono equation the spacing is proportional to $\epsilon^{2/3}$;
    \item 
    for the KdV equation the spacing is proportional to $\epsilon\ln(\epsilon^{-1})$.
\end{itemize}
Hence one quantitative difference between the two cases is that the soliton train near the edge is more densely packed for KdV than for Benjamin-Ono.  This can be seen as a rigorous explanation for the following heuristic.  In the case of initial data $u_0$ that is a positive bump, the KdV direct scattering problem reduces to the spectral problem for a semiclassical Schr\"odinger operator with a potential well, in which case the Weyl law predicts a density of states that is finite and positive for the most extreme eigenvalues.  These eigenvalues correspond to the fastest-moving solitons and they would be expected to emerge from the dispersive shock wave at the soliton edge.  For the Benjamin-Ono equation the corresponding density of eigenvalues was calculated using trace formul\ae\ by Matsuno \cite{Matsuno81}, and it vanishes for the most extreme eigenvalues, which again correspond to the fastest-moving solitons.  Hence solitons near the edge of the dispersive shock wave would be expected to be more rare than for the KdV equation with the same initial data.

Both the Benjamin-Ono formula \eqref{eq:u-soliton-edge-sum} and the KdV  formula \eqref{eq:KdV-soliton-edge} describe a similar important effect, whereby the train of solitons becomes more densely packed as the observation point penetrates further into the dispersive shock wave.  For the Benjamin-Ono equation, this effect is due to the increasing density of the Airy zeros $a_n\sim -n^{2/3}$ as $n$ increases.  However for KdV the effect is more subtle:  indeed at the leading order in $\epsilon$ the solitons are equally spaced as $a_n^\mathrm{KdV}(\epsilon)\sim -n$ for large $n$.  So one has to look at the correction term in $a_n^\mathrm{KdV}(\epsilon)$ inversely proportional to $\ln(\epsilon^{-1})$.  Because of the factorial in the Hermite constants $h_{n-1}$, for sufficiently large $n$ the correction term is positive and increasing with $n$, hence shaving off more and more of the basic spacing of $1$ as $n$ increases.  This effect is the reason why one must only use a finite sum of solitons in \eqref{eq:KdV-soliton-edge};  for each fixed $\epsilon>0$ however small, the soliton centers $a_n^\mathrm{KdV}(\epsilon)$ (in the coordinate $X_\srm^\mathrm{KdV}$) eventually become positive and tend to $+\infty$ as $n$ increases\footnote{This subtle point was unfortunately missed in \cite{ClaeysG10b} where the formula \eqref{eq:KdV-soliton-edge} is written with an infinite series.}.  This is of course quite different from the Benjamin-Ono case where on the scale $X$, the sequence of soliton centers (Airy zeros $a_n$) tends monotonically to $-\infty$ with decreasing spacing, and the asymptotic formula \eqref{eq:u-soliton-edge-sum} can indeed be written with an infinite series. 

\subsection{Asymptotic behavior near the harmonic edge}
\subsubsection{Formulation of main result}
When $(t_\cc,x_\cc)$ is a generic gradient catastrophe point and $x=x_\hrm(t)$ for $t-t_\cc$ positive and sufficiently small, there are again two distinct real solutions $y$ of the characteristic equation $x=2u_0(y)t+y$ with one solution denoted $y=y_\hrm(t)$ being a double root and the other denoted $y=y_2(t)$ being a simple root.  In this case the double root lies to the right of the simple root:  $y_\hrm(t)>y_2(t)$.

\begin{theorem}[Benjamin-Ono universality at the harmonic edge]
\label{thm:HarmEdge}
Let $u_0$ be a rational initial condition of the form~\eqref{eq:rationalIC}, and assume that $(t_\cc,x_\cc)$ is a generic gradient catastrophe point and moreover that all complex roots $y$ of the characteristic equation $x_\cc=2u_0(y)t_\cc+y$ are simple. 
Set $u^\Brm_\hrm(t) := u_0(y_\hrm(t))$, $u^\Brm_2(t) := u_0(y_2(t))$ and $k_\hrm(t) := 2 t [\frac{1}{2}u_0''(y_\hrm(t))]^{1/3} > 0$. 
Then, there exists $\delta>0$ sufficiently small that for $t\in K:=(t_\cc,t_\cc+\delta)$,
\begin{equation}
    u(t, x; \epsilon) = u^\Brm_2(t) + \epsilon^{1/6} \sqrt{\frac{2 \big( u_2^\Brm(t) - u_\hrm^\Brm(t) \big)}{\pi k_\hrm(t)}} \frac{\cos \left( \epsilon^{-1}\Theta(t, x)+\varphi(t,X_\hrm) \right)}{|\Ai(\ee^{\ii \pi/3} X_\hrm)|} + O(\epsilon^{1/3}), \quad \epsilon \to 0,
    \label{eq:BOHarmonicApprox}
\end{equation}
where the leading term in the phase is defined by
\begin{equation}
    \Theta(t,x):=\left(u_2^\Brm(t) - u_\hrm^\Brm(t)\right)(x-x_\hrm(t))+\left( u^\Brm_2(t)^2 -  u^\Brm_\hrm(t)^2 \right)t  +  \int_{y_\hrm(t)}^{y_2(t)} u_0(y)\,  \dd y, 
    \label{eq:BOHarmonicPhase-0}
\end{equation}
and a sub-leading term by
\begin{equation}
    \varphi(t,X_\hrm):=\arg \left( \Ai(\ee^{\ii \pi/3} X_\hrm) \right) - \frac{5\pi}{12} + \Lambda(t),
    \label{eq:BOHarmonicPhase-1}
\end{equation}
in which
\begin{equation}
\Lambda(t):=\frac{1}{2\pi} \int_0^{+\infty} \ln \left( \frac{g_\hrm(y_{\hrm}(t) + s, t)}{g_\hrm(y_{\hrm}(t) - s, t)} \right) \frac{\dd s}{s}-\frac{1}{2\pi} \int_0^{+\infty} \ln \left( \frac{g_\hrm(y_{2}(t) + s, t)}{g_\hrm(y_{2}(t) - s, t)} \right) \frac{\dd s}{s},
\label{eq:BOHarmonicPhaseh2}
\end{equation}
where $g_\hrm(y, t)$ is the positive rational function defined for $y \in \R$ and $t \in K$ by
\begin{equation}
    g_\hrm(y, t) := \frac{y + 2t u_0(y) - x_\hrm(t)}{\big( y + 2t u^\Brm_2(t) - x_\hrm(t) \big) \big( y + 2t u^\Brm_\hrm(t) -x_\hrm(t) \big)^2}.
\label{eq:BOHarmonicgdef}
\end{equation}
The approximation formula \eqref{eq:BOHarmonicApprox} holds uniformly for $X_\hrm$ bounded in $\mathbb{R}$ where
\begin{equation}
    X_\hrm := \frac{x-x_\hrm(t)}{k_\hrm(t) \epsilon^{2/3}}.
    \label{eq:Xh-coordinate}
\end{equation}
and $t$ in compact subsets of $K$.
\label{thm:harmonic-edge}
\end{theorem}

\begin{remark}
    A similar generalization as described in Remark~\ref{rem:OtherSolitonEdges} applies to Theorem~\ref{thm:HarmEdge} as well. Indeed, the definition of the harmonic edge can be generalized to the locus of points $(t, x)$ where the characteristic equation $x=2u_0(y)t+y$ has exactly one real double root to the right of one real simple root. This locus can be decomposed into a finite union of curves $\{ (t, x^i_\hrm(t)) \, | \, t \in I_\hrm^i \}$ where $I_\hrm^i$ is an open (possibly semi-infinite) interval of $\R_+$ and $x^i_\hrm$ is a differentiable function on $I_\hrm^i$. One can then check that the asymptotic formula in Theorem \ref{thm:HarmEdge} holds uniformly for $t \in K$, a compact subset of $I^i_\hrm$ for a particular harmonic edge $x = x^i_\hrm(t)$ on which all of the non-real roots $y$ of $x^i_\hrm(t)=2u_0(y)t+y$ remain simple.
\end{remark}

\begin{figure}[h]
\begin{center}
    \includegraphics[width=0.49\linewidth]{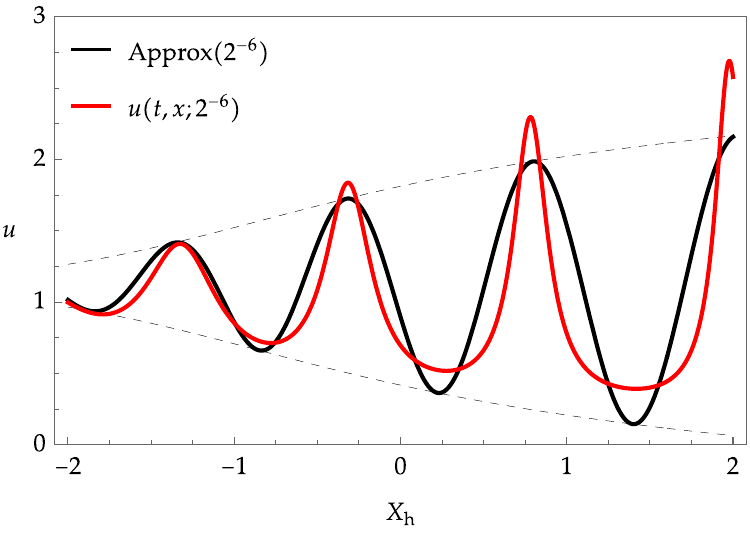}\hfill%
    \includegraphics[width=0.49\linewidth]{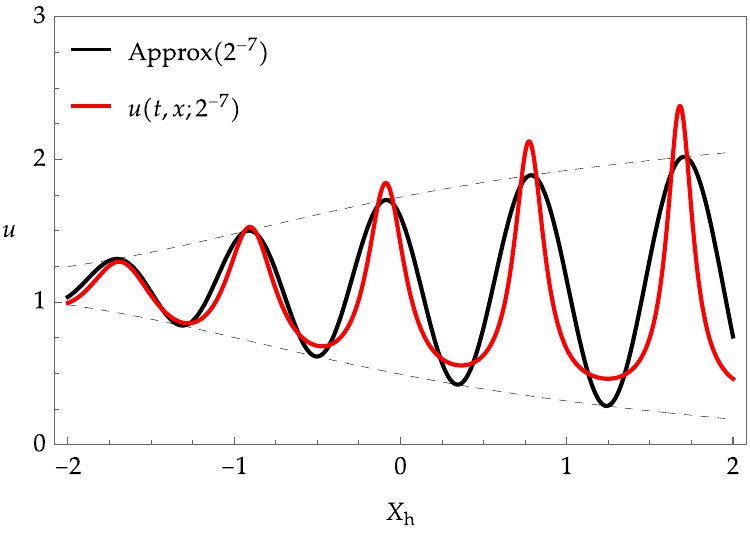}\\
    \includegraphics[width=0.49\linewidth]{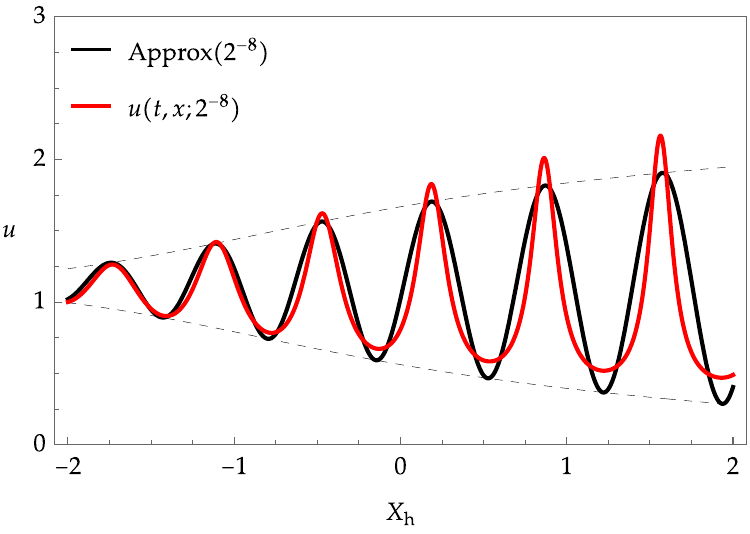}\hfill%
    \includegraphics[width=0.49\linewidth]{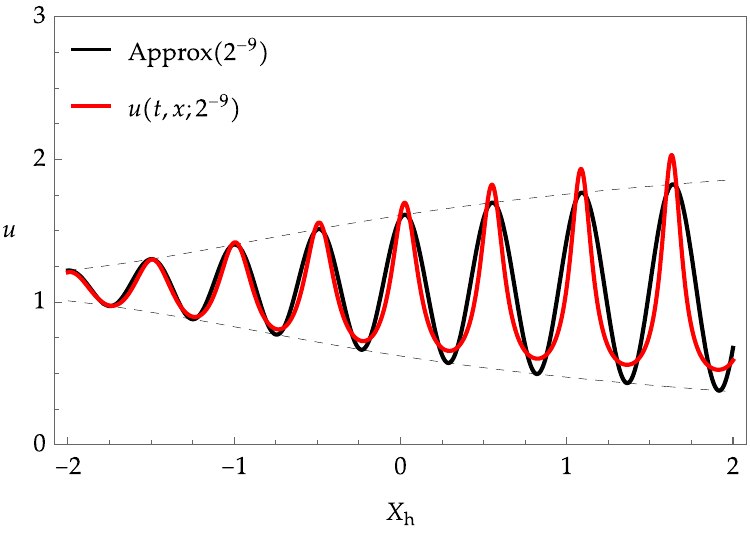}
\end{center}
\caption{Comparison between the exact solution $u(t,x;\epsilon)$ (black lines) with initial data $u_0(x)=2/(1+x^2)$ for $t=2>t_\cc$ and the explicit terms on the right-hand side of \eqref{eq:BOHarmonicApprox} (red lines) and the envelope (replacing the cosine with $\pm 1$, dashed lines) for $\epsilon=2^{-6}$ (upper left), $\epsilon=2^{-7}$ (upper right), $\epsilon=2^{-8}$ (lower left), and $\epsilon=2^{-9}$ (lower right).}
\label{fig:harmonic-edge}
\end{figure}

The accuracy of Theorem~\ref{thm:harmonic-edge} is illustrated in Figure~\ref{fig:harmonic-edge}.
One can judge from these plots that when $\epsilon$ decreases by a factor of $8$, the maximum error on the fixed interval of $X_\hrm$ shown in the plot decreases approximately by a factor of $2$, consistent with the $O(\epsilon^{1/3})$ error estimate in \eqref{eq:BOHarmonicApprox}.

Observe that 
\begin{equation}
    \kappa(t):=\frac{\partial \Theta}{\partial x}(t,x_\hrm(t))=u_2^\Brm(t)-u_\hrm^\Brm(t)>0.
\label{eq:kappa-define}
\end{equation}
Similarly, by direct calculation using also $u_2^\Brm(t):=u_0(y_2(t))$ and $u_\hrm^\Brm(t):=u_0(y_\hrm(t))$,
\begin{multline}
    \omega(t):=-\frac{\partial\Theta}{\partial t}(t,x_\hrm(t))=\left(u_2^\Brm(t)-u_\hrm^\Brm(t)\right)x_\hrm'(t) -u_2^\Brm(t)^2+u_\hrm^\Brm(t)^2 \\{}-2t\left(u_2^\Brm(t)u_2^{\Brm\prime}(t)-u_\hrm^\Brm(t)u_\hrm^{\Brm\prime}(t)\right) -u_2^\Brm(t)y_2'(t)+u_\hrm^\Brm(t)y_\hrm'(t).
    \label{eq:omega-define}
    \end{multline}
But, using also $y_2(t)=x_\hrm(t)-2tu_2^\Brm(t)$ and $y_\hrm(t)=x_\hrm(t)-2tu_\hrm^\Brm(t)$ we can rewrite this as
\begin{equation}\label{eq:omega-simple}
    \omega(t)=u_2^\Brm(t)^2-u_\hrm^\Brm(t)^2.
\end{equation}
Now fix $t=t_0\in K=(t_\cc,t_\cc+\delta)$.  The Benjamin-Ono equation admits constant solutions and its linearization about the constant $u_2^\Brm(t_0)$ (i.e., replace $u$ in \eqref{eq:BO} with $u_2^\Brm(t_0)+\Delta u$ for $\Delta u=\Delta u(t,x;\epsilon)$ and neglect the quadratic term in $\Delta u$) reads
\begin{equation}
    \partial_t\Delta u+ 2u_2^\Brm(t_0)\partial_x\Delta u= \epsilon\partial_x|D_x|\Delta u.
    \label{eq:linearizedBO}
\end{equation}
To obtain the dispersion relation we substitute a plane wave $\Delta u=\ee^{\ii (\kappa x-\omega t)/\epsilon}$ and get (for $\epsilon>0$)
\begin{equation}
    -\omega + 2u_2^\Brm(t_0)\kappa = \kappa|\kappa|.
    \label{eq:BO-linearized-dispersion}
\end{equation}
Therefore, we easily see that for given $t_0\in K$, the derived quantities $\kappa(t_0)$ and $\omega(t_0)$ defined in \eqref{eq:kappa-define} and \eqref{eq:omega-simple} respectively satisfy \eqref{eq:BO-linearized-dispersion}.  This means that the explicit terms in \eqref{eq:BOHarmonicApprox} have the form of a modulated wavepacket solution of the linearized Benjamin-Ono equation \eqref{eq:linearizedBO}, where the spatial modulation enters via dependence on $X_\hrm$.  Moreover, according to the linearized theory, the modulation envelope of the packet should move by translation at the group velocity $c_\mathrm{g}:=\dd\omega/\dd\kappa$.  We calculate from \eqref{eq:BO-linearized-dispersion} that $c_\mathrm{g}=2(u_2^\Brm(t_0)-\kappa)$ at wavenumber $\kappa>0$.  Hence
\begin{equation}
    \kappa=\kappa(t_0)\implies c_\mathrm{g}=2u_\hrm^\Brm(t_0).
\end{equation}
This matches exactly with the condition that a coordinate $X_\hrm$ in the range of the map $(t,x)\mapsto X_\hrm$ is held constant.  Indeed, since $y_\hrm(t)=y(t,x_\hrm(t))$ where $y(t,x)$ satisfies both the characteristic equation $x=y+2tu_0(y)$ and its $y$-derivative $0=1+2tu_0'(y)$, differentiation of the former evaluated for $x=x_\hrm(t)$ and $y=y_\hrm(t)$ gives $x_\hrm'(t)=(1+2tu_0'(y_\hrm(t)))y_\hrm'(t) + 2u_0(y_\hrm(t))=2u_\hrm^\Brm(t)$.  So $X_\hrm$ is stationary at $t=t_0$ when $\dd x/\dd t(t_0) = x_\hrm'(t_0)=2u_\hrm^\Brm(t_0)=c_\mathrm{g}$.

\subsubsection{Comparison with KdV}
For KdV, as mentioned above, the harmonic edge curve $x=x_\hrm^\mathrm{KdV}(t)$ is a different curve than $x=x_\hrm(t)$ in the Benjamin-Ono case, and it is specified along with two auxiliary functions $U_\hrm(t)$ and $V_\hrm(t)$ by the following system of equations:
\begin{equation}
    \begin{split}
        2tU_\hrm(t)+f_\mathrm{L}(U_\hrm(t))&=x_\hrm^\mathrm{KdV}(t)\\
        2t +\theta_\hrm(V_\hrm(t),U_\hrm(t))&=0\\
        \partial_V\theta_\hrm(V_\hrm(t),U_\hrm(t))&=0,
    \end{split}
\end{equation}
where for $U>V$,
\begin{equation}
    \theta_\hrm(V,U):=\frac{1}{2\sqrt{U-V}}\int_V^U\frac{f_\mathrm{L}'(\xi)\,\dd\xi}{\xi-V}.
\end{equation}
See \cite[Eqns.\@ (1.11)--(1.14)]{ClaeysG10a}.  The function $U_\hrm(t)$ is the continuation from the left of the solution of the inviscid Burgers equation, evaluated on the curve $x=x_\hrm^\mathrm{KdV}(t)$, so it is a KdV analogue of $u_2^\Brm(t)$ in the Benjamin-Ono case.

In the paper \cite{ClaeysG10a}, Claeys and Grava studied the initial-value problem for KdV under the same assumptions mentioned in Section~\ref{sec:soliton-comparison} and analyzed the solution $u(t,x;\epsilon)$ near the harmonic edge $x=x_\hrm^\mathrm{KdV}(t)$.  Their asymptotic formula involves a particular global solution of the homogeneous Painlev\'e-II equation
\begin{equation}
    q''(s) =sq(s)+2q(s)^3
\end{equation}
known as the \emph{Hastings-McLeod solution}, denoted $q_\mathrm{HM}(s)$.  It is smooth and real-valued for $s\in\mathbb{R}$, and is uniquely determined by the asymptotic condition that 
\begin{equation}
    q_\mathrm{HM}(s)=\Ai(s)(1+o(1)),\quad s\to+\infty.
\end{equation}
It is known to have square-root asymptotics in the opposite direction:
\begin{equation}
    q_\mathrm{HM}(s)=\sqrt{-\frac{s}{2}}(1+o(1)),\quad s\to-\infty,
\end{equation}
which is an example of a connection formula for Painlev\'e transcendents.

To describe the result for KdV, set
\begin{equation}
    p(t):=-27\sqrt{\frac{U_\hrm(t)-V_\hrm(t)}{3}}\partial^2_V\theta_\hrm(V(t),U(t))>0,
\end{equation}
and taking the positive cube root, define a rescaled and recentered spatial variable by
\begin{equation}
    X^\mathrm{KdV}_\hrm:=\frac{1}{p(t)^{1/3}}\sqrt{\frac{3}{U_\hrm(t)-V_\hrm(t)}}\frac{x-x^\mathrm{KdV}_\hrm(t)}{\epsilon^{2/3}}.
\end{equation}
Finally, define a phase by
\begin{multline}
    \Theta^\mathrm{KdV}(t,x):=2\sqrt{\frac{U_\hrm(t)-V_\hrm(t)}{3}}\left(x-x_\hrm^\mathrm{KdV}(t)\right)\\+2\int_{V_\hrm(t)}^{U_\hrm(t)}\left(f_\mathrm{L}'(\xi)+2t\right)\sqrt{\xi-V_\hrm(t)}\,\dd\xi.
\end{multline}
Then, the main result of \cite{ClaeysG10a} is that the solution of the KdV initial-value problem with initial data $u_0$ satisfies
\begin{equation}
    u(t,x;\epsilon)=U_\hrm(t)-\frac{12\epsilon^{1/3}}{p(t)^{1/3}}q_\mathrm{HM}(-X_\hrm^\mathrm{KdV})\cos\left(\epsilon^{-1}\Theta^\mathrm{KdV}(t,x)\right) + O(\epsilon^{2/3}),\quad \epsilon\to 0,
    \label{eq:KdV-harmonic-edge}
\end{equation}
uniformly for $X_\hrm^\mathrm{KdV}$ bounded.

Like Theorem~\ref{thm:harmonic-edge}, this result describes the wave field $u(t,x;\epsilon)$ as a small and highly-oscillatory fluctuation of the background inviscid Burgers solution with a well-defined and slowly varying (relatively) amplitude modulation.  One difference is that for the Benjamin-Ono equation, the size of the fluctuation is much larger (proportional to $\epsilon^{1/6}$) than for the KdV equation (proportional to $\epsilon^{1/3}$).  Another important distinction is that while for KdV the amplitude modulation is described by the Painlev\'e transcendent $q_\mathrm{HM}(-\diamond)$, for the Benjamin-Ono equation one finds a more classical special function, namely $1/|\Ai(\ee^{\ii\pi/3}\diamond)|$. Both of these are real-analytic functions, with different asymptotic behavior.  As noted earlier,
\begin{equation}
    q_\mathrm{HM}(-X)=\begin{cases}\displaystyle\frac{\ee^{-\frac{2}{3}(-X)^{3/2}}}{2\sqrt{\pi}(-X)^{1/4}}(1+o(1)),&X\to-\infty,\\
    \displaystyle\sqrt{\frac{X}{2}}(1+o(1)),& X\to+\infty,
    \end{cases}
\end{equation}
while
\begin{equation}
    \frac{1}{|\Ai(\ee^{\ii\pi/3}X)|}=\begin{cases}\displaystyle 2\sqrt{\pi}(-X)^{1/4}\ee^{-\frac{2}{3}(-X)^{3/2}}(1+o(1)),& X\to-\infty,\\
    2\sqrt{\pi}X^{1/4}(1+o(1)),& X\to+\infty.
    \end{cases}
\end{equation}
Both functions therefore decay exponentially (with reciprocal power-law corrections) as $X\to-\infty$, which is the direction moving from the harmonic edge $x=x_\hrm(t)$ away from the dispersive shock wave.  In the opposite direction, the Benjamin-Ono amplitude grows more slowly into the dispersive shock wave (proportional to $X^{1/4}$) than does the KdV amplitude (proportional to $X^{1/2}$).  Both of these growth rates are consistent with the rates of decay of the amplitude of the oscillations in the dispersive shock wave going toward the harmonic edge.  

Another interesting difference between the asymptotic formul\ae\ \eqref{eq:BOHarmonicApprox} and \eqref{eq:KdV-harmonic-edge} is that in the KdV case there is evidently no need for a lower-order phase correction, while for Benjamin-Ono the phase correction $\varphi(t,X_\hrm)$ is not negligible. It describes the slow variation of a complex amplitude $1/\Ai(\ee^{-\ii\pi/3}X_\hrm)$.

\subsection{Asymptotic behavior near the critical point}

\subsubsection{Formulation of main result}

\begin{figure}[h]
\begin{center}
    \includegraphics[width=0.65\linewidth]{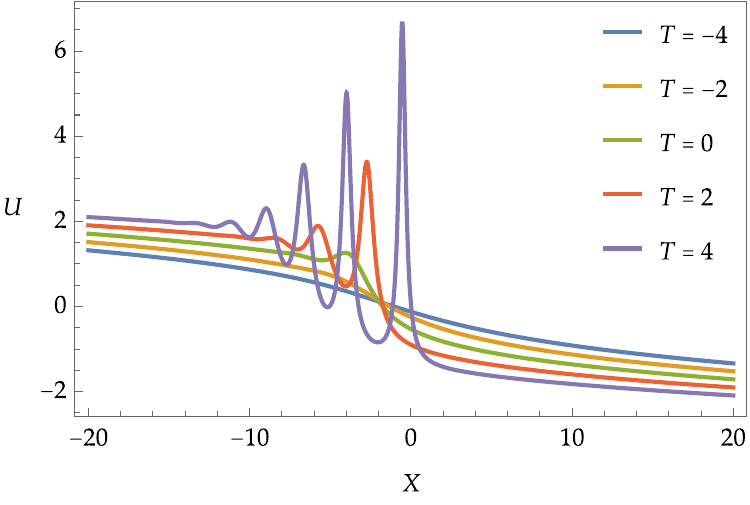}
\caption{The function $U(T,X)$ as a function of $X$ for various values of $T$.}
\end{center}
\label{fig:U-Plot}
\end{figure}

Let a function be defined on $\mathbb{R}^2$ by
\begin{equation}
    U(T,X):=-2\mathrm{Im}\left(\frac{\partial}{\partial X}\log(\tau(T,X))\right),
\label{eq:U-profile}
\end{equation}
where $\tau$ is an entire function of $(T,X)\in\mathbb{C}^2$ defined by the Pearcey integral
\begin{equation}
    \tau(T,X):=\int_{\infty\ee^{-7\ii\pi/8}}^{\infty\ee^{-3\ii\pi/8}}\ee^{\ii (\zeta^4-T\zeta^2-X\zeta)}\,\dd\zeta.
\label{eq:tau-define}
\end{equation}
In particular, it follows from classical steepest descent that
\begin{equation}
\begin{split}
    U(T,X)&=-\mathrm{sgn}(X)\left(\frac{1}{4}|X|\right)^{1/3}+O(|X|^{-1/3})\\
    U_X(T,X)&=-\frac{1}{12}\left(\frac{1}{4}|X|\right)^{-2/3}+O(|X|^{-4/3})
\end{split}
\label{eq:UandUX-asymptotics}
\end{equation}
as $X\to\pm\infty$ for each fixed $T\in\mathbb{R}$.  See Figure~\ref{fig:U-Plot}.
At a generic gradient catastrophe point $(t_\cc,x_\cc)$ the unique real solution $y=y_\cc$ of the characteristic equation $x_\cc=2u_0(y)t_\cc+y$ is a triple root, so one also has the identities $t_\mathrm{c}=-1/(2u_0'(y_\cc))$ (from the first derivative, and note that $t_\cc>0$ implies $u_0'(y_\cc)<0$) and $u_0''(y_\cc)=0$ (from the second derivative), while the condition that $y_\cc$ is not a root of order greater than three is equivalent to the inequality $u_0'''(y_\cc)>0$.
\begin{theorem}[Benjamin-Ono universality at the gradient catastrophe point]
Let $u_0$ be a rational initial condition of the form~\eqref{eq:rationalIC}, and assume that $(t_\cc,x_\cc)$ is a generic gradient catastrophe point and moreover that all complex roots $y$ of the characteristic equation $x_\cc=2u_0(y)t_\cc+y$ are simple. 
Set $u_\mathrm{c}=u_0(y_\mathrm{c})=u^\mathrm{B}(t_\mathrm{c},x_\mathrm{c})$ and $k:=(\frac{2}{3})^{1/4}t_\mathrm{c}u_0'''(y_\mathrm{c})^{1/4}>0$.  Then, 
\begin{equation}
    u(t,x;\epsilon)= u_\mathrm{c} +\frac{\epsilon^{1/4}}{k}U(T,X) + O(\epsilon^{1/2}),\quad\epsilon\to 0 
\label{eq:u-critical-statement}
\end{equation}
holds uniformly for $(T,X)$ in compact subsets of $\mathbb{R}^2\setminus Z$, where
\begin{equation}
    T:=\frac{t-t_\mathrm{c}}{k^2\epsilon^{1/2}},\quad X:=\frac{x-x_\mathrm{c}-2u_\mathrm{c}(t-t_\mathrm{c})}{k\epsilon^{3/4}}
    \label{eq:catastrophe-coordinates}
\end{equation}
and $Z$ denotes the zero divisor of $\tau(T,X)$.
\label{thm:BO-catastrophe}
\end{theorem}

\begin{remark}
    Once again, for a given initial condition $u_0$ the definition of a generic gradient catastrophe point $(t_\cc, x_\cc)$ can be extended to allow for any point $(t_\cc,x_\cc)\in\R_+\times\R$ having an open neighborhood $N\subset\R^2$ such that $x=2u_0(y)t+y$ has only one real simple root on $N \cap \{ (t, x) \, | \, 0 \leq t < t_\cc \}$ and $x_\cc=2u_0(y)t_\cc+y$ has exactly one real triple root. 
    The proof of Theorem \ref{thm:BO-catastrophe} then shows that it applies at any other generic gradient catastrophe points allowed under this extended definition. As an example, the unlabeled cusp point near $x\approx 13$ on the caustic locus $C$ shown in Figure~\ref{fig:edge} would satisfy this extended definition.
\end{remark}

\begin{figure}[h]
\begin{center}
    \includegraphics[width=0.49\linewidth]{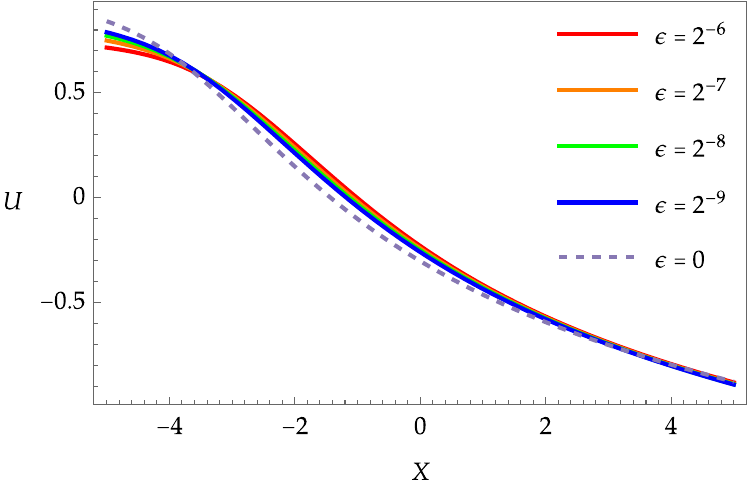}\hfill%
    \includegraphics[width=0.49\linewidth]{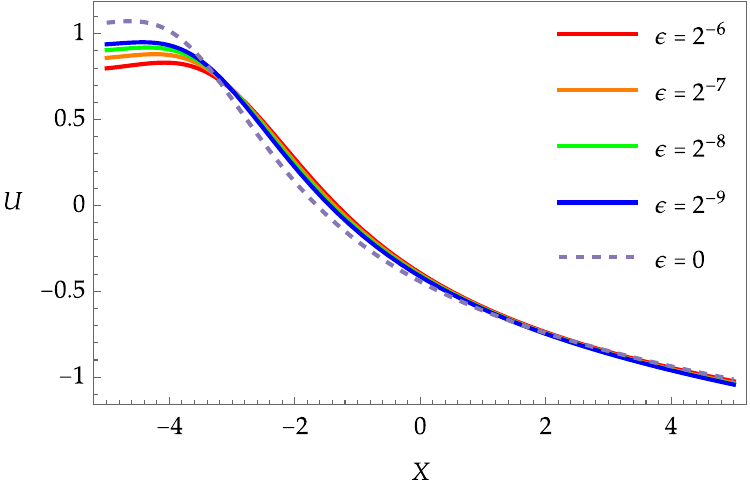}\\
    \includegraphics[width=0.49\linewidth]{CritConvergenceTm0p5.pdf}\hfill%
    \includegraphics[width=0.49\linewidth]{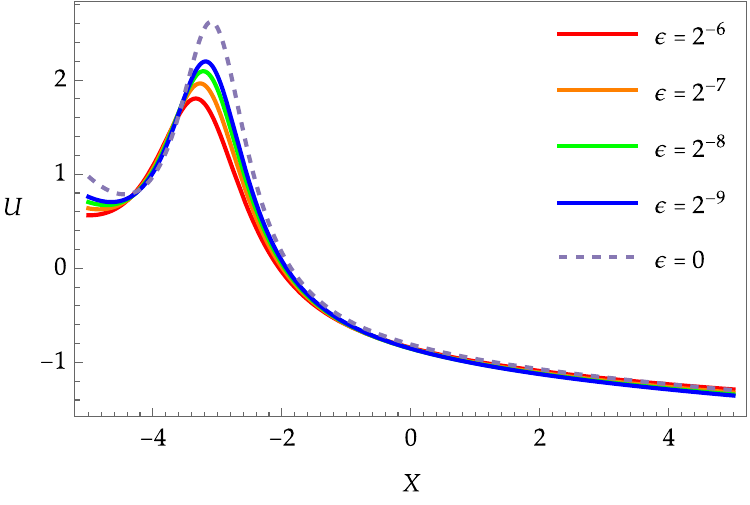}%
\end{center}
\caption{For $u_0(x)=2/(1+x^2)$, the function $X\mapsto k\epsilon^{-1/4}(u(t_\cc+k^2\epsilon^{1/2}T,x_\cc+2u_\cc k^2\epsilon^{1/2}T+k\epsilon^{3/4}X;\epsilon)-u_\cc)$ for $\epsilon=2^{-6},2^{-7},2^{-8},2^{-9}$ (solid curves) along with $U(T,X)$ (dashed curves) for different fixed $T$.  Top row:  left: $T=-1.5$, right:  $T=-0.5$.  Bottom row:  left: $T=0.5$, right: $T=1.5$. }
\label{fig:CritConvergence}
\end{figure}

The convergence asserted in Theorem~\ref{thm:BO-catastrophe} is illustrated in Figure~\ref{fig:CritConvergence}.  One can see from these plots that when $\epsilon$ decreases by a factor of $16$, the error measured by $k\epsilon^{-1/4}(u(t,x;\epsilon)-u_0)-U(T,X)$ decreases by approximately half, which is consistent with the error estimate claimed in \eqref{eq:u-critical-statement}.

Note that since for each $T\in\mathbb{R}$, $X\mapsto \tau(T,X)$ is an entire function with nonzero asymptotics for large $X$, it has only finitely many zeros for fixed $T$.  Numerical evidence for bounded values of $X$ supports the following conjecture, the proof of which we leave for future work.
\begin{conjecture}
    $Z\cap\mathbb{R}^2=\emptyset$, i.e., $\tau(T,X)\neq 0$ holds for all $(T,X)\in\mathbb{R}^2$.
\label{conj:tau-nonzero}
\end{conjecture}

\subsubsection{Comparison to KdV}
In a 2006 paper \cite{Dubrovin06}, Dubrovin studied a class of Hamiltonian, weakly dispersive, perturbations to the inviscid Burgers equation, and used canonical transformations and formal perturbation theory to study the way that solutions of such dispersive models behave in a double-scaling limit that the dispersion coefficient $\epsilon$ is small and that the independent variables are localized near a generic gradient catastrophe point $(t_\cc,x_\cc,u_\cc)$ for the unperturbed equation.  The generic nature of the catastrophe is captured by the condition that 
\begin{equation} 
k:=-27f_\mathrm{L}'''(u_\mathrm{c})=27\frac{u_0'''(y_\cc)}{u_0'(y_\cc)^4}>0, 
\end{equation}
where as above $f_\mathrm{L}(u)$ denotes the inverse function of the decreasing part of the initial data $u_0$.
Dubrovin's formal analysis led to the conjecture that the dispersive solution is modeled near the gradient catastrophe point by a particular solution $\mathcal{D}=\mathcal{D}(X;T)$ of the second equation in the Painlev\'e-I hierarchy
\begin{equation}
    X=T\mathcal{D}-\left[\frac{1}{6}\mathcal{D}^3 + \frac{1}{24}\left(\left(\frac{\dd\mathcal{D}}{\dd X}\right)^2+2\mathcal{D}\frac{\dd^2\mathcal{D}}{\dd X^2}\right)+\frac{1}{240}\frac{\dd^4\mathcal{D}}{\dd X^4}\right],
    \label{eq:PI-2}
\end{equation}
which is a fourth-order ordinary differential equation with a parameter $T\in\mathbb{R}$.  The particular solution of interest is singled out for each $T\in\mathbb{R}$ by the property that it behaves as $X\to\pm\infty$ as $\mathcal{D}(X;T)=\mp (6|X|)^{1/3}+O(|X|^{-1/3})$ and is real-valued and smooth for $X\in\mathbb{R}$. The existence and uniqueness of this solution was proven later by Claeys and Vanlessen \cite{ClaeysV07}.  The Dubrovin conjecture is remarkable in its scope:  it asserts universality of the limiting profile of the dispersive solution with respect to both initial data and the choice of dispersive Hamiltonian perturbing terms.  The KdV equation belongs to the class of dispersive perturbations of the inviscid Burgers equation to which Dubrovin's theory applies, and in 2009, Claeys and Grava proved \cite{ClaeysG09} that Dubrovin's conjecture holds for the initial-value problem for the KdV equation.
Defining rescaled and recentered variables by
\begin{equation}
    T:=\frac{6(t-t_\cc)}{(4k^3\epsilon^4)^{1/7}}\quad\text{and}\quad X:=\frac{x-x_\cc-2u_\cc(t-t_\cc)}{(8k\epsilon^6)^{1/7}},
\end{equation}
the main result of \cite{ClaeysG09} is that, under the assumptions on $u_0$ stated in Section~\ref{sec:soliton-comparison}, the solution $u(t,x;\epsilon)$ of the initial-value problem for the KdV equation in the form \eqref{eq:KdV} with initial data $u_0$, for which the solution of the inviscid Burgers equation exhibits a generic gradient catastrophe point, satisfies  
\begin{equation}
    u(t,x;\epsilon)=u_\cc + 3\left(\frac{2\epsilon^2}{k^2}\right)^{1/7}\mathcal{D}(X;T) + O(\epsilon^{4/7}),\quad \epsilon\to 0
\end{equation}
uniformly for coordinates $(T,X)$ in compact subsets of $\mathbb{R}^2$.

While it is a weakly-dispersive and Hamiltonian perturbation of the inviscid Burgers equation, the Benjamin-Ono equation \eqref{eq:BO} does not lie in Dubrovin's universality class, so it is not surprising that our result in Theorem~\ref{thm:BO-catastrophe} involves a different limiting function $U(T,X)$. It is perhaps remarkable that the function $U(T,X)$ is far more explicit than is $\mathcal{D}(X;T)$, a result that continues the theme of Theorems~\ref{thm:soliton-edge} and \ref{thm:harmonic-edge}.  

In \cite{MasoeroRA15}, a formal multiscale analysis parallel to that of Dubrovin \cite{Dubrovin06} was carried out for a wider variety of dispersive nonlocal Hamiltonian perturbations of the inviscid Burgers equation that includes nonlocal equations such as the Benjamin-Ono equation.  That analysis suggested that among the nonlocal equations considered, the Benjamin-Ono equation was the only one that did not ultimately fall into Dubrovin's universality class.  The scaling analysis presented in \cite{MasoeroRA15} for the Benjamin-Ono case correctly predicts the $(\Delta u,\Delta t,\Delta x-2u_\cc\Delta t)\sim (\epsilon^{1/4},\epsilon^{1/2},\epsilon^{3/4})$ scaling evident in Theorem~\ref{thm:BO-catastrophe}.  There was no concrete prediction given in \cite{MasoeroRA15} of the limiting wave profile $U(T,X)$ given in \eqref{eq:U-profile}--\eqref{eq:tau-define}; however it was predicted that the profile should satisfy an analogue of \eqref{eq:PI-2} which in our scaling reads 
\begin{equation}
    X-2UT +4U^3 - 6U|D_X|U-3|D_X|U^2-4\partial_X^2U=0.
    \label{eq:string-equation}
\end{equation}
As in \eqref{eq:PI-2}, $T$ is considered a parameter.  This is a nonlocal ordinary differential equation that was proposed in \cite{MasoeroRA15} as a new type of Painlev\'e equation.  This leads us to formulate the following conjecture.
\begin{conjecture}
Let $U(T,X)$ be defined explicitly by \eqref{eq:U-profile}--\eqref{eq:tau-define}, and assume that Conjecture~\ref{conj:tau-nonzero} holds.  Then 
$U(T,X)$ is a global solution of \eqref{eq:string-equation}.
\label{conj:string-equation}
\end{conjecture}

We have verified Conjecture~\ref{conj:string-equation} numerically.  To do this, we chose specific points $(T,X)$, and evaluated the terms not involving the nonlocal operator $|D_X|$ by numerical contour integration.  For the terms involving $|D_X|$, we use the singular integral representation \eqref{eq:AbsDx} which converges pointwise in $X$ for the terms in question even though the functions to which the Hilbert transform is applied are not in $L^2(\mathbb{R})$ according to the asymptotic behavior given in \eqref{eq:UandUX-asymptotics}.  The principal value integrals can be computed accurately by integrating on a large interval $[-L,L]$ containing the singularity $X$ using \texttt{NIntegrate} in \texttt{Mathematica} (version 14.0) with the option \texttt{Method->"PrincipalValue"}; we then fill in the integration over $(-\infty,-L]\cup[L,+\infty)$ by omitting the principal-value option and replacing the functions $\partial_XU$ and $2U\partial_XU$ in the integrand by their asymptotics \eqref{eq:UandUX-asymptotics} to avoid evaluating the contour integrals for arguments that are too extreme.  For $L=100$ and values of $(T,X)$ not too large, this establishes that the left-hand side of \eqref{eq:string-equation} vanishes to about 5 digits of accuracy.

While we do not have a proof of Conjecture~\ref{conj:string-equation}, we can easily derive a coupled system of ordinary differential equations involving $U(T,X)$ and an auxiliary function.  Using the differential equation $4\partial_X^3\tau+2T\partial_X\tau-\ii X\tau=0$ satisfied by the Pearcey integral \eqref{eq:tau-define}, upon introducing the complex-valued function $\Lambda(T,X):=2\ii\partial_X\log(\tau(T,X))$, it is straightforward to show that $X+\Lambda T-\frac{1}{2}\Lambda^3-3\ii \Lambda\partial_X\Lambda+2\partial_X^2\Lambda=0$. Noting that $U(T,X)=\mathrm{Re}(\Lambda(T,X))$ and letting $V(T,X):=\mathrm{Im}(\Lambda(T,X))$, this leads to the coupled system:
\begin{equation}
    \begin{split}
        X+UT-\frac{1}{2}U^3+\frac{3}{2}UV^2+3\partial_X(UV)+2\partial_X^2U&=0\\
        VT-\frac{3}{2}U^2V+\frac{1}{2}V^3-3U\partial_XU+3V\partial_XV+2\partial_X^2V&=0.
    \end{split}
\end{equation}
Evidently, it should be possible to eliminate $V(T,X)$ in favor of $U(T,X)$ and nonlocal quantities such as $|D_X|U$ and $|D_X|U^2$ using suitable analytic properties of $\Lambda(T,X)$ in a half-plane of $X$ to obtain \eqref{eq:string-equation}.  We leave this for future work.

\subsection*{Acknowledgements}  
E. Blackstone was partially supported by the National Science Fou\-ndation under grant DMS-1812625.  P. D. Miller and M. D. Mitchell were partially supported by the National Science Foundation under grant DMS-2204896.  The authors would like to thank Louise Gassot and Sheehan Olver for useful discussions.

\section{Setup}
\subsection{Solution formula for the Benjamin-Ono equation}
Our proofs are based on a recently obtained formula \cite{RationalBO} for the exact solution of the initial-value problem \eqref{eq:BO}--\eqref{eq:rationalIC} in terms of the ratio of two determinants of dimension $(N+1)\times (N+1)$, whose entries are explicit contour integrals.  Let $\log$ denote the principal branch complex logarithm such that $|\operatorname{Im}(\log(\cdot))|<\pi$, and define for $y\in\mathbb{R}$ 
\begin{equation}
h(y;t,x):=\frac{1}{4t}(y-x)^2+\sum_{n=1}^N\left[c_n\log(y-p_n)+c_n^*\log(y-p_n^*)\right].
\label{eq:h-def}
\end{equation}
We next analytically continue $h$ to a maximal domain of analyticity in the complex plane and call the resulting function $z\mapsto h(z;t,x)$.  We have to omit a branch cut emanating from each of the logarithmic singularities, and to formulate the solution we assume that these cuts are pairwise disjoint.  The branch cuts $\Gamma_1,\dots,\Gamma_N$ emanating from the points $p_1,\dots,p_N$ in the upper half-plane are all assumed to tend to infinity in the direction $\arg(z)=\frac{3}{4}\pi$.  The branch cuts $\bar{\Gamma}_1,\dots,\bar{\Gamma}_N$ emanating from the points $p_1^*,\dots,p_N^*$ in the lower half-plane are less important to specify, but we assume that they all lie below the real line and to the lower left of a diagonal line $\mathrm{Re}(z)+\mathrm{Im}(z)=K$.  Then we define a contour $W_0$ in the domain of analyticity of $h(z;t,x)$ that originates at $z=\infty\ee^{3\pi\ii/4}$ to the left of all branch cuts $\Gamma_1,\dots,\Gamma_N$ and terminates at $z=\infty\ee^{-\ii\pi/4}$ to the right of all branch cuts $\bar{\Gamma}_1,\dots,\bar{\Gamma}_N$.  We also define for each $n=1,\dots,N$ a contour $W_n$ in the domain of analyticity of $h(z;t,x)$ that begins and ends at $z=\infty\ee^{3\pi\ii/4}$ and that encircles only the branch cut $\Gamma_n$.  

For such contours, the following contour integrals are convergent for $\epsilon>0$ and $t>0$:
\begin{equation}
\begin{split}
A_{j1}(t,x;\epsilon)&:=\int_{W_{j-1}}u_0(z)\ee^{-\ii h(z;t,x)/\epsilon}\,\dd z, \quad j=1,\dots,N+1,
	\\
B_{j1}(t,x;\epsilon)&:=\int_{W_{j-1}}\ee^{-\ii h(z;t,x)/\epsilon}\,\dd z, \quad j=1,\dots,N+1
\end{split}
\label{eq:first-columns}
\end{equation}
and
\begin{equation}
A_{jk}(t,x;\epsilon)=B_{jk}(t,x;\epsilon):=\int_{W_{j-1}}\frac{\ee^{-\ii h(z;t,x)/\epsilon}\,\dd z}{z-p_{k-1}},\quad j=1,\dots,N+1,\quad k=2,\dots,N+1.
\label{eq:other-columns}
\end{equation}
Then, the contour integrals define the elements of two square matrices of dimension $(N+1)\times (N+1)$ denoted by $\mathbf{A}(t,x;\epsilon)$ and $\mathbf{B}(t,x;\epsilon)$.
We set
\begin{equation}
    N(t,x;\epsilon):=\det(\mathbf{A}(t,x;\epsilon))\quad\text{and}\quad
    D(t,x;\epsilon):=\det(\mathbf{B}(t,x;\epsilon)).
    \label{eq:num-and-demon-defs}
\end{equation}
These determinants are both entire functions of $x$ for $\epsilon>0$ and $t>0$.
The following theorem, proven in \cite[Theorem 1.3]{RationalBO}, then provides the exact solution of the initial-value problem for \eqref{eq:BO} with rational initial data \eqref{eq:rationalIC}.
\begin{theorem}
Suppose that $\epsilon>0$ is such that for each index $n=1,\dots,N$, $-\ii c_n\not\in\epsilon\mathbb{N}$.  Then for all $t>0$, $D(t,\diamond;\epsilon)$ has no real zeros, and the solution of the initial-value problem for the Benjamin-Ono equation \eqref{eq:BO} with rational initial data \eqref{eq:rationalIC} is given by
\begin{equation}
    u(t,x;\epsilon)=2\mathrm{Re}\left(\frac{N(t,x;\epsilon)}{D(t,x;\epsilon)}\right).
    \label{eq:IVP-solution}
\end{equation}
\label{thm:IVP-solution}
\end{theorem}
In \cite{RationalBO}, this theorem is proved in the more general setting allowing for the possibility that also $-\ii c_n\in\epsilon\mathbb{N}$ for one or more of the indices $n=1,\dots,N$, which then requires a different choice of contours to define the determinants $N(t,x;\epsilon)$ and $D(t;x,\epsilon)$.  One can see that by elementary row operations that do not affect $u(t,x;\epsilon)$, the contours $W_0,\dots,W_N$ can be replaced by others, and in fact the proof of Theorem~\ref{thm:IVP-solution} in \cite{RationalBO} uses a different system of equivalent contours denoted $C_0,\dots,C_N$ in that paper.

\subsection{The limit $\epsilon\to 0$}
It is then clear that to study the limit $\epsilon\to 0$ in the explicit formula \eqref{eq:IVP-solution} requires the analysis of the contour integrals \eqref{eq:first-columns}--\eqref{eq:other-columns} by the classical method of steepest descent.  The exponent function for all of these integrals is exactly the same:  $z\mapsto -\ii h(z;t,x)$.  As is well-known, the steepest descent analysis of integrals with exponent $-\ii h(z;t,x)$ hinges on the critical points of the exponent and the way they depend on auxiliary parameters such as the independent variables $(t,x)$.

Since $h'(z;t,x)=0$ is equivalent to the equation $x=z+2tu_0(z)$ according to \eqref{eq:rationalIC}, the real critical points $z=y\in\mathbb{R}$ are precisely the intercepts $y$ of the characteristic lines in the $(t,x)$-plane for the solution of the inviscid Burgers equation $\partial_tu^\Brm+2u^\Brm\partial_xu^\Brm=0$.    In general for $u_0$ rational of the form \eqref{eq:rationalIC} there are also complex-conjugate pairs of critical points, and while all of the real and complex critical points of $h$ play a role when $\epsilon>0$ is small, those that are real or have small imaginary parts are the most important.  

In the paper \cite{SmallDispersionBO-1}, the steepest descent method is applied to deduce the asymptotic behavior of $u(t,x;\epsilon)$ for rational initial data \eqref{eq:rationalIC} from \eqref{eq:IVP-solution}, assuming that all critical points of $h(z;t,x)$ are simple.  The most technical part of that analysis is the algorithmic specification, given $(t,x)\in\mathbb{R}^2$ near $(t_\cc, x_\cc)$, of suitable branch cuts $\Gamma_1,\dots,\Gamma_N$ and contours $W_0,\dots,W_N$ for which the integrals over each contour $W_n$ are dominated by precisely one critical point of $h$ in the upper half-plane, or by precisely one pair of real critical points, and the dominating critical point or real pair is uniquely associated to the contour $W_n$.  This ensures that the leading term of the determinants $N$ and $D$ is in each case the corresponding determinant of leading terms obtained by approximating the integrals in the matrix elements using contributions from the dominating critical point(s).  Indeed, while the ratio of determinants is insensitive to the choice of contours, the asymptotic analysis can require the extraction of terms beyond all orders if an improper choice of contours is made.

It should be possible to modify the analysis in \cite{SmallDispersionBO-1} to allow for the coalescence of complex critical points of $h(z;t,x)$, however it breaks down in an essential way if there are any real critical points that are not simple.  Since the real critical points are the intercepts $y$ of characteristic lines through $(t,x)\in\mathbb{R}^2$ for the solution of the inviscid Burgers equation $\partial_t u^\Brm +2u^\Brm\partial_x u^\Brm=0$ with initial data $u_0$, the set of points $(t,x)\in\mathbb{R}^2$ for which there is at least one real critical point of $h(z;t,x)$ of higher multiplicity coincides precisely with the caustic locus $C$ described in the introduction.  To be precise, the caustic locus $C$ consists of those points $(t,x)\in\mathbb{R}^2$ for which both the characteristic equation and its derivative with respect to $y$ simultaneously hold:
\begin{equation}
    (t,x)\in C\iff \exists y\in\mathbb{R}: \begin{cases} x=y+2tu_0(y),\\0=1+2tu_0'(y).\end{cases}
    \label{eq:Caustics-Definition}
\end{equation}
Clearly, there are no points of $C$ with $t=0$.
At a gradient catastrophe point $(t_\cc,x_\cc)\in C$, the second derivative with respect to $y$ also vanishes, which determines corresponding value(s) of $y$ denoted $y_\cc\in\mathbb{R}$ by $u_0''(y_\cc)=0$.  The coordinates $(t_\cc,x_\cc)$ are then determined explicitly in terms of $y_\cc$ from the two equations in \eqref{eq:Caustics-Definition} at $y=y_\cc$:  $t_\cc=-1/(2u_0'(y_\cc))$ and $x_\cc=y_\cc -u_0(y_\cc)/u_0'(y_\cc)$.  Assuming that the gradient catastrophe point is generic in the sense that $u_0'''(y_\mathrm{c})>0$, we expand $u_0$ in a Taylor series about $y=y_\cc$:
\begin{equation*}
    u_0(y)=u_\cc -(y-y_\cc)/(2t_\cc) + \frac{1}{6}u_0'''(y_\cc)(y-y_\cc)^3 + O((y-y_\cc)^4).
\end{equation*}
Here we defined $u_\cc:=u_0(y_\cc)$ and eliminated $u_0'(y_\cc)$ using the second equation from \eqref{eq:Caustics-Definition} at $y=y_\cc$ and $t=t_\cc$.  Substituting the expansion into the second equation from \eqref{eq:Caustics-Definition} for general $y\approx y_\cc$ and $t\approx t_\cc$ yields
\begin{equation}
    y-y_\cc = \pm\frac{(t-t_\cc)^{1/2}}{t_\cc u_0'''(y_\cc)^{1/2}} + o((t-t_\cc)^{1/2}),\quad t\downarrow t_\cc.
    \label{eq:y-near-yc}
\end{equation}
Using the Taylor series of $u_0(y)$ in the first equation from \eqref{eq:Caustics-Definition} then gives
\begin{equation}
    x=y+2tu_\cc -\frac{t}{t_\cc}(y-y_\cc)+\frac{1}{3}tu_0'''(y_\cc)(y-y_\cc)^3 + O((y-y_\cc)^4),\quad (t,x)\in C.
    \label{eq:first-equation-expanded}
\end{equation}
Eliminating $y$ from \eqref{eq:first-equation-expanded} using \eqref{eq:y-near-yc} gives the relationship between $x$ and $t$ near the catastrophe point:
\begin{equation}
    x-x_\cc -2u_\cc\cdot (t-t_\cc)=\mp\frac{2(t-t_\cc)^{3/2}}{3t_\cc^2u_0'''(y_\cc)^{1/2}} + o((t-t_\cc)^{3/2}),\quad (t,x)\in C,\quad t\downarrow t_\cc.
    \label{eq:Caustic-approximations}
\end{equation}
Alternatively, one can use \eqref{eq:Caustic-approximations} to eliminate $x$ from \eqref{eq:first-equation-expanded} to obtain an equation for all of the critical points $y\approx y_\cc$.  The most important terms give a cubic equation for $y$:
\begin{equation}
    \frac{1}{3}t_\cc u_0'''(y_\cc)(y-y_\cc)^3-\frac{1}{t_\cc}(t-t_\cc)(y-y_\cc) \pm\frac{2(t-t_\cc)^{3/2}}{3t_\cc^2u_0'''(y_\cc)^{1/2}}\approx 0.
\end{equation}
It is easy to see that this cubic has a repeated factor (this is the very meaning of the fact that \eqref{eq:Caustic-approximations} approximates the points of $C$ near the gradient catastrophe point) and factors as
\begin{equation}
    \frac{1}{3}t_\cc u_0'''(y_\cc)\left(y-y_\cc\mp\frac{(t-t_\cc)^{3/2}}{t_\cc u_0'''(y_\cc)^{1/2}}\right)^2\left(y-y_\cc \pm \frac{2(t-t_\cc)^{3/2}}{t_\cc u_0'''(y_\cc)^{1/2}}\right)\approx 0.
\end{equation}
These calculations show that 
the condition $u_0'''(y_\cc)>0$ guarantees that precisely two distinct curves of $C$ emanate from the catastrophe point $(t_\cc,x_\cc)$ in the direction of increasing $t$ forming a \emph{cusped caustic}, see Figure 1 of either \cite{SmallDispersionBO-1,BGM23}.  The branches of $C$ near the catastrophe point therefore are parametrized by two branches of the cusped caustic that we denote by $x_\hrm(t)<x_\srm(t)$ with $x_\hrm(t)$ and $x_\srm(t)$ corresponding to the choice of the minus and plus signs in \eqref{eq:Caustic-approximations} respectively.  Moreover,
\begin{itemize}
    \item when $x=x_\hrm(t)$, there is a double critical point to the right of $y=y_\cc$ that we denote $y=y_\hrm(t)$ and a simple critical point to the left of $y=y_\cc$ that we denote $y=y_2(t)$;
    \item when $x=x_\srm(t)$, there is a double critical point to the left of $y=y_\cc$ that we denote $y=y_\srm(t)$ and a simple critical point to the right of $y=y_\cc$ that we denote $y=y_0(t)$.
\end{itemize}
In the limit $t\downarrow t_\cc$ along either branch, the simple critical point coalesces with the double critical point to form a triple critical point exactly at the gradient catastrophe, as $x_\hrm(t_\cc)=x_\srm(t_\cc)=x_\cc$.

Fix $t>t_\cc$ with $t-t_\cc$ sufficiently small.  Then since $t_\cc>0$ is the earliest time of gradient catastrophe by assumption, for each $x\in\mathbb{R}$ there are at most three real critical points of $h(z;t,x)$, counted with multiplicity. 

\subsubsection{Critical points and contours for $(t,x)$ near the soliton edge}
\label{sec:soliton-edge-contours}
The \emph{soliton edge} refers to the curve $x=x_\srm(t)$.  Given any $\mu>0$, if $x-x_\srm(t)$ is sufficiently small, there will be a simple real critical point $z=y_0(t,x)$ of $h(z;t,x)$ close to $z=y_0(t)$ and a pair of critical points (real and distinct if $x<x_\srm(t)$, a complex-conjugate pair if $x>x_\srm(t)$, and a real double point if $x=x_\srm(t)$) close to $y=y_\srm(t)<y_0(t)$, all of which lie in the horizontal strip $|\mathrm{Im}(z)|<\mu$. The remaining $2N-2$ critical points are all complex-conjugate pairs, lying in the complementary region $|\mathrm{Im}(z)|\ge \mu$; by assuming that these are all simple critical points at the gradient catastrophe point $(t_\cc,x_\cc)$ we ensure that they remain simple for $t\in K=(t_\cc,t_\cc+\delta)$ for $\delta>0$ sufficiently small. For $t\in K$, we may therefore consistently enumerate the critical points in the upper half-plane as $z_2(t,x),\dots,z_N(t,x)$ all of which are real-analytic functions of $x\approx x_\srm(t)$.  The algorithm in \cite{SmallDispersionBO-1} puts these points into a $1-1$ correspondence with contours $W_2,\dots,W_N$ encircling the branch cuts $\Gamma_2,\dots,\Gamma_N$ respectively such that the integral along $W_n$ of $\ee^{-\ii h(z;t,x)/\epsilon}f(z)$ for $f$ analytic on $W_n$ can be expanded by the classical method of steepest descent with the expansion involving only the germ of $f$ at the simple critical point $z=z_n(t,x)$.  

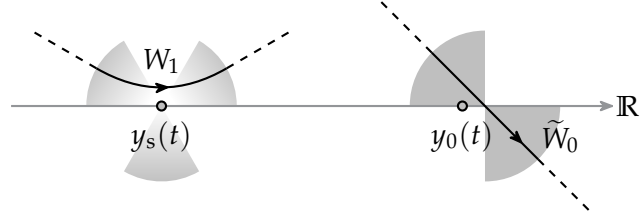
\begin{figure}[h]
\begin{center}
\begin{tikzpicture}
\begin{scope}[shift={(-2cm,0cm)}]
    \fill[lightgray, path fading=east,fading transform={rotate=-30}]
                (0cm,0cm) -- (180:1cm) arc[start angle=180, end angle=120, radius=1cm] -- cycle;
    \fill[lightgray,path fading=east,fading transform={rotate=-150}] (0cm,0cm) -- (0:1cm) arc[start angle=0,end angle=60,radius=1cm] -- cycle;
    \fill[lightgray,path fading=east,fading transform={rotate=90}] (0cm,0cm) -- (-60:1cm) arc[start angle=-60,end angle=-120,radius=1cm] -- cycle;
\end{scope}

\begin{scope}[shift={(2.3cm,0cm)}]
\fill[lightgray] (0cm,0cm) -- (180:1cm) arc[start angle=180, end angle=90,radius=1cm] --cycle;
\fill[lightgray] (0cm,0cm) -- (0:1cm) arc[start angle=0,end angle=-90,radius=1cm] -- cycle;
\end{scope}

    \draw[thick,gray,postaction={frac arrow=1}] (-4cm,0cm)--(4cm,0cm);
\begin{scope}[shift={(-2cm,0cm)}]
\draw[thick,black,postaction={frac arrow=0.55}] (150:1cm) to[out=-30,in=-150] (30:1cm);
\draw[thick,black,dashed] (150:1cm) -- (150:2cm);
\draw[thick,black,dashed] (30:1cm) -- (30:2cm);
\end{scope}

\begin{scope}[shift={(2.3cm,0cm)}]
\draw[thick,black,postaction={frac arrow=0.85}] (135:1cm) -- (-45:1cm);
\draw[thick,black,dashed] (135:1cm) -- (135:2cm);
\draw[thick,black,dashed] (-45:1cm) -- (-45:2cm);
\end{scope}

    \node (ys) at (-2cm,0cm) [thick,black,circle,fill=lightgray,draw=black,inner sep=1.2pt]{};
    \node (y0) at (2cm,0cm) [thick,black,circle,fill=lightgray,draw=black,inner sep=1.2pt]{};

\node at (-2cm,-0.4cm) {$y_\srm(t)$};
\node at (2cm,-0.4cm) {$y_0(t)$};

\node at (-2cm,0.6cm) {$W_1$};
\node at (3.3cm,-0.4cm) {$\widetilde{W}_0$};

\node at (4.2cm,0cm) {$\mathbb{R}$};
\end{tikzpicture}
\end{center}
\caption{The contours $\widetilde{W}_0$ and $W_1$ near the points $y_0(t)$ and $y_\srm(t)$ when $x\approx x_\srm(t)$.  $\widetilde{W}_0$ passes a real simple critical point $z=y_0(t,x)$ that converges to $y_0(t)\in\mathbb{R}$ as $x\to x_\srm(t)$.  The regions where $\mathrm{Re}(-\ii h(z;t,x))<0$ near the points $z=y_\srm(t)$ and $z=y_0(t)$ are indicated with shading, and we deliberately faded the shading near $z=y_\srm(t)$ because the details near that point depend on the sign of $x-x_\srm(t)$.  Those details are taken care of in the Chester-Friedman-Ursell method by finding a cubic normal form for the exponent near $z=y_\srm(t)$.}
\label{fig:soliton-edge-contours}
\end{figure}

The algorithm in \cite{SmallDispersionBO-1} specifies the contour $W_1$ differently depending on whether $x<x_\srm(t)$ or $x>x_\srm(t)$, however in both cases the contour enters and exits with increasing $\mathrm{Re}(z)$ the half-plane $\mathrm{Im}(z)<\mu$ near $z=y_\srm(t)$ and the dominant contributions to integrals over $W_1$ come from the pair of critical points in a neighborhood of $y_\srm(t)$.  In our analysis for $x\approx x_\srm(t)$ we will not attempt to thread $W_1$ over any particular critical points, but instead deal with the expansion of integrals of $\ee^{-\ii h(z;t,x)/\epsilon}f(z)$ over $W_1$ using a method of Chester, Friedman, and Ursell \cite{ChesterFriedmanUrsell} that allows for coalescing critical points.  The contour $W_0$ specified by the algorithm in \cite{SmallDispersionBO-1} also enters and exits with increasing $\mathrm{Re}(z)$ the half-plane $\mathrm{Im}(z)<\mu$ near $z=y_\srm(t)$, after which it crosses the real line one final time in the diagonal direction $\arg(\dd z)=-\frac{1}{4}\pi$ at the simple real critical point denoted $y_0(t,x)$ lying to the right of $y_\srm(t)$.  By the same elementary row operations on the first two rows of both determinants $N$ and $D$, we may replace $W_0$ with $\widetilde{W}_0:=W_0-W_1$ (understood as concatenation of contours; the actual row operations may involve an additional multiplier on the second row if the critical points lie on different sheets of the multivalued logarithm in the definition \eqref{eq:h-def} of $h(z;t,x)$ by analytic continuation).  Then as for the contours $W_2,\dots,W_N$, integrals over the contour $\widetilde{W}_0$ of $\ee^{-\ii h(z;t,x)/\epsilon}f(z)$ can be expanded by the classical steepest descent method with dominant contributions from the germ of $f$ at $z=y_0(t,x)$.  Allowing $x-x_\srm(t)$ to be small of either sign, the landscape of $\mathrm{Re}(-\ii h(z;t,x))$ traversed by the contour $W_1$ is somewhat ambiguous near $z=y_\srm(t)$, but it is precisely this ambiguity that is resolved by the Chester-Friedman-Ursell method.  Because when $x=x_\srm(t)$ we have $h'(y_\srm(t);t,x_\srm(t))=h''(y_\srm(t);t,x_\srm(t))=0$ and $h'''(y_\srm(t);t,x_\srm(t))<0$, the only important properties of $W_1$ near $y_\srm(t)$ are that it should enter the neighborhood of this point from the direction $\arg(z-y_\srm(t))\approx\frac{5}{6}\pi$ and exit in the direction $\arg(z-y_\srm(t))\approx\frac{1}{6}\pi$.  See Figure~\ref{fig:soliton-edge-contours}.

\subsubsection{Critical points and contours for $(t,x)$ near the harmonic edge}
\label{sec:harmonic-edge-contours}
The \emph{harmonic edge} refers to the curve $x=x_\hrm(t)$.  When $x-x_\hrm(t)$ is small, there is a simple real critical point of $h(z;t,x)$ denoted $z=y_2(t,x)$ near $z=y_2(t)$ and a pair of critical points near $z=y_\hrm(t)$ forming a complex-conjugate pair for $x<x_\hrm(t)$ and a real pair for $x>x_\hrm(t)$.  All other critical points of $h$ are in complex-conjugate pairs away from the real line.  Assuming that the complex critical points are all simple for $(t,x)=(t_\cc,x_\cc)$ ensures that they remain simple for $t\in K=(t_\cc,t_\cc+\delta)$ for $\delta>0$ sufficiently small, so we enumerate them as $z_2(t,x),\dots,z_N(t,x)$, and use the algorithm in \cite{SmallDispersionBO-1} to specify corresponding contours $W_2,\dots,W_N$ as described in Section~\ref{sec:soliton-edge-contours}.

\begin{figure}[h]
\begin{center}
\begin{tikzpicture}
\begin{scope}[shift={(-1.7cm,0cm)}]
\fill[lightgray] (0cm,0cm) -- (180:1cm) arc[start angle=180, end angle=90,radius=1cm] --cycle;
\fill[lightgray] (0cm,0cm) -- (0:1cm) arc[start angle=0,end angle=-90,radius=1cm] -- cycle;
\end{scope}

\draw[thick,gray,postaction={frac arrow=1}] (-4cm,0cm)--(4cm,0cm);

\begin{scope}[shift={(2cm,0cm)}]
    \fill[lightgray, path fading=east,fading transform={rotate=-30+60}]
                (0cm,0cm) -- (180+60:1cm) arc[start angle=180+60, end angle=120+60, radius=1cm] -- cycle;
    \fill[lightgray,path fading=east,fading transform={rotate=-150+60}] (0cm,0cm) -- (0+60:1cm) arc[start angle=0+60,end angle=60+60,radius=1cm] -- cycle;
    \fill[lightgray,path fading=east,fading transform={rotate=90+60}] (0cm,0cm) -- (-60+60:1cm) arc[start angle=-60+60,end angle=-120+60,radius=1cm] -- cycle;
\end{scope}

\begin{scope}[shift={(-1.7cm,0cm)}]

\draw[thick,black] (135:1cm) -- (-45:1cm);
\draw[thick,dashed,black] (135:1cm) -- (135:2cm);
\draw[thick,black,postaction={frac arrow=1}] (-45:1cm) to[out=-45,in=180] (1.5cm,-1.0cm);
\end{scope}

\begin{scope}[shift={(2.0cm,0cm)}]
\draw[thick,black]
(-2.2cm,-1.0cm) to[out=0,in=-150] (-150:1cm);
\draw[thick,black] (-150:1cm) to[out=30,in=-90] (90:1cm);
\draw[thick,black,dashed] (90:1cm) -- (90:2cm);

\draw[thick,black,postaction={frac arrow=0.8}] (90:1cm) to[out=-90,in=150] (-30:1cm);
\draw[thick,black,dashed] (-30:1cm) -- (-30:2cm);
\end{scope}

    \node (yd) at (2cm,0cm) [thick,black,circle,fill=lightgray,draw=black,inner sep=1.2pt]{};
    \node (y2) at (-2cm,0cm) [thick,black,circle,fill=lightgray,draw=black,inner sep=1.2pt]{};

\node at (2cm,-0.4cm) {$y_\hrm(t)$};
\node at (-2cm,-0.4cm) {$y_2(t)$};

\node at (-0.3cm,-0.7cm) {$W_1$};
\node at (3.3cm,-0.4cm) {$\widetilde{W}_0$};

\node at (4.2cm,0cm) {$\mathbb{R}$};

\end{tikzpicture}
\end{center}
\caption{The contours $\widetilde{W}_0$ and $W_1$ near the points $z=y_2(t)$ and $z=y_\hrm(t)$ when $x\approx x_\hrm(t)$.}    
\label{fig:harmonic-edge-contours}
\end{figure}
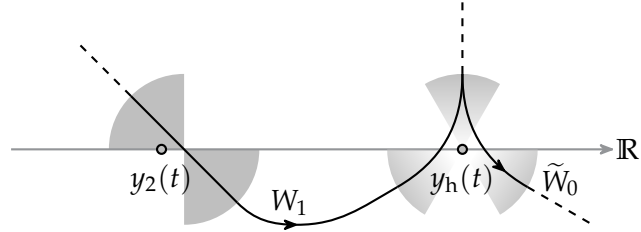

The difference between this case and the case $x\approx x_\srm(t)$ considered in Section~\ref{sec:soliton-edge-contours} is that the simple real critical point $z=y_2(t,x)$ now lies to the left of the pair of critical points close to $y_\hrm(t)$.  The algorithm in \cite{SmallDispersionBO-1} then constructs the contour $W_1$ differently depending on whether or not the critical points close to $y_\hrm(t)$ are real or not.  For the purposes of having a setup that we may take to be qualitatively independent of the sign of $x-x_\hrm(t)$, we will take $W_1$ to resemble that built by the algorithm when the critical points near $y_\hrm(t)$ are both real.  In that situation, $W_1$ first enters the half-plane $\mathrm{Im}(z)<\mu$ with increasing $\mathrm{Re}(z)$ just to the left of the left-most simple saddle point $z=y_2(t,x)$, which it traverses with the angle $\arg(\dd z)=-\frac{1}{4}\pi$; then it traverses the central simple saddle point near $y_\hrm(t)$ with angle $\arg(\dd z)=\frac{1}{4}\pi$ after which it exits the half-plane $\mathrm{Im}(z)<\mu$, never to return.  Since the coalescence of the critical points near $y_\hrm(t)$ will be handled by the Chester-Friedman-Ursell method, we need not worry about the details of how $W_1$ behaves in a neighborhood of this point, except that it should enter and exit this neighborhood in directions conducive to exponential decay of $\ee^{-\ii h(z;t,x)/\epsilon}$ as $\epsilon\downarrow 0$.  Based on the landscape of $\mathrm{Re}(-\ii h(z;t,x))$ when $x=x_\hrm(t)$ which is locally determined from the properties that $h'(y_\hrm(t);t,x_\hrm(t))=h''(y_\hrm(t);t,x_\hrm(t))=0$ and $h'''(y_\hrm(t);t,x_\hrm(t))>0$, $W_1$ should enter the neighborhood of $y_\hrm(t)$ from the direction $\arg(z-y_\hrm(t))=-\frac{5}{6}\pi$ and exit in the direction $\arg(z-y_\hrm(t))=\frac{1}{2}\pi$.  The contour $W_0$ specified by the algorithm in \cite{SmallDispersionBO-1} follows a similar path as $W_1$ into the neighborhood of $y_\hrm(t)$, but it exits that neighborhood in the direction $\arg(z-y_\hrm(t))=-\frac{1}{6}\pi$ and ultimately terminates at $z=\infty\ee^{-\ii\pi/4}$.  Then exactly as in Section~\ref{sec:soliton-edge-contours}, we use elementary row operations to replace $W_0$ with $\widetilde{W}_0=W_0-W_1$, so that now $\widetilde{W}_0$ only visits the neighborhood of $z=y_\hrm(t)$ which it enters from the direction $\arg(z-y_\hrm(t))=\frac{1}{2}\pi$ and exits in the direction $\arg(z-y_\hrm(t))=-\frac{1}{6}\pi$.  See Figure~\ref{fig:harmonic-edge-contours}.

\subsubsection{Critical points and contours for $(t,x)$ near the gradient catastrophe point}
\label{sec:catastrophe-contours}
\begin{figure}[h]
\begin{center}
    \begin{tikzpicture}
    \fill[lightgray, path fading=east,fading transform={rotate=-22.5}]
                (0cm,0cm) -- (-4cm,0cm) -- (-4cm,4cm) -- cycle;
    \fill[lightgray, path fading=east,fading transform={rotate=-112.5}] (0cm,0cm) -- (0cm,4cm) -- (4cm,4cm) -- cycle;
    \fill[lightgray, path fading=east,fading transform={rotate=180-22.5}] (0cm,0cm) -- (4cm,0cm) -- (4cm,-4cm) -- cycle;
    \fill[lightgray, path fading=east,fading transform={rotate=90-22.5}] (0cm,0cm) -- (0cm,-4cm) -- (-4cm,-4cm) -- cycle;
    
    \draw[thick,gray,postaction={frac arrow=1}] (-4cm,0cm)--(4cm,0cm);
\node at (0cm,-0.4cm) {$y_\cc$};
\node at (4.2cm,0cm) {$\mathbb{R}$};
\draw[thick,black,postaction={frac arrow=0.5}] (180-22.5:3cm) to[out=-22.5,in=-90-22.5] (90-22.5:3cm);
\draw[thick,black,dashed] (180-22.5:4cm) -- (180-22.5:3cm);
\draw[thick,black,dashed] (90-22.5:4cm) -- (90-22.5:3cm);
\draw[thick,black,postaction={frac arrow=0.7}] (180-22.5:3cm) -- (-22.5:3cm);
\draw[thick,black,dashed] (180-22.5:3cm) -- (180-22.5:4cm);
\draw[thick,black,dashed] (-22.5:3cm) -- (-22.5:4cm);
    \node (yc) at (0cm,0cm) [thick,black,circle,fill=lightgray,draw=black,inner sep=1.2pt]{};
    \node at (-0.8cm,1.5cm) {$W_1$};
    \node at (1cm,-0.8cm) {$W_0$};
\end{tikzpicture}
\end{center}
\caption{The contours $W_0$ and $W_1$ in the neighborhood of the point $y_\cc$ when $(t,x)\approx (t_\cc,x_\cc)$.}
\label{fig:catastrophe-contours}
\end{figure}
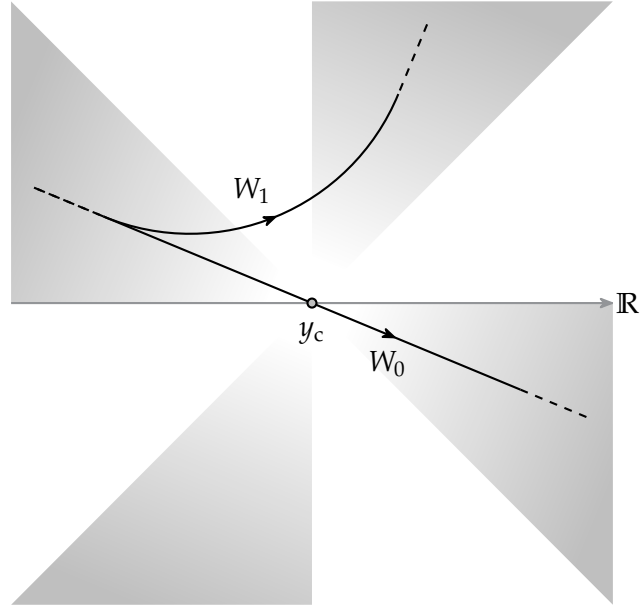

Finally, we consider the situation that $(t,x)\approx (t_\cc,x_\cc)$, which is the vertex of the cusped caustic formed by $x=x_\srm(t)$ and $x=x_\hrm(t)$ for $t>t_\cc$.  In this situation, $h(z;t,x)$ has three critical points very close to $z=y_\cc\in\mathbb{R}$, and all other critical points are in complex-conjugate pairs outside the strip $|\mathrm{Im}(z)|<\mu$.  We assume that the latter critical points are all simple when $(t,x)=(t_\cc,x_\cc)$, and hence it remains the case for nearby $(t,x)$.  Then as in Sections~\ref{sec:soliton-edge-contours} and \ref{sec:harmonic-edge-contours}, the algorithm from \cite{SmallDispersionBO-1} assigns contours $W_2,\dots,W_N$ each enclosing a branch cut $\Gamma_n$ and passing respectively over the complex critical points $z_2(t,x),\dots,z_N(t,x)$ with $\mathrm{Im}(z_n(t,x))\ge \mu$ that we may consistently enumerate for $(t,x)\approx (t_\cc,x_\cc)$.  The remaining contours $W_1$ and $W_0$ are difficult to precisely specify near $z=y_\cc$, but since we are going to study the integrals over these contours using a generalization of the Chester-Friedman-Ursell method, these details are not so important.  Based on the landscape of $\mathrm{Re}(-\ii h(z;t,x))$ near $z=y_\cc$ when $(t,x)=(t_\cc,x_\cc)$, where $h'(y_\cc;t_\cc,x_\cc)=h''(y_\cc;t_\cc,x_\cc)=h'''(y_\cc;t_\cc,x_\cc)=0$ but $h''''(y_\cc;t_\cc,x_\cc)>0$, we only insist that $W_1$ should enter the neighborhood of $z=y_\cc$ from the direction $\arg(z-y_\cc)=\frac{7}{8}\pi$ and exit in the direction $\arg(z-y_\cc)=\frac{3}{8}\pi$, and that $W_0$ should enter from the direction $\arg(z-y_\cc)=\frac{7}{8}\pi$ and exit in the opposite direction $\arg(z-y_\cc)=-\frac{1}{8}\pi$.  See Figure~\ref{fig:catastrophe-contours}.

\subsection{Structure of the proofs}
\label{sec:structure-of-proofs}
The proofs of Theorems~\ref{thm:soliton-edge}, \ref{thm:harmonic-edge}, and \ref{thm:BO-catastrophe} consist of three roughly parallel steps.  After starting from the exact solution formula given in Theorem~\ref{thm:IVP-solution} with contours specified as described in Sections~\ref{sec:soliton-edge-contours}, \ref{sec:harmonic-edge-contours}, and \ref{sec:catastrophe-contours} respectively, the steps are the following.
\begin{enumerate}
    \item Construct a suitable conformal map on a neighborhood of the coalescing critical points of $h(z;t,x)$ with the aim of representing $h(z;t,x)$ in an appropriate normal form.  
    \label{step:normal-form}
    \item Expand the matrix elements of $\mathbf{A}(t,x;\epsilon)$ and $\mathbf{B}(t,x;\epsilon)$ using the classical steepest descent method for the contributions from simple critical points of $h(z;t,x)$ and by a generalization based on the normal-form exponent in the case of coalescing critical points.  In the case that only two critical points coalesce, this is the method of Chester-Friedman-Ursell \cite{ChesterFriedmanUrsell}, and when three critical points coalesce (the gradient catastrophe case) a further generalization is required.
    \label{step:expand}
    \item The previous step results yields leading approximations of the matrix elements with error terms that are controlled in the small-dispersion limit $\epsilon\downarrow 0$.  We then use multilinearity and special determinant identities to further simplify the ratio of determinants $N/D$.
    \label{step:algebra}
\end{enumerate}

\subsubsection{Conformal mappings and normal forms of exponents}

Here we establish a key lemma that will be used in Step~\ref{step:normal-form}.
\begin{lemma}[Normal forms]\label{lem:normal-form}
    Let $\mathbf{x}\in\mathbb{R}^m$ be a vector of parameters, and suppose that $z\mapsto H(z;\mathbf{x})$ is a function analytic in $z$ on a neighborhood of $z_0\in\mathbb{C}$ and $C^\infty$ in $\mathbf{x}$.  If there is a point $\mathbf{x}_0\in\mathbb{R}^m$ such that 
    \begin{equation}
        \frac{\dd^kH}{\dd z^k}(z_0;\mathbf{x}_0)=0,\quad k=1,\dots,n,\quad\text{and}\quad\frac{\dd^{n+1}H}{\dd z^{n+1}}(z_0;\mathbf{x}_0)\neq 0,
        \label{eq:degeneration-description}
    \end{equation}
    then for each $\nu\neq 0$ satisfying
    \begin{equation}
        \nu^{n+1}=\frac{1}{n!}\frac{\dd^{n+1}H}{\dd z^{n+1}}(z_0;\mathbf{x}_0),
        \label{eq:nu-power}
    \end{equation}
    there exists a unique conformal mapping $w=\mathcal{W}(z;\mathbf{x})$ univalent on a neighborhood $D$ of $z_0$ and of class $C^\infty$ on a neighborhood $\mathcal{X}$ of $\mathbf{x}_0$ in the parameter space, and coefficients $z_0(\mathbf{x}),\dots,z_{n-1}(\mathbf{x})$ of class $C^\infty(\mathcal{X})$ such that
    \begin{equation}
        H(z;\mathbf{x})-H(z_0;\mathbf{x})=P(w;\mathbf{x}):=\frac{w^{n+1}}{n+1}+\sum_{k=0}^{n-1}z_k(\mathbf{x})w^k,\quad w=\mathcal{W}(z;x)
        \label{eq:desired-identity}
    \end{equation}
    holds whenever $z\in D$ and $\mathbf{x}\in \mathcal{X}$, and $\mathcal{W}'(z_0;\mathbf{x}_0)=\nu$.
\end{lemma}

The polynomial $P(w;\mathbf{x})$ on the right-hand side is a local \emph{normal form} for $h$.  The problem of finding a normal form requires determining both the coefficients $z_0(\mathbf{x}),\dots,z_{n-1}(\mathbf{x})$ and the conformal mapping $w=\mathcal{W}(z;\mathbf{x})$ so that \eqref{eq:desired-identity} holds.  The coefficients $z_0(\mathbf{x}),\dots,z_{n-1}(\mathbf{x})$ are said to unfold a \emph{cuspoid catastrophe of codimension $n-1$}, see \cite[\S 36.2]{DLMF}.  In the case $n=2$, Lemma~\ref{lem:normal-form} was proved and applied in the paper of Chester, Friedman, and Ursell \cite{ChesterFriedmanUrsell}.  In the more general setting $n\ge 2$, this result was proved by Levinson \cite[Theorem 1.2]{Levinson} in the category of analytic dependence on the parameters $\mathbf{x}$ using an induction argument.  There is no loss of generality in the assumption that the coefficient of $w^n$ of $P(w;\mathbf{x})$ on the right-hand side of \eqref{eq:desired-identity} is zero, since otherwise it can be removed by a constant shift $\mathcal{W}(z;\mathbf{x})\mapsto \mathcal{W}(z;\mathbf{x})+c(\mathbf{x})$.  We now give a proof of Lemma~\ref{lem:normal-form} that is different in character from those in \cite{ChesterFriedmanUrsell} and \cite{Levinson}. 

\begin{proof}[Proof of Lemma~\ref{lem:normal-form}]
Imposing the conditions $z_k(\mathbf{x}_0)=0$, $k=0,\dots, n-1$, is reasonable in light of \eqref{eq:degeneration-description}.  Doing so and setting $\mathbf{x}=\mathbf{x}_0$ in \eqref{eq:desired-identity} yields 
\begin{equation}
\begin{split}
    \mathcal{W}(z;\mathbf{x}_0)^{n+1}&=(n+1)\left(H(z;\mathbf{x}_0)-H(z_0;\mathbf{x}_0)\right)\\
    &=\frac{1}{n!}\frac{\dd^{n+1}H}{\dd z^{n+1}}(z_0;\mathbf{x}_0)(z-z_0)^{n+1} + \frac{n+1}{(n+2)!}\frac{\dd^{n+2}H}{\dd z^{n+2}}(z_0;\mathbf{x}_0)(z-z_0)^{n+2}\\
    &\qquad\qquad\qquad\qquad{}+O((z-z_0)^{n+3})\\
    &=\nu^{n+1}(z-z_0)^{n+1}\left(1+\frac{n+1}{\nu^{n+1}(n+2)!}\frac{\dd^{n+2}H}{\dd z^{n+2}}(z_0;\mathbf{x}_0)(z-z_0)+O((z-z_0)^2)\right),
\end{split}
\end{equation}
where the error term indicates a convergent power series in $z-z_0$. Therefore, for each $\nu\neq 0$ solving \eqref{eq:nu-power}, by taking an analytic root, we get an analytic solution
\begin{equation}
\begin{split}
    \mathcal{W}(z;\mathbf{x}_0)&=\nu(z-z_0)\left(1+\frac{1}{\nu^{n+1}(n+2)!}\frac{\dd^{n+2}H}{\dd z^{n+2}}(z_0;\mathbf{x}_0)(z-z_0)+O((z-z_0)^2)\right)\\
    &=\nu(z-z_0) + c(z-z_0)^2 + O((z-z_0)^3),\quad c:=\frac{1}{\nu^n(n+2)!}\frac{\dd^{n+2}H}{\dd z^{n+2}}(z_0;\mathbf{x}_0),
\end{split}
\label{eq:W-general-crit}
\end{equation}
which is univalent near $z_0$.  There are then obviously $n+1$ conformal maps such that \eqref{eq:desired-identity} holds for $\mathbf{x}=\mathbf{x}_0$ with $z_k(\mathbf{x}_0)=0$, $k=0,\dots,n-1$.  The rest of the proof will show that there are $n+1$ corresponding conformal maps uniquely defined for $\mathbf{x}\in \mathcal{X}$ such that \eqref{eq:desired-identity} holds for $n+1$ corresponding normal forms $P(w;\mathbf{x})$.

Suppose next that for $\mathbf{x}\neq \mathbf{x}_0$ there is a univalent solution $\mathcal{W}(z;\mathbf{x})$ of \eqref{eq:desired-identity} for suitable coefficients $z_0(\mathbf{x}),\dots,z_{n-1}(\mathbf{x})$ vanishing at $\mathbf{x}=\mathbf{x}_0$ and depending smoothly on $\mathbf{x}$ in $\mathcal{X}$.  Let $\partial=d_1(\mathbf{x})\partial_{x_1} + \cdots + d_m(\mathbf{x})\partial_{x_m}$ denote a smooth vector field on $\mathbb{R}^m$.  Fixing $z\in D$, applying $\partial$ to \eqref{eq:desired-identity} gives
\begin{equation}
    \partial H(z;\mathbf{x})-\partial H(z_0;\mathbf{x})=P'(\mathcal{W}(z;\mathbf{x});\mathbf{x})\partial\mathcal{W}(z;\mathbf{x}) +\sum_{k=0}^{n-1}\partial z_k(\mathbf{x})\cdot \mathcal{W}(z;\mathbf{x})^{k}.
    \label{eq:derivative-z}
\end{equation}
Using the inverse mapping $z=\mathcal{W}^{-1}(w;\mathbf{x})$, which exists by univalence, this becomes
\begin{equation}
    \partial H(\mathcal{W}^{-1}(w;\mathbf{x});\mathbf{x})-\partial H(z_0;\mathbf{x})=P'(w;\mathbf{x})\partial \mathcal{W}(\mathcal{W}^{-1}(w;\mathbf{x});\mathbf{x}) +\sum_{k=0}^{n-1}w^k\partial z_k(\mathbf{x}).
    \label{eq:derivative-w}
\end{equation}
Let $w_1(\mathbf{x}),\dots,w_n(\mathbf{x})$ denote the roots of $w\mapsto P'(w;\mathbf{x})$, i.e., the critical points of the normal form $P(w;\mathbf{x})$.  Assuming they are all simple, which holds on a subset of $\mathcal{X}\subset\mathbb{R}^m$ of codimension 1, we evaluate \eqref{eq:derivative-w} at each of these roots and obtain an invertible Vandermonde system for the derivatives of the coefficients $z_0(\mathbf{x}),\dots,z_{n-1}(\mathbf{x})$:
\begin{equation}
    \begin{bmatrix}
        1 & w_1(\mathbf{x}) & \cdots & w_1(\mathbf{x})^{n-1}\\
        \vdots & \vdots & \ddots & \vdots\\
        1 & w_n(\mathbf{x}) & \cdots & w_n(\mathbf{x})^{n-1}
    \end{bmatrix}
    \begin{bmatrix}\partial z_0(\mathbf{x})\\\vdots \\ \partial z_{n-1}(\mathbf{x})
    \end{bmatrix} = \begin{bmatrix}\partial H(\mathcal{W}^{-1}(w_1(\mathbf{x});\mathbf{x});\mathbf{x})\\\vdots \\
    \partial H(\mathcal{W}^{-1}(w_n(\mathbf{x});\mathbf{x});\mathbf{x})\end{bmatrix} -\partial H(z_0;\mathbf{x})\begin{bmatrix}1\\\vdots \\ 1\end{bmatrix}.
\label{eq:Vandermonde-system}
\end{equation}
Using univalence of $\mathcal{W}(z;\mathbf{x})$ once again, we can express the first term on the right-hand side as a contour integral:
\begin{equation}
    \begin{bmatrix}\partial H(\mathcal{W}^{-1}(w_1(\mathbf{x});\mathbf{x});\mathbf{x})\\\vdots \\
    \partial H(\mathcal{W}^{-1}(w_n(\mathbf{x});\mathbf{x});\mathbf{x})\end{bmatrix} = 
    \frac{1}{2\pi\ii}\oint_{L}\partial H(z;\mathbf{x}) \mathcal{W}'(z;\mathbf{x})\begin{bmatrix} (\mathcal{W}(z;\mathbf{x})-w_1(\mathbf{x}))^{-1}\\\vdots\\(\mathcal{W}(z;\mathbf{x})-w_n(\mathbf{x}))^{-1}\end{bmatrix}\,\dd z,
\end{equation}
where $L$ is a positively-oriented Jordan curve in $D$ that encloses all of the preimages of the $w_j(\mathbf{x})$ under $z\mapsto \mathcal{W}(z;\mathbf{x})$.  We may fix $L$, and this condition will hold for $\|\mathbf{x}-\mathbf{x}_0\|$ sufficiently small.  Therefore, by linearity it is enough to invert the Vandermonde matrix on two vectors:  $[1,\dots,1]^\top$ and $[(w-w_1(\mathbf{x}))^{-1},\dots,(w-w_n(\mathbf{x}))^{-1}]^\top$, $w=\mathcal{W}(z;\mathbf{x})$. By Cramer's rule (or just by inspection), it is easy to see that 
\begin{equation}
  \begin{bmatrix}
        1 & w_1(\mathbf{x}) & \cdots & w_1(\mathbf{x})^{n-1}\\
        \vdots & \vdots & \ddots & \vdots\\
        1 & w_n(\mathbf{x}) & \cdots & w_n(\mathbf{x})^{n-1}
    \end{bmatrix}^{-1}\begin{bmatrix}1\\\vdots\\ 1\end{bmatrix}=\begin{bmatrix}1 \\0\\\vdots\\0\end{bmatrix}.  
\end{equation}
Similarly, denoting 
\begin{equation}
    R_k(w;\mathbf{x}):=\left(\begin{bmatrix}
        1 & w_1(\mathbf{x}) & \cdots & w_1(\mathbf{x})^{n-1}\\
        \vdots & \vdots & \ddots & \vdots\\
        1 & w_n(\mathbf{x}) & \cdots & w_n(\mathbf{x})^{n-1}
    \end{bmatrix}^{-1}\begin{bmatrix}(w-w_1(\mathbf{x}))^{-1}\\\vdots\\(w-w_n(\mathbf{x}))^{-1}\end{bmatrix}\right)_{k+1},\quad k=0,\dots,n-1,
\end{equation}
Cramer's rule yields
\begin{equation}
    R_k(w;\mathbf{x})=\frac{\displaystyle\begin{vmatrix} 1 & \cdots & w_1(\mathbf{x})^{k-1} & (w-w_1(\mathbf{x}))^{-1} & w_1(\mathbf{x})^{k+1}&\cdots & w_1(\mathbf{x})^{n-1}\\
    \vdots & \ddots & \vdots & \vdots & \vdots & \ddots & \vdots\\
    1 & \cdots & w_n(\mathbf{x})^{k-1} & (w-w_n(\mathbf{x}))^{-1} & w_n(\mathbf{x})^{k+1} & \cdots & w_n(\mathbf{x})^{n-1}\end{vmatrix}}{\displaystyle\begin{vmatrix} 1 & w_1(\mathbf{x}) & \cdots & w_1(\mathbf{x})^{n-1}\\
    \vdots & \vdots & \ddots & \vdots \\
    1 & w_n(\mathbf{x}) & \cdots & w_n(\mathbf{x})^{n-1}\end{vmatrix}}.
\end{equation}
By multiplying $P'(w;\mathbf{x})=(w-w_1(\mathbf{x}))\cdots (w-w_n(\mathbf{x}))$ into column $k+1$ of the determinant in the numerator it is clear that $R_k(w;\mathbf{x})P'(w;\mathbf{x})$ is a polynomial in $w$.  Also, by expanding the determinant in the numerator for large $w$ using geometric series $(w-w_j(\mathbf{x}))^{-1}=w^{-1}+w_j(\mathbf{x})w^{-2} + w_j(\mathbf{x})^2w^{-3}+\cdots$, one checks that for $k=0,\dots,n-1$
\begin{equation}
    R_k(w;\mathbf{x})=w^{-(k+1)} + O(w^{-(n+1)}),\quad w\to\infty.
\end{equation}
Since $P'(w;\mathbf{x})$ is a polynomial of degree $n$, as $w\to\infty$ we have
\begin{equation}
    R_k(w;\mathbf{x})P'(w;\mathbf{x})=w^{-(k+1)}P'(w;\mathbf{x}) + O(w^{-1})=w^{n-k-1} + \sum_{\ell=1}^{n-1}\ell z_\ell(\mathbf{x})w^{\ell-k-2} + O(w^{-1}),
\end{equation}
where we substituted the definition of $P(w;\mathbf{x})$ from \eqref{eq:desired-identity}.  However, since the left-hand side is a polynomial in $w$, we have the exact identity
\begin{equation}
    R_k(w;\mathbf{x})P'(w;\mathbf{x})=w^{n-k-1} + \sum_{\ell=k+2}^{n-1}\ell z_\ell(\mathbf{x})w^{\ell-k-2},\quad k=0,\dots,n-1.
\end{equation}
Therefore, also
\begin{equation}
    R_k(w;\mathbf{x})=\frac{\displaystyle w^{n-k-1} + \sum_{\ell=k+2}^{n-1}\ell z_\ell(\mathbf{x})w^{\ell-k-2}}{\displaystyle w^n + \sum_{\ell=1}^{n-1}\ell z_\ell(\mathbf{x})w^{\ell-1}}.
    \label{eq:Rk-rational}
\end{equation}
Combining these results shows that the solution of \eqref{eq:Vandermonde-system} is
\begin{equation}
\begin{split}
    \partial z_0(\mathbf{x})&=\frac{1}{2\pi\ii}\oint_L\partial H(z;\mathbf{x})\mathcal{W}'(z;\mathbf{x})R_0(\mathcal{W}(z;\mathbf{x});\mathbf{x})\,\dd z -\partial H(z_0;\mathbf{x})\\
    &=\frac{1}{2\pi\ii}\oint_L\partial H(z;\mathbf{x})\mathcal{W}'(z;\mathbf{x})\frac{\displaystyle \mathcal{W}(z;\mathbf{x})^{n-1} +\sum_{\ell=2}^{n-1}\ell z_\ell(\mathbf{x})\mathcal{W}(z;\mathbf{x})^{\ell-2}}{\displaystyle \mathcal{W}(z;\mathbf{x})^n+\sum_{\ell=1}^{n-1}\ell z_\ell(\mathbf{x})\mathcal{W}(z;\mathbf{x})^{\ell-1}}\,\dd z -\partial H(z_0;\mathbf{x})
    \end{split}
    \label{eq:partial-a0}
\end{equation}
and, for $k=1,\dots,n-1$,
\begin{equation}
    \begin{split}
        \partial z_k(\mathbf{x})&=\frac{1}{2\pi\ii}\oint_L\partial H(z;\mathbf{x})\mathcal{W}'(z;\mathbf{x})R_k(\mathcal{W}(z;\mathbf{x});\mathbf{x})\,\dd z\\
        &=\frac{1}{2\pi\ii}\oint_L\partial H(z;\mathbf{x})\mathcal{W}'(z;\mathbf{x})\frac{\displaystyle \mathcal{W}(z;\mathbf{x})^{n-k-1} + \sum_{\ell=k+2}^{n-1}\ell z_\ell(\mathbf{x})\mathcal{W}(z;\mathbf{x})^{\ell-k-2}}{\displaystyle \mathcal{W}(z;\mathbf{x})^n+\sum_{\ell=1}^{n-1}\ell z_\ell(\mathbf{x})\mathcal{W}(z;\mathbf{x})^{\ell-1}}\,\dd z.
    \end{split}
    \label{eq:partial-ak}
\end{equation}
These formul\ae\ for derivatives of the coefficients $z_0(\mathbf{x}),\dots,z_{n-1}(\mathbf{x})$, while derived assuming simplicity of the roots of $P'(w;\mathbf{x})$, clearly admit continuation to the discriminant locus of $P'(w;\mathbf{x})$ within $\mathcal{X}$, since the roots have been eliminated in favor of the coefficients themselves.  In particular, all of the roots coincide when $\mathbf{x}=\mathbf{x}_0$ and hence $z_k(\mathbf{x}_0)=0$ for $k=0,\dots,n-1$.  In this situation, we obtain
\begin{equation}
    \partial z_0(\mathbf{x}_0)=\frac{1}{2\pi\ii}\oint_L\frac{\partial H(z;\mathbf{x}_0)\mathcal{W}'(z;\mathbf{x}_0)}{\mathcal{W}(z;\mathbf{x}_0)}\,\dd z-\partial H(z_0;\mathbf{x}_0)
    \label{eq:partial-a0-crit}
\end{equation}
and for $k=1,\dots,n-1$,
\begin{equation}
    \partial z_k(\mathbf{x}_0)=\frac{1}{2\pi\ii}\oint_L\frac{\partial H(z;\mathbf{x}_0)\mathcal{W}'(z;\mathbf{x}_0)}{\mathcal{W}(z;\mathbf{x}_0)^{k+1}}\,\dd z.
    \label{eq:partial-ak-crit}
\end{equation}
Since $\mathcal{W}(z;\mathbf{x}_0)$ is fully determined once a solution $\nu$ of \eqref{eq:nu-power} is selected and all coefficients of its power series about $z=z_0$ can be found systematically, the derivatives \eqref{eq:partial-a0-crit}--\eqref{eq:partial-ak-crit} at $\mathbf{x}=\mathbf{x}_0$ can be easily evaluated by residues at $z=z_0$.

These calculations all assumed existence of $\mathcal{W}(z;\mathbf{x})$ and $z_0(\mathbf{x}),\dots,z_{n-1}(\mathbf{x})$ depending smoothly on $\mathbf{x}\in \mathcal{X}$ and such that \eqref{eq:desired-identity} is an identity.  To prove that this is the case, we return to \eqref{eq:derivative-z} and divide by $P'(\mathcal{W}(z;\mathbf{x});\mathbf{x})$ to obtain 
\begin{equation}
\partial \mathcal{W}(z;\mathbf{x})=\frac{\displaystyle \partial H(z;\mathbf{x})-\partial H(z_0;\mathbf{x})-\sum_{k=0}^{n-1}\partial z_k(\mathbf{x})\cdot \mathcal{W}(z;\mathbf{x})^k}{P'(\mathcal{W}(z;\mathbf{x});\mathbf{x})}.    
\label{eq:partial-W}
\end{equation}
Note that if $\partial z_0(\mathbf{x}),\dots,\partial z_{n-1}(\mathbf{x})$ are eliminated using \eqref{eq:partial-a0}--\eqref{eq:partial-ak}, the apparent singularities at the roots of $P'(\mathcal{W}(z;\mathbf{x});\mathbf{x})$ are seen to be removable in the rewritten expression.  Hence the rewritten right-hand side of \eqref{eq:partial-W} defines an analytic function of $z\in D$ that depends explicitly on the unknowns $\mathcal{W}(\diamond;\mathbf{x}),z_0(\mathbf{x}),\dots,z_{n-1}(\mathbf{x})$.  We may then view the rewritten \eqref{eq:partial-W} along with \eqref{eq:partial-a0}--\eqref{eq:partial-ak} as a first order differential equation on the Banach space $\mathcal{B}:=A_\infty(D)\times \mathbb{C}^n$ ($A_\infty(D)$ is the space of bounded analytic functions on $D$) with norm
\begin{equation}
    \|(\mathcal{W}(\diamond),z_0,\dots,z_{n-1})\|:=\sup_{z\in D}|\mathcal{W}(z)| +\sum_{k=0}^{n-1}|z_k|.
    \label{eq:norm}
\end{equation}
Given a solution $\nu$ of \eqref{eq:nu-power}, we already have an initial condition for this differential equation as the point $(\mathcal{W}(\diamond;\mathbf{x}_0),0,\dots,0)\in \mathcal{B}$.  It is not difficult to show that the right-hand side of the differential equation is Lipschitz on $\mathcal{B}$ and $C^\infty$ in $\mathbf{x}\in \mathcal{X}$ provided that $\mathcal{X}\subset\mathbb{R}^m$ is taken to be a sufficiently small neighborhood of $\mathbf{x}_0$ and $D$ is taken to be a sufficiently small neighborhood of $z_0\in\mathbb{C}$.  The argument exploits the maximum modulus principle to measure the first term in \eqref{eq:norm} on the boundary of $D$ away from the removable singularities in the interior.  Invoking the Picard theorem in this setting then provides existence and uniqueness of the solution.

    To properly define $\mathcal{W}(\diamond;\mathbf{x})$ and $z_0(\mathbf{x}),\dots,z_{n-1}(\mathbf{x})$ at a given point $\mathbf{x}\in \mathcal{X}$ given $\nu$, one should pick a specific smooth curve in $\mathcal{X}\subset\mathbb{R}^m$ connecting the initial point $\mathbf{x}_0$ and $\mathbf{x}$, and define the vector field $\partial$ to be tangent to the curve.  However, since the differential equations in $\mathcal{B}$ were obtained by differentiation of the assumed equation \eqref{eq:desired-identity} in the arbitrary direction $(d_1(\mathbf{x}),\dots,d_m(\mathbf{x}))$, it is irrelevant which curve one selects to define the solution, which proves that, given one of the $n+1$ values of $\nu$ from \eqref{eq:nu-power} that determines the initial condition, the quantities $\mathcal{W}(\diamond;\mathbf{x})$ and $z_0(\mathbf{x}),\dots,z_{n-1}(\mathbf{x})$ are well defined on $\mathcal{X}$.

        Furthermore, since $z\mapsto \mathcal{W}(z;\mathbf{x}_0)$ is univalent by construction, the flow ensures that $\mathcal{W}(z;\mathbf{x})$ remains univalent in $D$ for $\|\mathbf{x}-\mathbf{x}_0\|$ sufficiently small.  Hence a conformal mapping $z\mapsto \mathcal{W}(z;\mathbf{x})$ is indeed defined by the procedure.

    Once the solution is determined, one can differentiate \eqref{eq:partial-a0}--\eqref{eq:partial-ak} and \eqref{eq:partial-W} again and eliminate the first derivatives appearing on the right-hand side to obtain expressions for higher-order partial derivatives (possibly with respect to different vector fields $\partial$) that one can check are continuous for $\mathbf{x}$ near $\mathbf{x}_0$.  Hence the conformal mapping $\mathcal{W}(\diamond;\mathbf{x})$ and the coefficients $z_0(\mathbf{x}),\dots,z_{n-1}(\mathbf{x})$ are all infinitely differentiable functions on $\mathcal{X}\subset\mathbb{R}^m$.
\end{proof}

\subsubsection{The classical steepest descent method}
The classical steepest descent method mentioned in Step~\ref{step:expand} is the following result (see, e.g., \cite{AAA}).
\begin{lemma}[Steepest descent method] 
Suppose that $L$ is a piecewise-smooth contour along which $h(z)$ and $f(z)$ are analytic, that $\zeta\in L$ is a simple critical point of $h(z)$, and that $\mathrm{Im}(-\ii h(z))$ is constant along $L$ while $\mathrm{Re}(-\ii h(z))$ is maximized along $L$ at $z=\zeta$.  Then in the limit $\epsilon\downarrow 0$,
\begin{equation}
    \int_L \ee^{-\ii h(z)/\epsilon}f(z)\,\dd z = \ee^{\ii\varphi}\sqrt{\frac{2\pi\epsilon}{|h''(\zeta)|}}\ee^{-\ii h(\zeta)/\epsilon}\left(f(\zeta)+\mathcal{O}(\epsilon)\right),
    \label{eq:steepest-descent}
\end{equation}
provided the integral is convergent for $\epsilon>0$ sufficiently small, 
where $\varphi=\arg(\dd z)$ at $z=\zeta$ is the angle relative to the direction $\arg(z-\zeta)=0$ formed by the oriented tangent to $L$ at $z=\zeta$.  The same result holds for other contours $L'$ homotopic to $L$ within the domain of analyticity of the functions $h(z)$ and $f(z)$, and for $L$ it is only necessary that $\mathrm{Im}(-\ii h(z))$ be constant locally near $z=\zeta$.
\label{lem:steepest-descent}
\end{lemma}
We will need to apply this result in the setting that $h(z)=h(z;t,x)$ depends parametrically on $(t,x)\in\mathbb{R}^2$, and hence also the critical point and contour angle depend on these parameters:  $\zeta=\zeta(t,x)$ and $\varphi=\varphi(t,x)$.  In this situation, the leading term in \eqref{eq:steepest-descent} can be expanded near a fixed point $(t,x)=(t_0,x_0)$ as follows:
\begin{multline}
    \int_L\ee^{-\ii h(z;t,x)/\epsilon}f(z)\,\dd z = \ee^{\ii\varphi(t,x)}\sqrt{\frac{2\pi\epsilon}{|h''(\zeta(t,x);t,x)|}}\ee^{-\ii h(\zeta(t,x);t,x)/\epsilon}\\
    \times\left(f(\zeta(t_0,x_0)) + O(t-t_0) + O(x-x_0)+O(\epsilon)\right).
\label{eq:steepest-descent-reexpand}    
\end{multline}
Here we intentionally leave the factors on the first line unexpanded as they are independent of $f$ and hence common to all of the integrals appearing in any given row of the matrices $\mathbf{A}(t,x;\epsilon)$ and $\mathbf{B}(t,x;\epsilon)$ so will factor out of both determinants $N$ and $D$.

\subsubsection{Some Cauchy determinant identities and the pole-critical point relation}
\label{sec:Cauchy-Determinant}
Due to the matrix determinant structure of the BO solution formula \eqref{eq:IVP-solution}, we will need to evaluate two specific Cauchy-like matrix determinants. 
\begin{lemma} \label{lem:Cauchy-determinants}
    Let $p_1,\dots,p_N,w_0,\dots,w_N$ be $2N+1$ distinct complex numbers.  Then the following determinant identities hold:
    \begin{align}
        C^{(0)}(\mathbf{w},\mathbf{p}):=\begin{vmatrix}
             1 &  (w_0 - p_1)^{-1} & \cdots &  (w_0 - p_N)^{-1} \\
             1 &  (w_1 - p_1)^{-1} & \cdots &  (w_1 - p_N)^{-1} \\
            \vdots & \vdots & \ddots & \vdots \\
             1 &  (w_N-p_1)^{-1} & \cdots &  (w_N-p_N)^{-1}
        \end{vmatrix} &= \Delta(\mathbf{w},\mathbf{p}), \label{eq:CauchyDetFirstCoeff} \\
        C^{(1)}(\mathbf{w},\mathbf{p}):=\begin{vmatrix}
             w_0  &  (w_0 - p_1)^{-1} & \cdots &  (w_0 - p_N)^{-1} \\
             w_1 &  (w_1 - p_1)^{-1} & \cdots &  (w_1 - p_N)^{-1} \\
            \vdots & \vdots & \ddots & \vdots \\
             w_N &  (w_N-p_1)^{-1} & \cdots &  (w_N-p_N)^{-1}
        \end{vmatrix} &= \Delta(\mathbf{w},\mathbf{p}) \left( \sum_{i = 0}^N w_i - \sum_{n = 1}^N p_n \right),
        \label{eq:CauchyDetSecondCoeff}
    \end{align}
    where $\Delta(\mathbf{w},\mathbf{p})$ is defined by
    \begin{equation}
        \Delta(\mathbf{w},\mathbf{p}) := (-1)^N \frac{\displaystyle\prod_{i = 2}^N \prod_{j = 1}^{i-1} (w_i - w_j)(p_j - p_i)}{\displaystyle \prod_{i = 0}^N \prod_{j = 1}^N (w_i - p_j)} \prod_{i = 1}^N (w_i - w_0).
        \label{eq:BODeltaDef}
    \end{equation}
\end{lemma}

\begin{proof}
We introduce an auxiliary variable $p_0\in\mathbb{C}$ that is sufficiently large in modulus to not coincide with any of the $w_0,\dots,w_N$ or $p_1,\dots,p_N$, and define the $N+1\times N+1$
Cauchy determinant:
\begin{equation}
    C(\mathbf{w},\mathbf{p}, p_0) := \begin{vmatrix}
         (w_0 - p_0)^{-1} &  (w_0 - p_1)^{-1} & \cdots &  (w_0 - p_N)^{-1} \\
         (w_1 - p_0)^{-1} &  (w_1 - p_1)^{-1} & \cdots &  (w_1 - p_N)^{-1} \\
        \vdots & \vdots & \ddots & \vdots \\
         (w_N - p_0)^{-1} &  (w_N - p_1)^{-1} & \cdots &  (w_N - p_N)^{-1}
    \end{vmatrix}.
\end{equation}
We expand $C(\mathbf{w},\mathbf{p},p_0)$ in a Laurent series for large $p_0$ by expanding the elements in the first column in geometric series and using linearity of the determinant in the first column:
\begin{equation}
\begin{split}
    C(\mathbf{w},\mathbf{p}, p_0) &= -\frac{1}{p_0} \begin{vmatrix}
        (1 - w_0/p_0)^{-1} & (w_0 - p_1)^{-1} & \cdots & (w_0 - p_N)^{-1} \\
        (1 - w_1/p_0)^{-1} & (w_1 - p_1)^{-1} & \cdots & (w_1 - p_N)^{-1} \\
        \vdots & \vdots & \ddots & \vdots \\
         (1 - w_N/p_0)^{-1} &  (w_N-p_1)^{-1} & \cdots &  (w_N-p_N)^{-1} \\
    \end{vmatrix} \\
    &= -\frac{1}{p_0} \begin{vmatrix}
         1 + w_0p_0^{-1} + O\left( p_0^{-2} \right)  &  (w_0 - p_1)^{-1} & \cdots &  (w_0 - p_N)^{-1} \\
         1 + w_1p_0^{-1} + O\left( p_0^{-2} \right) &  (w_1 - p_1)^{-1} & \cdots &  (w_1 - p_N)^{-1} \\
        \vdots & \vdots & \ddots & \vdots \\
         1 +w_Np_0^{-1} + O\left( p_0^{-2} \right) &  (w_N-p_1)^{-1} & \cdots &  (w_N-p_N)^{-1} \\
    \end{vmatrix} \\
    &= -\frac{1}{p_0} C^{(0)}(\mathbf{w},\mathbf{p}) - \frac{1}{p_0^2} C^{(1)}(\mathbf{w},\mathbf{p}) + O\left( \frac{1}{p_0^3} \right).
\end{split}
\end{equation}
On the other hand, the value of $C(\mathbf{w},\mathbf{p}, p_0)$ can be computed exactly from the known Cauchy determinant formula:
\begin{equation}
\begin{split}
    C(\mathbf{w},\mathbf{p}, p_0) &= \frac{\displaystyle\prod_{i = 1}^N \prod_{j = 0}^{i-1} (w_i - w_j)(p_j - p_i)}{\displaystyle\prod_{i = 0}^N \prod_{j = 0}^N (w_i - p_j)} \\
    &= \frac{\displaystyle\prod_{i = 2}^N \prod_{j = 1}^{i-1} (w_i - w_j)(p_j - p_i)}{\displaystyle\prod_{i = 0}^N \prod_{j = 1}^N (w_i - p_j)} \frac{\displaystyle\prod_{i = 1}^N (w_i - w_0)(p_0 - p_i)}{\displaystyle\prod_{i = 0}^N (w_i - p_0)} \\
    &= -\Delta(\mathbf{w},\mathbf{p}) \frac{1}{p_0} \frac{\displaystyle\prod_{n = 1}^N (1 - p_n/p_0)}{\displaystyle\prod_{i = 0}^N (1 - w_i/p_0)}  \\
    &= -\Delta(\mathbf{w},\mathbf{p}) \left( \frac{1}{p_0} + \left( \sum_{i = 0}^N w_i - \sum_{n = 1}^N p_n \right) \frac{1}{p_0^2} + \bigO{\frac{1}{p_0^3}} \right),
\end{split}
\end{equation}
where the $p_0$-independent factors have been gathered into the expression $\Delta(\mathbf{w},\mathbf{p})$ defined in \eqref{eq:BODeltaDef}.
Comparing the unique coefficients of the two expansions, we complete the proofs of \eqref{eq:CauchyDetFirstCoeff} and \eqref{eq:CauchyDetSecondCoeff}.  
\end{proof}

\begin{remark}
    An alternate formulation of this result can be found in \cite[Lemma 8]{Gerard23}, with a different proof.  We note that the above proof also allows for the explicit computation of determinants $C^{(n)}(\mathbf{w},\mathbf{p})$ with higher powers of $w_j$ in the first column.
\end{remark}

Additionally, we define the \emph{pole-critical point relation} via the following lemma.

\begin{lemma} \label{lem:PoleSPRelation}
    Fix $(t,x)\in\mathbb{R}^2$.  Let $\mathcal{Y}$ denote the multiset of real critical points of $h(\diamond;t,x)$ and let $\mathcal{Z}$ denote the multiset of complex critical points of $h(\diamond;t,x)$ that lie in the open upper half-plane.  In both cases, we include non-simple critical points with repetition according to their multiplicity. Then
    \begin{equation}
        x + \sum_{n = 1}^N (p_n + p_n^*) = \sum_{y_\alpha\in\mathcal{Y}} y_\alpha + \sum_{z_\beta\in\mathcal{Z}} (z_\beta + z_\beta^*),
        \label{eq:PoleSPRelation}
    \end{equation}
    that is, $x$ plus the sum of positions of the poles of the rational initial data $u_0$ is always equal to the sum of all the positions of the critical points. 
\end{lemma}

The proof of this follows from clearing fractions in the rational function $h'(\diamond; t, x)$, noting the roots of the numerator polynomial are exactly the critical points, and then equating coefficients after expressing the numerator polynomial in two different ways. The proof is identical to that of \cite[Lemma 6.3]{SmallDispersionBO-1} after making the appropriate adjustments to account for higher order critical points.

\begin{remark}
    The identity \eqref{eq:PoleSPRelation} was also derived by G\'erard \cite[Proof of Lemma 8]{Gerard23} in a different way.
\end{remark}

\section{Analysis near the soliton edge $x=x_\srm(t)$.  Proof of Theorem~\ref{thm:soliton-edge}}
\label{sec:SolitonEdgeProof}
The arguments in this section are adapted from \cite[Section 4.2]{MitchellThesis}.  In the following we give the proof of Theorem~\ref{thm:soliton-edge}, following the three steps indicated in Section~\ref{sec:structure-of-proofs}.
\subsection{Step~\ref{step:normal-form}:  Conformal map and normal form of the exponent}
\label{sec:soliton-edge-step1}
To analyze the contour integrals in the second rows of $\mathbf{A}(t,x;\epsilon)$ and $\mathbf{B}(t,x;\epsilon)$ involving the contour $W_1$ illustrated in Figure~\ref{fig:soliton-edge-contours}, we first apply Lemma~\ref{lem:normal-form} with $m=1$ and $n=2$ to the function $H(z;\mathbf{x}):=h(z;t,x)$ for fixed $t\in K=(t_\cc,t_\cc+\delta)$ with $z_0=y_\srm(t)$ and $x_0=x_\srm(t)$.  Hence $\mathcal{X}$ is an open interval containing $x_0$.  According to the discussion in Section~\ref{sec:soliton-edge-contours}, the parameter $\nu=\nu_\srm(t)\neq 0$ in \eqref{eq:nu-power} then satisfies 
\begin{equation}
    \nu_\srm(t)^3 = \frac{1}{2}h'''(y_\srm(t);t,x_\srm(t))=\frac{1}{2}u_0''(y_\srm(t))<0.
    \label{eq:nu-s-cubed}
\end{equation}
We choose the negative real cube root $\nu_\srm(t)<0$ so that the corresponding conformal map denoted $w=\mathcal{W}_\srm(z;t,x)$ maps the real line near $y_\srm(t)$ to a real neighborhood of $w=0$ in a monotone-decreasing fashion.  Then, on a neighborhood $D$ of $y_\srm(t)$, Lemma~\ref{lem:normal-form} asserts that under the conformal mapping $w=\mathcal{W}_\srm(z;t,x)$,
\begin{equation}
    h(z;t,x)=h(y_\srm(t);t,x)+\frac{1}{3}w^3 +q_{\srm}(t,x)w+r_{\srm}(t,x),\quad w=\mathcal{W}_\srm(z;t,x),
\label{eq:cubic-exponent}
\end{equation}
where $q_{\srm}$ and $r_{\srm}$ are $C^\infty$ functions of $x\in \mathcal{X}$ that satisfy $q_{\srm}(t,x_\srm(t))=r_\srm(t,x_\srm(t))=0$.  It is easy to see that they are real-valued for $t\in K$ and $x\in \mathcal{X}$.  Since for fixed $z$, 
\begin{equation}
\partial_x h(z;t,x)= -\frac{z-x}{2t},    
\end{equation}
using \eqref{eq:partial-a0-crit} and \eqref{eq:partial-ak-crit} respectively, we
have 
\begin{equation}
    \partial_xr_{\srm}(t,x_\srm(t))=-\frac{1}{2t}\frac{1}{2\pi\ii}\oint_L\frac{(z-x_\srm(t))\mathcal{W}_\srm'(z;t,x_\srm(t))}{\mathcal{W}_\srm(z;t,x_\srm(t))}\,\dd z +\frac{y_\srm(t)-x_\srm(t)}{2t}
\end{equation}
and
\begin{equation}
    \partial_x q_{\srm}(t,x_\srm(t))=-\frac{1}{2t}\frac{1}{2\pi\ii}\oint_L\frac{(z-x_\srm(t))\mathcal{W}_\srm'(z;t,x_\srm(t))}{\mathcal{W}_\srm(z;t,x_\srm(t))^2}\,\dd z.
\end{equation}
Using the Taylor expansion $\mathcal{W}_\srm(z;t,x_\srm(t))=\nu_\srm(t)(z-y_\srm(t))(1+O(z-y_\srm(t)))$ near $z=y_\srm(t)$ then gives, by a residue calculation at $z=y_\srm(t)$,
\begin{equation}
    \partial_xr_{\srm}(t,x_\srm(t))=0\quad\text{and}\quad 
    \partial_x q_{\srm}(t,x_\srm(t))=-\frac{1}{2t\nu_\srm(t)}.
\label{eq:a0-a1-derivs-c}
\end{equation}

\subsection{Step~\ref{step:expand}:  Expansions of integrals}
\label{sec:soliton-edge-step2}
Following Chester-Friedman-Ursell \cite{ChesterFriedmanUrsell}, we use the normal form to evaluate the contour integrals in the second row of the matrices $\mathbf{A}(t,x;\epsilon)$ and $\mathbf{B}(t,x;\epsilon)$, which have the form
\begin{equation}
    I(t,x;\epsilon):=\int_{W_1}\ee^{-\ii h(z;t,x)/\epsilon}f(z)\,\dd z 
\end{equation}
in which either $f(z)=u_0(z)$, $f(z)\equiv 1$, or $f(z)=(z-p_n)^{-1}$ for $n=1,\dots,N$.  These functions are all analytic and independent of $(t,x)$ along the contour $W_1$.  Factoring out the value of the exponential at $z=y_\srm(t)$ and multiplying by the oscillatory exponential $\ee^{\ii r_{\srm}(t,x)/\epsilon}$ gives
\begin{equation}
    \ee^{\ii [h(y_\srm(t);t,x)+r_{\srm}(t,x)]/\epsilon}I(t,x;\epsilon)=\int_{W_1}\ee^{-\ii [h(z;t,x)-h(y_\srm(t);t,x)-r_{\srm}(t,x)]/\epsilon}f(z)\,\dd z.
\end{equation}
Contributions to the integral coming from arcs of $W_1$ outside of a small neighborhood $D$ of $z=y_\srm(t)$ are exponentially small as $\epsilon\downarrow 0$ by choice of the contour $W_1$, and along the arc of $W_1$ within $D$ we may use the conformal map $w=\mathcal{W}_\srm(z;t,x)$ and \eqref{eq:cubic-exponent} to obtain
\begin{multline}
    \ee^{\ii [h(y_\srm(t);t,x)+r_{\srm}(t,x)]/\epsilon}I(t,x;\epsilon)\\
    =\int_{\mathcal{W}_\srm(W_1\cap D;t,x)}\ee^{-\ii[w^3/3 +q_{\srm}(t,x)w]/\epsilon}f(\mathcal{W}_\srm^{-1}(w;t,x))\frac{\dd \mathcal{W}^{-1}_\srm}{\dd w}(w;t,x)\,\dd w + O(\epsilon^\infty),\quad\epsilon\downarrow 0.
    \label{eq:soliton-I-conformal}
\end{multline}
Since $\nu_\srm(t)<0$, when $|x-x_\srm(t)|$ is sufficiently small, the arc $\mathcal{W}_\srm(W_1\cap D;t,x)$ enters a neighborhood of the origin $w=0$ from the direction $\arg(w)=-\frac{1}{6}\pi$ and exits in the direction $\arg(w)=-\frac{5}{6}\pi$.    By Taylor approximation of $f(\mathcal{W}_\srm^{-1}(w;t,x))\mathcal{W}_\srm^{-1\prime}(w;t,x)$ about $w=0$ with remainder estimate:
\begin{equation}
    f(\mathcal{W}_\srm^{-1}(w;t,x))\frac{\dd \mathcal{W}^{-1}_\srm}{\dd w}(w;t,x)=\sum_{m=0}^M\Tc^{(m)}_\srm[f](t,x)w^m + O(w^{M+1}),
    \label{eq:Taylor-of-f}
\end{equation}
we have
\begin{multline}
    \ee^{\ii [h(y_\srm(t);t,x)+r_{\srm}(t,x)]/\epsilon}I(t,x;\epsilon)\\
    =\sum_{m=0}^M\Tc^{(m)}_\srm[f](t,x)\int_{\mathcal{W}_\srm(W_1\cap D;t,x)}\ee^{-\ii[w^3/3+q_{\srm}(t,x)w]/\epsilon}w^m\,\dd w \\
    + O\left(\int_{\mathcal{W}_\srm(W_1\cap D;t,x)}\left|\ee^{-\ii [w^3/3+q_{\srm}(t,x)w]/\epsilon}\right||w|^{M+1}\,|\dd w|\right) + O(\epsilon^\infty).
\end{multline}
Rescaling by $w=\ii \epsilon^{1/3}\zeta$ gives
\begin{multline}
    \ee^{\ii [h(y_\srm(t);t,x)+r_{\srm}(t,x)]/\epsilon}I(t,x;\epsilon)\\
    {}=\sum_{m=0}^M\ii^{m+1}\epsilon^{(m+1)/3}\Tc^{(m)}_\srm[f](t,x)\int_{-\ii\epsilon^{-1/3}\mathcal{W}_\srm(W_1\cap D;t,x)}\ee^{Y\zeta-\zeta^3/3}\zeta^m\,\dd \zeta \\{}+ O\left(\epsilon^{(M+2)/3}\int_{-\ii\epsilon^{-1/3}\mathcal{W}_\srm(W_1\cap D;t,x)}\left|\ee^{Y\zeta-\zeta^3/3}\right||\zeta|^{M+1}\,|\dd \zeta|\right) + O(\epsilon^\infty),
\end{multline}
where
\begin{equation}
    Y:=\epsilon^{-2/3}q_{\srm}(t,x).
\end{equation}
Assuming that $Y$ remains bounded as $\epsilon\downarrow 0$, we may extend the contour of integration to an infinite contour from $\infty\ee^{-2\pi\ii/3}$ to $\infty\ee^{2\pi\ii/3}$ at the cost of an error small beyond all orders. Thus we obtain
\begin{multline}
        \ee^{\ii [h(y_\srm(t);t,x)+r_{\srm}(t,x)]/\epsilon}I(t,x;\epsilon)\\
    {}=\sum_{m=0}^M\ii^{m+1}\epsilon^{(m+1)/3}\Tc^{(m)}_\srm[f](t,x)\int_{\infty\ee^{-2\pi\ii/3}}^{\infty\ee^{2\pi\ii/3}}\ee^{Y\zeta-\zeta^3/3}\zeta^m\,\dd\zeta + O(\epsilon^{(M+2)/3}),
\end{multline}
or recalling the definition of $\Ai(\diamond)$ as a contour integral, we may finally write for any $M=0,1,2,\dots$ that
\begin{equation}
    \ee^{\ii [h(y_\srm(t);t,x)+r_{\srm}(t,x)]/\epsilon}I(t,x;\epsilon)=2\pi\ii\sum_{m=0}^M\ii^{m+1}\epsilon^{(m+1)/3}\Tc^{(m)}_\srm[f](t,x)\Ai^{(m)}(Y) + O(\epsilon^{(M+2)/3})
\end{equation}
as $\epsilon\downarrow 0$. This expansion only involves the function $\Ai(Y)$ and its derivative $\Ai'(Y)$ due to the differential equation $\Ai''(Y)=Y\Ai(Y)$.  Now if $Y$ is bounded, then we also have from $q_{\srm}(t,x_\srm(t))=0$ and \eqref{eq:a0-a1-derivs-c} that
\begin{equation}
    Y=X_\srm + O(\epsilon^{2/3}),
\end{equation}
where $X_\srm$ is as defined in \eqref{eq:Xs-coordinate}.  Therefore, keeping just the first two terms ($M=1$) we have
\begin{multline}
    \ee^{\ii [h(y_\srm(t);t,x)+r_{\srm}(t,x)]/\epsilon}I(t,x;\epsilon)\\=-2\pi\epsilon^{1/3}\Tc^{(0)}_\srm[f](t,x)\Ai(X_\srm) -2\pi\ii\epsilon^{2/3}\Tc^{(1)}_\srm[f](t,x)\Ai'(X_\srm) + O(\epsilon),\quad\epsilon\downarrow 0.
\end{multline}
Furthermore, since $x=x_\srm(t)+O(\epsilon^{2/3})$ by assumption, the differentiability of the Taylor coefficients $\Tc^{(j)}_\srm[f](t,x)$ implies that $\Tc^{(j)}_\srm[f](t,x)=\Tc^{(j)}_\srm[f](t,x_\srm(t))+O(\epsilon^{2/3})$, and so for an error term of the same order,
\begin{multline}
    \ee^{\ii [h(y_\srm(t);t,x)+r_{\srm}(t,x)]/\epsilon}I(t,x;\epsilon)\\=-2\pi\epsilon^{1/3}\Tc^{(0)}_\srm[f](t,x_\srm(t))\Ai(X_\srm) -2\pi\ii\epsilon^{2/3}\Tc^{(1)}_\srm[f](t,x_\srm(t))\Ai'(X_\srm) + O(\epsilon),\quad\epsilon\downarrow 0.
\label{eq:W1-s-expand}
\end{multline}
It remains to compute the Taylor coefficients $\Tc^{(0)}_\srm[f](t,x_\srm(t))$ and $\Tc^{(1)}_\srm[f](t,x_\srm(t))$.  For this we start with the Taylor expansion of $\mathcal{W}_\srm(z;t,x_\srm(t))$ (see \eqref{eq:W-general-crit}):
\begin{equation}
w=    \mathcal{W}_\srm(z;t,x_\srm(t))=\nu_\srm(t)(z-y_\srm(t)) + c_\srm(t)(z-y_\srm(t))^2 + O((z-y_\srm(t))^2),
\end{equation}
where
\begin{equation}
    c_\srm(t):=\frac{1}{24\nu_\srm(t)^2}h^{(4)}(y_\srm(t);t,x_\srm(t)),
    \label{eq:c-s}
\end{equation}
from which we easily get the expansion of the inverse mapping
\begin{equation}
    z=\mathcal{W}_\srm^{-1}(w;t,x_\srm(t))=y_\srm(t) + \frac{1}{\nu_\srm(t)}w-\frac{c_\srm(t)}{\nu_\srm(t)^3}w^2+O(w^3).
\end{equation}
Therefore using this result in \eqref{eq:Taylor-of-f} along with the Taylor expansion of $f(\diamond)$ we obtain
\begin{equation}
\begin{split}
    \Tc^{(0)}_\srm[f](t,x_\srm(t))&=\frac{1}{\nu_\srm(t)}f(y_\srm(t)),\\
    \Tc^{(1)}_\srm[f](t,x_\srm(t))&=\frac{1}{\nu_\srm(t)^2}\left[f'(y_\srm(t))-\frac{2c_\srm(t)}{\nu_\srm(t)}f(y_\srm(t))\right].
    \end{split}
\label{eq:ks0ks1}
\end{equation}
It follows that \eqref{eq:W1-s-expand} can be written in the form
\begin{multline}
    -\frac{\nu_\srm(t)}{2\pi\epsilon^{1/3}}\ee^{\ii[h(y_\srm(t);t,x)+r_\srm(t,x)]/\epsilon}I(t,x;\epsilon)=f(y_\srm(t))\Ai(X_\srm)\\ + \frac{\ii\epsilon^{1/3}}{\nu_\srm(t)}\left[f'(y_\srm(t))-\frac{2c_\srm(t)}{\nu_\srm(t)}f(y_\srm(t))\right]\Ai'(X_\srm) + O(\epsilon^{2/3}),\quad\epsilon\downarrow 0.
    \label{eq:W1-s-expand-rewrite}
\end{multline}

The coefficient $2c_\srm(t)/\nu_\srm(t)$ can be computed explicitly.  Indeed, note that by combining \eqref{eq:nu-s-cubed} and \eqref{eq:c-s} we obtain
\begin{equation}
    \frac{2c_\srm(t)}{\nu_\srm(t)}=\frac{h^{(4)}(y_\srm(t);t,x_\srm(t))}{6h'''(y_\srm(t);t,x_\srm(t))} \, .
    \label{eq:TwoCOverNu}
\end{equation}

Now, according to \eqref{eq:h-def} and the discussion in Section~\ref{sec:soliton-edge-contours}, $z\mapsto h'(z;t,x_\srm(t))$ is a rational function having  simple poles only at $p_1,\dots,p_N$ and $p_1^*,\dots,p_N^*$, simple zeros only at $z=y_0(t)$ and the complex values $z=z_j(t,x_\srm(t))$ and $z=z_j(t,x_\srm(t))^*$ for $j=2,\dots,N$, and a double zero at $z=y_\srm(t)$, with no other zeros or poles. Therefore, since $h'(z;t,x_\srm(t))=z/(2t)+O(1)$ as $z\to\infty$, we have the representation
\begin{align}
    h'(z;t,x_\srm(t)) &= (z-y_\srm(t))^2 G(z;t), \\
    G(z;t) &:= \frac{1}{2t}(z-y_0(t))\frac{\displaystyle\prod_{j=2}^N (z-z_j(t,x_\srm(t)))(z-z_j(t,x_\srm(t))^*)}{\displaystyle \prod_{n=1}^N (z-p_n)(z-p_n^*)}.
\end{align}
Taking either two or three more derivatives with respect to $z$ using this representation and evaluating at the double zero $z=y_\srm(t)$ gives 
\begin{equation}
    h'''(y_\srm(t);t,x_\srm(t))=2G(y_\srm(t);t)\quad\text{and}\quad
    h^{(4)}(y_\srm(t);t,x_\srm(t))=6G'(y_\srm(t);t)
\end{equation}
because only the terms in which exactly two derivatives fall on the factor $(z-y_\srm(t))^2$ produce a nonzero result upon evaluation at $z=y_\srm(t)$.  
It follows that
\begin{equation}
\begin{split}
    \frac{2c_\srm(t)}{\nu_\srm(t)}&=\frac{1}{2}\frac{G'(y_\srm(t);t)}{G(y_\srm(t);t)}\\ & = \frac{1}{2}\cdot\frac{1}{y_\srm(t)-y_0(t)} + 
    \mathrm{Re}\left(\sum_{j=2}^N\frac{1}{y_\srm(t)-z_j(t,x_\srm(t))}-\sum_{n=1}^N\frac{1}{y_\srm(t)-p_n}\right),
    \end{split}
\label{eq:2covernu-s}
\end{equation}
where we used the fact that $y_\srm(t)\in\mathbb{R}$.

All of the contour integrals in all rows of $\mathbf{A}(t,x;\epsilon)$ and $\mathbf{B}(t,x;\epsilon)$ other than the second row take their dominant contributions from simple critical points of $h(z;t,x)$ and can be approximated using \eqref{eq:steepest-descent-reexpand}.  

\subsection{Step~\ref{step:algebra}:  Simplification of determinants}
\label{sec:soliton-edge-step3}

We will focus first on an expansion of $N(t, x; \epsilon)=\det(\mathbf{A}(t,x;\epsilon))$. Applying the colliding double critical point approximation \eqref{eq:W1-s-expand-rewrite} to the second row of the matrix $\mathbf{A}(t,x;\epsilon)$ (integrating over the contour $W_1$) and the simple critical point approximation \eqref{eq:steepest-descent-reexpand} to the other $N$ rows, we can identify an overall nonzero factor of $N(t,x;\epsilon)$ comprised of the common factors in the expansion of matrix elements from each row. This overall factor is:  
\begin{equation}
    \begin{split}
        E_\srm(t, x; \epsilon) &:= \ee^{-\ii \pi/4} \sqrt{ \frac{2 \pi \epsilon}{|h''(y_0(t,x); t, x)|} } \ee^{-\ii h(y_0(t,x); t, x)/\epsilon} \left(-\frac{2 \pi \epsilon^{1/3}}{\nu_\srm(t)}\right) \ee^{- \ii [h(y_\srm(t); t, x)+r_\srm(t,x)] / \epsilon} \\
        &\Quad[2] \times \prod_{j = 2}^N \ee^{\ii \varphi_j(t, x)} \sqrt{ \frac{2 \pi \epsilon}{|h''(z_j(t, x); t, x)|} } \ee^{-\ii h(z_j(t, x); t, x)/\epsilon},
    \end{split}
    \label{eq:SolExpFactor}
\end{equation}
where $\varphi_j(t,x)$ denotes the angle of steepest-descent passage of the integration contour $W_j$ over the critical point $z=z_j(t,x)$.
We use the linearity of the determinant in the rows to split the leading order contributions from the error terms iteratively in the last $N$ rows as well as the leading order, next-to-leading order and error terms in the first row to get the expansion
\begin{align}
    -\frac{N(t, x; \epsilon)}{E_\srm(t, x; \epsilon)} &= \Ai(X_\srm) N_\srm(t) + \ii \frac{\epsilon^{1/3}}{\nu_\srm(t)} \Ai'(X_\srm) \left( N_\srm^{(1)}(t) - 2\frac{c_\srm(t)}{\nu_\srm(t)} N_\srm(t) \right) + \bigO{\epsilon^{2/3}},
    \label{eq:SolNumFirstExpans}
\end{align}
where $N_\srm(t)$ is the new matrix determinant
\begin{equation}
    N_\srm(t) := \begin{vmatrix}
        u_0\big( y_\srm(t) \big) &  (y_\srm(t) - p_1)^{-1} & \cdots &  (y_\srm(t) - p_N)^{-1} \\
         u_0\big( y_0(t) \big) &  (y_0(t) - p_1)^{-1} & \cdots &  (y_0(t) - p_N)^{-1} \\
         u_0\big( z_{2,\srm}(t) \big) &  (z_{2,\srm}(t) - p_1)^{-1} & \cdots &  (z_{2,\srm}(t) - p_N)^{-1} \\
        \vdots & \vdots & \ddots & \vdots \\
         u_0\big( z_{N,\srm}(t) \big) &  (z_{N,\srm}(t) - p_1)^{-1} & \cdots &  (z_{N,\srm}(t) - p_N)^{-1}
    \end{vmatrix},
    \label{eq:USDef}
\end{equation}
while $N_\srm^{(1)}(t)$ is the same determinant with a derivative taken on the functions in the first row:
\begin{equation}
    N_\srm^{(1)}(t) := \begin{vmatrix}
         u_0'\big( y_\srm(t) \big) &  -(y_\srm(t) - p_1)^{-2} & \cdots &  -(y_\srm(t) - p_N)^{-2} \\
         u_0\big( y_0(t) \big) &  (y_0(t) - p_1)^{-1} & \cdots &  (y_0(t) - p_N)^{-1} \\
         u_0\big( z_{2,\srm}(t) \big) &  (z_{2,\srm}(t) - p_1)^{-1} & \cdots &  (z_{2,\srm}(t) - p_N)^{-1} \\
        \vdots & \vdots & \ddots & \vdots \\
         u_0\big( z_{N,\srm}(t) \big) &  (z_{N,\srm}(t) - p_1)^{-1} & \cdots &  (z_{N,\srm}(t) - p_N)^{-1}
    \end{vmatrix}.
    \label{eq:US1Def}
\end{equation}
The minus sign on the left-hand side of \eqref{eq:SolNumFirstExpans} is present because for convenience we have swapped the first two rows of the determinants on the right-hand side.
Here, we have introduced the shorthand notation $z_{j,\srm}(t) := z_j(t, x_\srm(t))$ for brevity of formulas. Since $y_\srm(t)$, $y_0(t)$, $z_{2,\srm}(t)$, \dots, $z_{N,\srm}(t)$ are all critical points of $h$, they satisfy 
\begin{align}
    0 &= \frac{y_\srm(t) - x_\srm(t)}{2t} + u_0(y_\srm(t)) \, , \label{eq:SolSSaddleCond} \\
    0 &= \frac{y_0(t) - x_\srm(t)}{2t} + u_0(y_0(t)) \, , \label{eq:Sol0SaddleCond} \\
    0 &= \frac{z_{j,\srm}(t) - x_\srm(t)}{2t} + u_0(z_{j,\srm}(t)), \quad j=2,\ldots,N. \label{eq:SolZSaddleCond}
\end{align}
This means we can replace the $u_0$ evaluations in \eqref{eq:USDef} and \eqref{eq:US1Def}. Moreover, since $y_\srm(t)$ is a double critical point, we have the derivative condition
\begin{equation}
    0 = \frac{1}{2 t} + u_0'(y_\srm(t)),
\end{equation}
and thus $u_0'(y_\srm(t)) = -1/2t$ in \eqref{eq:US1Def}. Therefore, using linearity to factor out $1/(2t)$ from the first column of each determinant yields
\begin{align}
    N_\srm(t) &= \frac{1}{2t}\begin{vmatrix}
        x_\srm(t) - y_\srm(t) &  (y_\srm(t) - p_1)^{-1} & \cdots &  (y_\srm(t) - p_N)^{-1} \\
        x_\srm(t) - y_0(t) &  (y_0(t) - p_1)^{-1} & \cdots &  (y_0(t) - p_N)^{-1} \\
        x_\srm(t) - z_{2,\srm}(t) &  (z_{2,\srm}(t) - p_1)^{-1} & \cdots &  (z_{2,\srm}(t) - p_N)^{-1} \\
        \vdots & \vdots & \ddots & \vdots \\
        x_\srm(t) - z_{N,\srm}(t) &  (z_{N,\srm}(t) - p_1)^{-1} & \cdots &  (z_{N,\srm}(t) - p_N)^{-1}
    \end{vmatrix} \, , \label{eq:USFullDet} \\
    N^{(1)}_\srm(t) &= \frac{1}{2t}\begin{vmatrix}
         -1 &  -(y_\srm(t) - p_1)^{-2} & \cdots &  -(y_\srm(t) - p_N)^{-2} \\
        x_\srm(t) - y_0(t) &  (y_0(t) - p_1)^{-1} & \cdots &  (y_0(t) - p_N)^{-1} \\
         x_\srm(t) - z_{2,\srm}(t) &  (z_{2,\srm}(t) - p_1)^{-1} & \cdots &  (z_{2,\srm}(t) - p_N)^{-1} \\
        \vdots & \vdots & \ddots & \vdots \\
         x_\srm(t) - z_{N,\srm}(t) &  (z_{N,\srm}(t) - p_1)^{-1} & \cdots &  (z_{N,\srm}(t) - p_N)^{-1}
    \end{vmatrix}.
    \label{eq:US1FullDet}
\end{align}
Further expanding via linearity in the first column and recalling the notation in \eqref{eq:CauchyDetFirstCoeff} and \eqref{eq:CauchyDetSecondCoeff} with vectors 
\begin{align}
    \mathbf{w}_\srm(t) &:= (y_\srm(t), y_0(t), z_{2,\srm}(t), \dots, z_{N,\srm}(t))^\top, \\
    \mathbf{p} &:= (p_1, p_2, \dots, p_N)^\top ,
    \label{eq:p-vector}
\end{align}
we have
\begin{align}
    N_\srm(t) &= \frac{x_\srm(t)}{2t} C^{(0)}(\mathbf{w}_\srm(t), \mathbf{p}) - \frac{1}{2t} C^{(1)}(\mathbf{w}_\srm(t), \mathbf{p}), \label{eq:US} \\
    N^{(1)}_\srm(t) &= \frac{x_\srm(t)}{2t} \frac{\partial C^{(0)}}{\partial w_0}(\mathbf{w}_\srm(t), \mathbf{p}) - \frac{1}{2t} \frac{\partial C^{(1)}}{\partial w_0}(\mathbf{w}_\srm(t), \mathbf{p}).
    \label{eq:US1}
\end{align}
Applying Lemma \ref{lem:Cauchy-determinants} then yields
\begin{align}
    N_\srm(t) &= \frac{\Delta(\mathbf{w}_\srm(t), \mathbf{p})}{2t} \left( x_\srm(t) - y_\srm(t) - y_0(t) - \sum_{j = 2}^N z_{j,\srm}(t) + \sum_{k = 1}^N p_k \right) \, , \label{eq:USDelta} \\
    N^{(1)}_\srm(t) &= \frac{1}{2t}\frac{\partial \Delta}{\partial w_0}(\mathbf{w}_\srm(t), \mathbf{p}) \left( x_\srm(t) - y_\srm(t) - y_0(t) - \sum_{j = 2}^N z_{j,\srm}(t) + \sum_{k = 1}^N p_k \right) - \frac{\Delta(\mathbf{w}_\srm(t), \mathbf{p})}{2t} \, ,
    \label{eq:US1Delta}
\end{align}
where we remind the reader that $\Delta(\mathbf{w},\mathbf{p})$ is given explicitly by a product formula in \eqref{eq:BODeltaDef}.
Finally, substituting \eqref{eq:USDelta} and \eqref{eq:US1Delta} into \eqref{eq:SolNumFirstExpans} and noting that $\Delta(\mathbf{w}_\srm(t),\mathbf{p})\neq 0$, we get
\begin{equation}
    \begin{split}
        -\frac{N(t, x; \epsilon)}{E_\srm(t, x; \epsilon)} &= \frac{1}{2 t} \left( x_\srm(t) - y_\srm(t) - y_0(t) - \sum_{j = 2}^N z_{j,\srm}(t) + \sum_{k = 1}^N p_k \right) \Delta(\mathbf{w}_\srm(t), \mathbf{p}) \\
        &\Quad[4] \times \left[ \Ai(X_\srm) + \ii \frac{\epsilon^{1/3}}{\nu_\srm(t)} \Ai'(X_\srm) \left( \frac{1}{\Delta(\mathbf{w}_\srm(t), \mathbf{p})}\frac{\partial \Delta}{\partial w_0}(\mathbf{w}_\srm(t), \mathbf{p}) - 2\frac{c_\srm(t)}{\nu_\srm(t)} \right) \right] \\
        &\Quad[2] - \ii \frac{\epsilon^{1/3}}{2 t \nu_\srm(t)} \Ai'(X_\srm) \Delta(\mathbf{w}_\srm(t), \mathbf{p}) + \bigO{\epsilon^{2/3}}.
    \end{split}
    \label{eq:SolNumFinalExpan}
\end{equation}

This is as far as we can go with the numerator, so now we work on the denominator matrix determinant $D(t,x;\epsilon)=\det(\mathbf{B}(t,x;\epsilon))$. We expand $D(t,x;\epsilon)$ in a similar fashion.  The first column of determinants analogous to \eqref{eq:USDef} and \eqref{eq:US1Def} is simply replaced by a column of all $1$, so the determinants are directly expressible in terms of $C^{(0)}(\mathbf{w}, \mathbf{p})$ and its derivative with respect to $w_0$ evaluated at $\mathbf{w}=\mathbf{w}_\srm(t)$ defined in \eqref{eq:CauchyDetFirstCoeff}:
\begin{equation}
    \begin{split}
        -\frac{D(t, x; \epsilon)}{E_\srm(t, x; \epsilon)} &= \Ai(X_\srm) C^{(0)}(\mathbf{w}_\srm(t), \mathbf{p}) \\
        &\Quad[2] + \ii \frac{\epsilon^{1/3}}{\nu_\srm(t)} \Ai'(X_\srm) \left( \frac{\partial C^{(0)}}{\partial w_0}(\mathbf{w}_\srm(t), \mathbf{p}) - 2\frac{c_\srm(t)}{\nu_\srm(t)} C^{(0)}(\mathbf{w}_\srm(t), \mathbf{p}) \right) \\
        &\Quad[2] + \bigO{\epsilon^{2/3}} \, .
    \end{split}
    \label{eq:SolDenomFirstExpans}
\end{equation}
Applying Lemma~\ref{lem:Cauchy-determinants} then yields
\begin{multline}
    -\frac{D(t, x; \epsilon)}{E_\srm(t, x; \epsilon)} = \Delta(\mathbf{w}_\srm(t), \mathbf{p}) \left[ \vphantom{\frac{\partial \Delta}{\partial w_0}}\Ai(X_\srm) \right.\\
    \left.+\ii \frac{\epsilon^{1/3}}{\nu_\srm(t)} \Ai'(X_\srm) \left( \frac{1}{\Delta(\mathbf{w}_\srm(t), \mathbf{p})}\frac{\partial \Delta}{\partial w_0}(\mathbf{w}_\srm(t), \mathbf{p}) - 2\frac{c_\srm(t)}{\nu_\srm(t)} \right) \right]  
     + \bigO{\epsilon^{2/3}}.
     \label{eq:SolDenomFinalExpans}
\end{multline}
Using \eqref{eq:SolDenomFinalExpans} shows that \eqref{eq:SolNumFinalExpan} can be written equivalently in the form
\begin{multline}
    -\frac{N(t,x;\epsilon)}{E_\srm(t,x;\epsilon)}=-\frac{1}{2t}\left( x_\srm(t) - y_\srm(t) - y_0(t) - \sum_{j = 2}^N z_{j,\srm}(t) + \sum_{k = 1}^N p_k \right)\frac{D(t,x;\epsilon)}{E_\srm(t,x;\epsilon)} \\-\ii\frac{\epsilon^{1/3}}{2t\nu_\srm(t)}\Ai'(X_\srm)\Delta(\mathbf{w}_\srm(t),\mathbf{p}) + \bigO{\epsilon^{2/3}}.
\end{multline}
Therefore, dividing by $-D(t,x;\epsilon)/E_\srm(t,x;\epsilon)$ and using \eqref{eq:SolDenomFinalExpans}, we have
\begin{multline}
    \frac{N(t,x;\epsilon)}{D(t,x;\epsilon)}=\frac{1}{2t}\left( x_\srm(t) - y_\srm(t) - y_0(t) - \sum_{j = 2}^N z_{j,\srm}(t) + \sum_{k = 1}^N p_k \right)\\
    -\frac{\displaystyle\frac{\ii}{2t\nu_\srm(t)}\Ai'(X_\srm)+\bigO{\epsilon^{1/3}}}{\displaystyle\epsilon^{-1/3}\Ai(X_\srm) +\frac{\ii}{\nu_\srm(t)}\Ai'(X_\srm)\left(\frac{1}{\Delta(\mathbf{w}_\srm(t),\mathbf{p})}\frac{\partial \Delta}{\partial w_0}(\mathrm{w}_\srm(t),\mathbf{p})-2\frac{c_\srm(t)}{\nu_\srm(t)}\right)+\bigO{\epsilon^{1/3}}},
    \label{eq:SolNumOverDenomSecondExpan}
\end{multline}
where in the second term we divided numerator and denominator by $\epsilon^{1/3}\Delta(\mathbf{w}_\srm(t),\mathbf{p})\neq 0$.

Next, we simplify the coefficient of $\Ai'(X_\srm)$ in the denominator.  By direct computation using \eqref{eq:BODeltaDef}, we obtain
\begin{equation}
    \begin{split}
        \frac{1}{\Delta(\mathbf{w}_\srm(t), \mathbf{p})}\frac{\partial \Delta}{\partial w_0}(\mathbf{w}_\srm(t), \mathbf{p}) &= \frac{\partial \log \Delta}{\partial w_0}(\mathbf{w}_\srm(t), \mathbf{p}) \\
        &= \frac{1}{y_\srm(t) - y_0(t)} + \sum_{j = 2}^N \frac{1}{y_\srm(t) - z_{j,\srm}(t)} - \sum_{j = 1}^N \frac{1}{y_\srm(t) - p_j} \, .
    \end{split}
    \label{eq:SolEdgeLogDerivDelta}
\end{equation}
Combining this with the calculation of $2c_\srm(t)/\nu_\srm(t)$ in \eqref{eq:2covernu-s} then gives
\begin{equation}
\begin{split}
    \frac{1}{\Delta(\mathbf{w}_\srm(t), \mathbf{p})}\frac{\partial \Delta}{\partial w_0}(\mathbf{w}_\srm(t), \mathbf{p}) - 2\frac{c_\srm(t)}{\nu_\srm(t)} &= \frac{1}{2} \frac{1}{y_\srm(t) - y_0(t)} + \ii \lambda(t)\\
    &=\frac{1}{4t}\frac{1}{u_0^\Brm(t)-u_\srm^\Brm(t)} +\ii\lambda(t),
\end{split}
    \label{eq:NextToLeadingOrderCoeff}
\end{equation}
where on the second line we used \eqref{eq:SolSSaddleCond} and \eqref{eq:Sol0SaddleCond} with $u_0(y_\srm(t))=u_\srm^\Brm(t)$ and $u_0(y_0(t))=u_0^\Brm(t)$, respectively, and
\begin{align}
    \lambda(t) := \im{ \sum_{j = 2}^N \frac{1}{y_\srm(t) - z_{j,\srm}(t)} - \sum_{j = 1}^N \frac{1}{y_\srm(t) - p_j} }
    \label{eq:SolGammaDef}
\end{align}
was defined to split the real and imaginary parts of~\eqref{eq:NextToLeadingOrderCoeff}. Since at least the real part of \eqref{eq:NextToLeadingOrderCoeff} is nonzero, and since $\Ai(\diamond)$ has only simple (real) zeros, the sum of the two explicit terms in the denominator of the fraction on the second line of \eqref{eq:SolNumOverDenomSecondExpan} is uniformly bounded away from zero as $X_\srm$ varies in bounded intervals with $t \in K$, and therefore
\begin{multline}
    \frac{N(t,x;\epsilon)}{D(t,x;\epsilon)}=\frac{1}{2t}\left( x_\srm(t) - y_\srm(t) - y_0(t) - \sum_{j = 2}^N z_{j,\srm}(t) + \sum_{k = 1}^N p_k \right)\\
    -\frac{\displaystyle\frac{\ii}{2t\nu_\srm(t)}\Ai'(X_\srm)}{\displaystyle\epsilon^{-1/3}\Ai(X_\srm) +\frac{\ii}{\nu_\srm(t)}\Ai'(X_\srm)\left(\frac{1}{4t} \frac{1}{u_0^\Brm(t) - u_\srm^\Brm(t)} + \ii \lambda(t)\right)}+\bigO{\epsilon^{1/3}}.
    \label{eq:SolNumOverDenomThirdExpan}
\end{multline}
Recalling that $k_\srm(t):=2t[-\frac{1}{2}u_0''(y_\srm(t))]^{1/3}=-2t\nu_\srm(t)$, and  dividing numerator and denominator of the fraction by $\ii \Ai'(X_\srm)/(2t\nu_\srm(t))$ we obtain
\begin{multline}
    \frac{N(t,x;\epsilon)}{D(t,x;\epsilon)}=\frac{1}{2t}\left( x_\srm(t) - y_\srm(t) - y_0(t) - \sum_{j = 2}^N z_{j,\srm}(t) + \sum_{k = 1}^N p_k \right)\\
    +\frac{1}{\displaystyle\frac{1}{2}\frac{1}{u_\srm^\Brm(t)-u_0^\Brm(t)}-\ii\left(\epsilon^{-1/3}k_\srm(t)\frac{\Ai(X_\srm)}{\Ai'(X_\srm)}+2t\lambda(t)\right)} + \bigO{\epsilon^{1/3}},
\end{multline}
where we understand that the entire fraction written on the second line vanishes when $X$ is a root of $\Ai'(\diamond)$.  Then, according to Theorem~\ref{thm:IVP-solution},
\begin{multline}
    u(t,x;\epsilon)=\frac{1}{t}\mathrm{Re}\left( x_\srm(t) - y_\srm(t) - y_0(t) - \sum_{j = 2}^N z_{j,\srm}(t) + \sum_{k = 1}^N p_k \right)\\
    +L\left(\epsilon^{-1/3}k_\srm(t)\frac{\Ai(X_\srm)}{\Ai'(X_\srm)}+2t\lambda(t);c(t)\right) + \bigO{\epsilon^{1/3}},
\label{eq:SolFinal}
\end{multline}
where $L(\xi;c)$ is defined by \eqref{eq:BOSolitonProfile} and $c(t)$ is defined by \eqref{eq:SolitonSpeed}.

Now we observe that
\begin{multline}
    \frac{1}{t}\mathrm{Re}\left( x_\srm(t) - y_\srm(t) - y_0(t) - \sum_{j = 2}^N z_{j,\srm}(t) + \sum_{k = 1}^N p_k \right)\\
    =\frac{1}{2t}\left(2x_\srm(t)-2y_\srm(t)-2y_0(t)-\sum_{j=2}^N(z_{j,\srm}(t)+z_{j,\srm}(t)^*) +\sum_{k=1}^N(p_k+p_k^*)\right),
\end{multline}
so applying Lemma~\ref{lem:PoleSPRelation} at $x=x_\srm(t)$ with real critical points $\mathcal{Y}=\{y_0(t),y_\srm(t),y_\srm(t)\}$ ($y_\srm(t)$ counted twice as a double critical point), and $\mathcal{Z}=\{z_{2,\srm}(t),\dots,z_{N,\srm}(t)\}$ we obtain simply
\begin{equation}
\begin{split}
    \frac{1}{t}\mathrm{Re}\left( x_\srm(t) - y_\srm(t) - y_0(t) - \sum_{j = 2}^N z_{j,\srm}(t) + \sum_{k = 1}^N p_k \right)
    &=\frac{1}{2t}\left(x_\srm(t)-y_0(t)\right)\\
    &=u_0(y_0(t))\\ &=u_0^\Brm(t),
\end{split}
\label{eq:SolConstantTerm}
\end{equation}
where on the second line we used the critical point condition \eqref{eq:Sol0SaddleCond}.

Using \eqref{eq:SolConstantTerm} in \eqref{eq:SolFinal} nearly completes the proof of Theorem~\ref{thm:soliton-edge}. It still remains to verify the identification of $\lambda(t)$ defined by \eqref{eq:SolGammaDef} as the expression given in \eqref{eq:IntroSolGammaDef}--\eqref{eq:IntroSolgDef}.  However, we note that the function $g_\srm(y,t)$ defined in \eqref{eq:IntroSolgDef} can be written in the equivalent form
\begin{equation}
    g_\srm(y,t)=\frac{\displaystyle\prod_{j=2}^N(y-z_{j,\srm}(t))(y-z_{j,\srm}(t)^*)}{\displaystyle\prod_{k=1}^N(y-p_k)(y-p_k^*)}.
\end{equation}
It follows that $y\mapsto L_\srm(y):=\ln((y^2+1)g_\srm(y,t))$ is a function in $L^2(\mathbb{R})$ that can be split into its Cauchy-Szeg\H{o} orthogonal projection $\Pi L_\srm(y)$ onto the Hardy space of functions analytic in $\mathbb{C}_+$ and its complement $L_\srm(y)-\Pi L_\srm(y)$ whose derivative is given by
\begin{equation}
\partial_y\left(    L_\srm(y)-\Pi L_\srm(y)\right)=\frac{1}{y-\ii}+\sum_{j=2}^N\frac{1}{y-z_{j,\srm}(t)} -\sum_{k=1}^N\frac{1}{y-p_k}.
\end{equation}
Hence, comparing with \eqref{eq:SolGammaDef}, 
\begin{equation}
    \lambda(t)=\mathrm{Im}\left(\left(\partial_y(L_\srm(y)-\Pi L_\srm(y)\right)_{y=y_\srm(t)}-\frac{1}{y_\srm(t)-\ii}\right).
\end{equation}
Using the fact that differentiation commutes with projection and the singular integral representation of $\Pi$
\begin{equation}
    \Pi f(y):=\frac{1}{2\pi\ii}\lim_{\delta\downarrow 0}\int_\mathbb{R}\frac{f(s)\,\dd s}{s-(y+\ii\delta)},
\end{equation}
one then finds that $\lambda(t)$ is equivalently given by \eqref{eq:IntroSolGammaDef}--\eqref{eq:IntroSolgDef}, as desired.  The proof of Theorem~\ref{thm:soliton-edge} is now complete.

\subsubsection*{Proof of Corollary~\ref{cor:soliton-edge}}
\paragraph{On Taylor expansion of $\Ai(X_\srm)/\Ai'(X_\srm)$ and its implications}
Recall from \eqref{eq:SolitonSpeed} the definition $c(t):=2(u_\srm^\Brm(t)-u_0^\Brm(t))$, and let $a_1>a_2>a_3>\cdots$ denote the zeros of the Airy function $\Ai(\diamond)$.  Consider the difference $\Delta_\srm(X_\srm;t,\epsilon)$ defined by
\begin{multline}
    \Delta_\srm(X_\srm;t,\epsilon):=L\left(\epsilon^{-1/3}k_\srm(t)\frac{\Ai(X_\srm)}{\Ai'(X_\srm)}+ 2t\lambda(t);c(t)\right)\\-\sum_{n=1}^\infty L\left(\epsilon^{-1/3}k_\srm(t)(X_\srm-a_n)+ 2t\lambda(t);c(t)\right).
\end{multline}
This expression makes sense, because as is well-known \cite[\S 9.9(iv)]{DLMF}, 
\begin{equation}
    a_n=-\left(\frac{3}{8}\pi(4n-1)\right)^{2/3}(1+O(n^{-2})), \quad n\to\infty,
\end{equation} 
so for large $n$ the summand is $O((\epsilon^{-1/3}a_n)^{-2})=O(\epsilon^{2/3}n^{-4/3})$ which is summable.  

For bounded $X_\srm$ such that for all $n$, $|X_\srm-a_n|>\delta$ with $\delta>0$, there is another constant $\delta'>0$ such that also $|\Ai(X_\srm)/\Ai'(X_\srm)|>\delta'$ because the roots of $\Ai(\diamond)$ are all simple.  
Then $\epsilon^{-1/3}k_\srm(t)\Ai(X_\srm)/\Ai'(X_\srm)\to\infty$ as $\epsilon\to 0$ for such $X_\srm$, so since $L(\xi;c)=O(\xi^{-2})$ as $\xi\to\infty$,
\begin{equation}
    L\left(\epsilon^{-1/3}k_\srm(t)\frac{\Ai(X_\srm)}{\Ai'(X_\srm)}+2t\lambda(t);c(t)\right)=O(\epsilon^{2/3}),\quad X_\srm=O(1),\quad |X_\srm-a_n|>\delta, \forall n.
\end{equation}
Under the same conditions, the terms in the infinite series can be uniformly estimated as $O(\epsilon^{2/3}n^{-4/3})$, so the sum of those terms is $O(\epsilon^{2/3})$.
Therefore for such $X_\srm$ we deduce the estimate $\Delta_\srm(X_\srm;t,\epsilon)=O(\epsilon^{2/3})$.

For the remaining values of $X_\srm=O(1)$, there is exactly one index $n=N$ such that $|X_\srm-a_N|\le\delta$. By Taylor expansion about $X_\srm=a_N$, we get
\begin{equation}
\begin{split}
    \frac{\Ai(X_\srm)}{\Ai'(X_\srm)}&=\frac{\Ai(a_N+(X_\srm-a_N))}{\Ai'(a_N+(X_\srm-a_N))}\\ &=\frac{\Ai(a_N) + \Ai'(a_N)(X_\srm-a_N) +\frac{1}{2}\Ai''(a_N)(X_\srm-a_N)^2 + O((X_\srm-a_N)^3)}{\Ai'(a_N) + \Ai''(a_N)(X_\srm-a_N) +O((X_\srm-a_N)^2)}\\
    &=\frac{\Ai'(a_N)(X_\srm-a_N) +O((X_\srm-a_N)^3)}{\Ai'(a_N) + O((X_\srm-a_N)^2)}\\
    &=X_\srm-a_N + O((X_\srm-a_N)^3),\quad X_\srm\to a_N
\end{split}
\end{equation}
because $\Ai(a_n)=0$ and $\Ai''(a_n)=a_n\Ai(a_n)=0$ but $\Ai'(a_n)\neq 0$ for all $n$.  Therefore, by the Mean Value Theorem,
\begin{multline}
    L\left(\epsilon^{-1/3}k_\srm(t)\frac{\Ai(X_\srm)}{\Ai'(X_\srm)}+2t\lambda(t);c(t)\right)=L\left(\epsilon^{-1/3}k_\srm(t)(X_\srm-a_N)+ 2t\gamma(t);c(t)\right) \\+ O(L'(\epsilon^{-1/3}k_\srm(t)\zeta+2t\lambda(t);c(t))\epsilon^{-1/3}(X_\srm-a_N)^3)
\label{eq:single-linearized-Lorentzian-minus}
\end{multline}
where $\zeta$ is between $\Ai(X_\srm)/\Ai'(X_\srm)$ and $X_\srm-a_N$.  We use the estimate $L'(\xi;c)=O(\langle \xi\rangle^{-3})$ where $\langle\xi\rangle:=\sqrt{1+\xi^2}$, which is uniform for $c$ in any compact subset of the open half-line $c>0$.
If $\epsilon^{1/3}\ll |X_\srm-a_N|\le\delta$ as $\epsilon\to 0$, then $\xi:=\epsilon^{-1/3}k_\srm(t)\zeta+2t\lambda(t)$ is large of size $\epsilon^{-1/3}|X_\srm-a_N|$ and hence $\langle\xi\rangle^{-3}=O(\epsilon/|X_\srm-a_N|^3)$, so the error term in \eqref{eq:single-linearized-Lorentzian-minus} is $O(\epsilon^{2/3})$.  On the other hand, if $|X_\srm-a_N|=O(\epsilon^{1/3})$, then we can certainly estimate $\langle\xi\rangle^{-3}=O(1)$ and again the error term in \eqref{eq:single-linearized-Lorentzian-minus} is $O(\epsilon^{2/3})$.  
Therefore, if $|X_\srm-a_N|\le\delta$ for some $N$, we conclude that
\begin{multline}
    L\left(\epsilon^{-1/3}k_\srm(t)\frac{\Ai(X_\srm)}{\Ai'(X_\srm)}+2t\lambda(t);c(t)\right)=L\left(\epsilon^{-1/3}k_\srm(t)(X_\srm-a_N)+ 2t\lambda(t);c(t)\right) \\+ O(\epsilon^{2/3}),
    \label{eq:single-linearized-Lorentzian}
\end{multline}
so also 
\begin{equation}
    \Delta_\srm(X_\srm;t,\epsilon)=-\mathop{\sum_{n=1}^\infty}_{n\neq N}
    L\left(\epsilon^{-1/3}k_\srm(t)(X_\srm-a_{n})+ 2t\lambda(t);c(t)\right) +O(\epsilon^{2/3}).
\end{equation}
As in the previous case, under the condition $n\neq N$ the summand is $O(\epsilon^{2/3}n^{-4/3})$, so we deduce that $\Delta_\srm(X_\srm;t,\epsilon)=O(\epsilon^{2/3})$ also holds for bounded $X_\srm$ when it is \emph{not} the case that $|X_\srm-a_n|>\delta$ for all $n$.

Since $\Delta_\srm(X_\srm;t,\epsilon)=O(\epsilon^{2/3})$ can be absorbed into the error term $O(\epsilon^{1/3})$ already present in Theorem~\ref{thm:soliton-edge}, we have proved that \eqref{eq:u-soliton-edge-sum} holds under the same conditions.

\paragraph{On approximation by individual solitons}
From \eqref{eq:single-linearized-Lorentzian} we see that if for some index $n$ we have $X_\srm-a_n=o(1)$, then Theorem~\ref{thm:soliton-edge} implies
\begin{equation}
    u(t,x;\epsilon)=u_0^\Brm(t)+L\left(\epsilon^{-1/3}k_\srm(t)(X_\srm-a_n)+ 2t\lambda(t);c(t)\right)+O(\epsilon^{1/3}),\quad \epsilon\to 0.
    \label{eq:u-1-Lorentzian}
\end{equation}
The first argument of $L$ can be written as
\begin{equation}
    \epsilon^{-1/3}k_\srm(t)(X_\srm-a_n)+2t\lambda(t)=\epsilon^{-1}\left(x-x_\srm(t)-\epsilon^{2/3}k_\srm(t)a_n\right)+2t\lambda(t),
\end{equation}
and the quantity in parentheses is assumed to be $o(\epsilon^{2/3})$.  Now we expand about a fixed value $t_0>t_\cc$ using the fact that $x_\srm(t)$ is $C^2$ (at least), and that $k_\srm(t)$ and $\lambda(t)$ are both $C^1$ (at least) to get
\begin{equation}
\begin{split}
    \epsilon^{-1}\left(x-x_\srm(t)-\epsilon^{2/3}k_\srm(t)a_n\right)&=\epsilon^{-1}\left(x-x_\srm(t_0)-\epsilon^{2/3}k_\srm(t_0)a_n -x_\srm'(t_0)(t-t_0) \right)\\ &{}\quad\quad\quad + O(\epsilon^{-1}(t-t_0)^2) + O(\epsilon^{-1/3}(t-t_0)),\\ 2t\lambda(t)&=2t_0\lambda(t_0) +O(t-t_0).
    \end{split}
\end{equation}
If we now assume that
\begin{equation}
    x-x_\srm(t_0)-\epsilon^{2/3}k_\srm(t_0)a_n=o(\epsilon^{2/3})\quad\text{and}\quad t-t_0=o(\epsilon^{2/3}),
    \label{eq:x-t-soliton-scales}
\end{equation}
then the error terms produce a perturbation of the first argument of $L$ of order $o(\epsilon^{1/3})$.  Since $L'(\xi;c)$ is bounded uniformly for $c$ in compact subsets of the half-line $c>0$ and is independent of $\epsilon$, $L$ is perturbed by a quantity of the same order, and this can be absorbed into the error term already present in \eqref{eq:u-1-Lorentzian}.  Furthermore, we may expand $u_\srm^\Brm(t)=u_\srm^\Brm(t_0)+o(\epsilon^{2/3})$ and $c(t)=c(t_0)+o(\epsilon^{2/3})$ both of which lead to error terms small compared to $O(\epsilon^{1/3})$.  Recalling the offset $x_0^{(n)}(t;\epsilon)$ defined by \eqref{eq:soliton-offset}
we therefore have
\begin{equation}
    u(t,x;\epsilon)=u_0^\Brm(t_0)+L\left(\epsilon^{-1}\left(x-x_0^{(n)}(t_0;\epsilon)-x_\srm'(t_0)(t-t_0)\right);c(t_0)\right) + O(\epsilon^{1/3})
    \label{eq:u-single-soliton-almost}
\end{equation}
in the limit $\epsilon\to 0$,
provided that \eqref{eq:x-t-soliton-scales} holds.  Finally, we calculate $x_\srm'(t_0)$.  Recall that when a point $(t,x=x_\srm(t))$ lies on the caustic curve, we have
\begin{equation}
x=y+2u_0(y)t\quad \text{and}\quad 0=1+2u_0'(y)t.
\end{equation}
The first equation determines $y=y(t,x)$, and the second equation then determines a relation between the variables $(t,x)$ that we write as $x=x_\srm(t)$.  Restricting $y(t,x)$ to this curve gives us $y_\srm(t):=y(t,x_\srm(t))$.  We differentiate the first identity $x_\srm(t)=y_\srm(t)+2u_0(y_\srm(t))t$ with respect to $t$ and obtain
\begin{equation}
    x_\srm'(t)=\left(1 + 2u_0'(y_\srm(t))t\right)y_\srm'(t) +2u_0(y_\srm(t)) = 2u_0(y_\srm(t))=2u_\srm^\Brm(t),
\end{equation}
where in the second equality we used the second equation restricted to the curve $x=x_\srm(t)$:  $0=1+2u_0'(y_\srm(t))t$.  Therefore \eqref{eq:u-single-soliton-almost} becomes \eqref{eq:u-single-soliton}, and the proof of Corollary~\ref{cor:soliton-edge} is finished.

\section{Analysis near the harmonic edge $x=x_\hrm(t)$.  Proof of Theorem~\ref{thm:harmonic-edge}}
\label{sec:HarmonicEdgeProof}
The material in this section is adapted from \cite[Section 4.3]{MitchellThesis}. In the following we give the proof of Theorem~\ref{thm:harmonic-edge}, again following the three steps from Section~\ref{sec:structure-of-proofs}.

\subsection{Step~\ref{step:normal-form}:
Conformal map and normal form of the exponent}
\label{sec:harmonic-edge-step1}

To analyze the contour integrals in the first two rows of $\mathbf{A}(t,x;\epsilon)$ and $\mathbf{B}(t,x;\epsilon)$ over the contours $\widetilde{W}_0$ (first row) and $W_1$ (second row), both of which pass through the neighborhood of a pair of critical points close to $z=y_\hrm(t)$ as shown in Figure~\ref{fig:harmonic-edge-contours}, we need to apply Lemma~\ref{lem:normal-form} with $m=1$ and $n=2$ to $H(z;\mathbf{x}):=h(z;t,x)$ for fixed $t\in K=(t_\cc,t_\cc+\delta)$ but now with $z_0=y_\hrm(t)$ and $x_0=x_\hrm(t)$.  From the discussion in Section~\ref{sec:harmonic-edge-contours}, $\nu=\nu_\hrm(t)\neq 0$ in \eqref{eq:nu-power} is a cube root of 
\begin{equation}
    \nu_\hrm(t)^3=\frac{1}{2}h'''(y_\hrm(t);t,x_\hrm(t))=\frac{1}{2}u_0''(y_\hrm(t))>0.
    \label{eq:HarmNuDef}
\end{equation}
We take the positive root $\nu_\hrm(t)>0$ and then the conformal map $w=\mathcal{W}_\hrm(z;t,x)$ maps the real line near $y_\hrm(t)$ to a real neighborhood of $w=0$ in a monotone-increasing fashion.  According to Lemma~\ref{lem:normal-form}, for $z$ in a neighborhood $D$ of $y_\hrm(t)$, the identity
\begin{equation}
    h(z;t,x)=h(y_\hrm(t);t,x)+\frac{1}{3}w^3 + q_\hrm(t,x)w + r_\hrm(t,x),\quad w=\mathcal{W}_\hrm(z;t,x)
\end{equation}
holds, where $q_\hrm$ and $r_\hrm$ are $C^\infty$ functions of $x\in \mathcal{X}$ satisfying $q_\hrm(t,x_\hrm(t))=r_\hrm(t,x_\hrm(t))=0$.
In a similar way as in Section~\ref{sec:soliton-edge-step1}, one checks that
\begin{equation}
    \partial_x r_\hrm(t,x_\hrm(t))=0\quad\text{and}\quad\partial_xq_\hrm(t,x_\hrm(t))=-\frac{1}{2t\nu_\hrm(t)}.
    \label{eq:CannonicalFormCoeffsDerivs}
\end{equation}

\subsection{Step~\ref{step:expand}:  Expansions of integrals}
\label{sec:harmonic-edge-step2}
The normal form of the exponent can be used to evaluate the contour integrals in the first two rows of $\mathbf{A}(t,x;\epsilon)$ and $\mathbf{B}(t,x;\epsilon)$ which have the form
\begin{equation}
    I_0(t,x;\epsilon):=\int_{\widetilde{W}_0}\ee^{-\ii h(z;t,x)/\epsilon}f(z)\,\dd z\quad\text{and}\quad
    I_1(t,x;\epsilon):=\int_{W_1}\ee^{-\ii h(z;t,x)/\epsilon}f(z)\,\dd z,
\end{equation}
where as in Section~\ref{sec:soliton-edge-step2}, either $f(z)=u_0(z)$, $f(z)\equiv 1$, or $f(z)=(z-p_n)^{-1}$ for $n=1,\dots,N$.  According to Figure~\ref{fig:harmonic-edge-contours}, both contours $\widetilde{W}_0$ and $W_1$ pass through the neighborhood $D$ of $z=y_\hrm(t)$, while $W_1$ additionally passes over the simple saddle point $z=y_2(t,x)$ that is close to $z=y_2(t)$ for $x\approx x_\hrm(t)$.

The expansion of $I_0(t,x;\epsilon)$ is very similar to that of $I(t,x;\epsilon)$ explained in Section~\ref{sec:soliton-edge-step2}.  In this case, the analogue of \eqref{eq:soliton-I-conformal} is
\begin{multline}
    \ee^{\ii [h(y_\hrm(t);t,x)+r_\hrm(t,x)]/\epsilon}I_0(t,x;\epsilon)\\
    =\int_{\mathcal{W}_\hrm(\widetilde{W}_0\cap D;t,x)}\ee^{-\ii [w^3/3+q_\hrm(t,x)w]/\epsilon}f(\mathcal{W}_\hrm^{-1}(w;t,x))\frac{\dd\mathcal{W}_\hrm^{-1}}{\dd w}(w;t,x)\,\dd w + O(\epsilon^{\infty}),\quad\epsilon\downarrow 0.
\end{multline}
Here, since $\nu_\hrm(t)>0$, for $x\approx x_\hrm(t)$ the arc $\mathcal{W}_\hrm(\widetilde{W}_0\cap D;t,x)$ enters a neighborhood of the origin $w=0$ from the direction $\arg(w)=\frac{1}{2}\pi$ and exits in the direction $\arg(w)=-\frac{1}{6}\pi$.  Following similar steps, we obtain the two-term expansion
\begin{multline}
    \ee^{\ii [h(y_\hrm(t);t,x)+r_\hrm(t,x)]/\epsilon}I_0(t,x;\epsilon)
    =2\pi\ee^{-\ii\pi/3}\epsilon^{1/3}\Tc_\hrm^{(0)}[f](t,x_\hrm(t))\Ai(\ee^{-\ii\pi/3}X_\hrm)\\+2\pi\ee^{5\pi\ii/6}\epsilon^{2/3}\Tc_\hrm^{(1)}[f](t,x_\hrm(t))\Ai'(\ee^{-\ii\pi/3}X_\hrm)+O(\epsilon)\quad\epsilon\downarrow 0,
\label{eq:W0-h-two-term}
\end{multline}
where $X_\hrm$ defined by \eqref{eq:Xh-coordinate} is assumed to be bounded and $\Tc_\hrm^{(m)}[f](t,x)$ are the Taylor coefficients at $w=0$ of $f(\mathcal{W}_\hrm^{-1}(w;t,x))\mathcal{W}_\hrm^{-1\prime}(w;t,x)$; in particular (cf., \eqref{eq:ks0ks1} and \eqref{eq:2covernu-s})
\begin{equation}
\begin{split}
    \Tc_\hrm^{(0)}[f](t,x_\hrm(t))&=\frac{1}{\nu_\hrm(t)}f(y_\hrm(t))\\
    \Tc_\hrm^{(1)}[f](t,x_\hrm(t))&=\frac{1}{\nu_\hrm(t)^2}\left[f'(y_\hrm(t)) - \frac{2c_\hrm(t)}{\nu_\hrm(t)}f(y_\hrm(t))\right],
\end{split}
\label{eq:kh0kh1}
\end{equation}
wherein 
\begin{equation}
    \frac{2c_\hrm(t)}{\nu_\hrm(t)}=\frac{1}{2}\cdot\frac{1}{y_\hrm(t)-y_2(t)} +
    \mathrm{Re}\left(\sum_{j=2}^N\frac{1}{y_\hrm(t)-z_j(t,x_\hrm(t))}-\sum_{n=1}^N\frac{1}{y_\hrm(t)-p_n}\right).
\label{eq:2covernu-h}
\end{equation}
Therefore, we may write \eqref{eq:W0-h-two-term} in the form
\begin{multline}
    \frac{\nu_\hrm(t)\ee^{\ii\pi/3}}{2\pi\epsilon^{1/3}}\ee^{\ii[h(y_\hrm(t);t,x)+r_\hrm(t,x)]/\epsilon}I_0(t,x;\epsilon)=f(y_\hrm(t))\Ai(\ee^{-\ii\pi/3}X_\hrm)\\
    +\frac{\ee^{-5\pi\ii/6}\epsilon^{1/3}}{\nu_\hrm(t)}\left[f'(y_\hrm(t))-\frac{2c_\hrm(t)}{\nu_\hrm(t)}f(y_\hrm(t))\right]\Ai'(\ee^{-\ii\pi/3}X_\hrm)+O(\epsilon^{2/3}),\quad\epsilon\downarrow 0.
\label{eq:W0-h-two-term-rewrite}
\end{multline}

To analyze the integrals of the form $I_1(t,x;\epsilon)$ involving the path $W_1$ that passes through $D$ as well as crossing the simple saddle point $z=y_2(t,x)$, we first divide the contour $W_1$ into an arc $W_{11}$ passing over $z=y_2(t,x)$ and a subsequent arc $W_{12}$ that passes through $D$.  We then have a corresponding decomposition of $I_1(t,x;\epsilon)=I_{11}(t,x;\epsilon)+I_{12}(t,x;\epsilon)$.  For $I_{12}(t,x;\epsilon)$, we follow the above steps, first obtaining  
\begin{multline}
    \ee^{\ii [h(y_\hrm(t);t,x)+r_\hrm(t,x)]/\epsilon}I_{12}(t,x;\epsilon)\\
    =\int_{\mathcal{W}_\hrm(W_{12}\cap D;t,x)}\ee^{-\ii [w^3/3+q_\hrm(t,x)w]/\epsilon}f(\mathcal{W}_\hrm^{-1}(w;t,x))\frac{\dd\mathcal{W}_\hrm^{-1}}{\dd w}(w;t,x)\,\dd w + O(\epsilon^{\infty}),\quad\epsilon\downarrow 0.
\end{multline}
The arc $\mathcal{W}_\hrm(W_{12}\cap D;t,x)$ enters a neighborhood of the origin from the direction $\arg(w)=-\frac{5}{6}\pi$ and exits in the direction $\arg(w)=\frac{1}{2}\pi$.  Therefore, just as above, we get a two-term expansion
\begin{multline}
    \ee^{\ii [h(y_\hrm(t);t,x)+r_\hrm(t,x)]/\epsilon}I_{12}(t,x;\epsilon)=2\pi\ee^{\ii\pi/3}\epsilon^{1/3}\Tc_\hrm^{(0)}[f](t,x_\hrm(t))\Ai(\ee^{\ii\pi/3}X_\hrm)\\
    +2\pi\ee^{\ii\pi/6}\epsilon^{2/3}\Tc_\hrm^{(1)}[f](t,x_\hrm(t))\Ai'(\ee^{\ii\pi/3}X_\hrm) + O(\epsilon),\quad\epsilon\downarrow 0,
\end{multline}
again assuming that $X_\hrm$ defined by \eqref{eq:Xh-coordinate} remains bounded.  For $I_{11}(t,x;\epsilon)$, we apply the formula \eqref{eq:steepest-descent-reexpand} with fixed angle of passage $\varphi(t,x)=-\frac{1}{4}\pi$ and simple critical point $\zeta(t,x)=y_2(t,x)$.  We also identify $t$ and $t_0$, and take $x_0=x_\hrm(t)$ so that $x-x_\hrm(t)=O(\epsilon^{2/3})$ holds for $X_\hrm$ bounded.  Thus,
\begin{equation}
\begin{split}
    I_{11}(t,x;\epsilon)&=\ee^{-\ii\pi/4}\sqrt{\frac{2\pi\epsilon}{|h''(y_2(t,x);t,x)|}}\ee^{-\ii h(y_2(t,x);t,x)/\epsilon}\left(f(y_2(t)) + O(\epsilon^{2/3})\right)\\
    &=\ee^{-\ii\pi/4}\sqrt{\frac{2\pi\epsilon}{|h''(y_2(t);t,x_\hrm(t))|}}\ee^{-\ii h(y_2(t,x);t,x_\hrm(t))/\epsilon}\left(f(y_2(t))+O(\epsilon^{2/3})\right),
\end{split}
\end{equation}
where on the second line we further used $y_2(t,x)=y_2(t)+O(\epsilon^{2/3})$ when $x-x_\hrm(t)=O(\epsilon^{2/3})$ and absorbed the error term into that already present on the first line.  Noting that $h''(y_2(t);t,x_\hrm(t))>0$ because $y_2(t,x)$ is the left-most real critical point of $h(z;t,x)$,
we combine these results using \eqref{eq:kh0kh1} to obtain the following approximation for $I_1(t,x;\epsilon)$:  
\begin{multline}
    \sqrt{\frac{h''(y_2(t);t,x_\hrm(t))}{2\pi\epsilon}}\ee^{\ii\pi/4}\ee^{\ii h(y_2(t,x);t,x)/\epsilon}I_1(t,x;\epsilon) = \\
\epsilon^{-1/6}\ee^{\ii\phi(t,x)/\epsilon}B_\hrm(t)f(y_\hrm(t))\ee^{7\pi\ii/12}\Ai(\ee^{\ii\pi/3}X_\hrm) + f(y_2(t))\\
    +\epsilon^{1/6}\ee^{\ii\phi(t,x)/\epsilon}B_\hrm(t)\frac{1}{\nu_\hrm(t)}\left[f'(y_\hrm(t))-\frac{2c_\hrm(t)}{\nu_\hrm(t)}f(y_\hrm(t))\right]\ee^{5\pi\ii/12}\Ai'(\ee^{\ii\pi/3}X_\hrm) + O(\epsilon^{1/2})
    \label{eq:W1-h-final}
\end{multline}
as $\epsilon\downarrow 0$, where we have defined
\begin{equation}
\begin{split}
    B_\hrm(t)&:=\frac{1}{\nu_\hrm(t)}\sqrt{2\pi h''(y_2(t);t,x_\hrm(t))}\\
    \phi(t,x)&:=h(y_2(t,x);t,x)-h(y_\hrm(t);t,x)-r_\hrm(t,x).
\end{split}
\label{eq:Bh-phi}
\end{equation}

We introduce the shorthand notation $z_{j,\hrm}(t):=z_j(t,x_\hrm(t))$ for the complex critical points of $h(\diamond;t,x_\hrm(t))$.  Then, since $h'(\diamond;t,x_\hrm(t))$ is a rational function having simple poles only at $p_1,\dots,p_N$ and $p_1^*,\dots,p_N^*$, simple zeros only at $z=y_2(t)$ and the complex values $z=z_{j,\hrm}(t)$ and $z=z_{j,\hrm}(t)^*$ for $j=2,\dots,N$, and a double zero at $z=y_\hrm(t)$, otherwise nonvanishing, and $h'(z;t,x_\hrm(t))=z/(2t)+O(1)$ as $z\to\infty$, we have
\begin{equation}
    h'(z;t,x_\hrm(t))=\frac{1}{2t}(z-y_\hrm(t))^2(z-y_2(t))\frac{\displaystyle\prod_{j=2}^N(z-z_{j,\hrm}(t))(z-z_{j,\hrm}(t)^*)}{\displaystyle\prod_{n=1}^N(z-p_n)(z-p_n^*)}.
    \label{eq:Harm-hprime-product}
\end{equation}
Taking one more derivative and evaluating at the simple root $z=y_2(t)\in\mathbb{R}$ gives
\begin{equation}
    h''(y_2(t);t,x_\hrm(t))=\frac{1}{2t}(y_2(t)-y_\hrm(t))^2\frac{\displaystyle\prod_{j=2}^N|y_2(t)-z_{j,\hrm}(t)|^2}{\displaystyle\prod_{n=1}^N|y_2(t)-p_n|^2},
    \label{eq:hdoubleprime-harmonic-edge-y2}
\end{equation}
which is obviously positive. 
Similarly, taking two more derivatives of \eqref{eq:Harm-hprime-product} and evaluating at the double root $z=y_\hrm(t)\in\mathbb{R}$ gives
\begin{equation}
    h'''(y_\hrm(t);t,x_\hrm(t))=\frac{1}{t}(y_\hrm(t)-y_2(t))\frac{\displaystyle\prod_{j=2}^N|y_\hrm(t)-z_{j,\hrm}(t)|^2}{\displaystyle\prod_{n=1}^N|y_\hrm(t)-p_n|^2}.
\label{eq:htripleprime-harmonic-edge-yh}
\end{equation}
The latter expression is equal to $2\nu_\hrm(t)^3$ according to \eqref{eq:HarmNuDef}.

\subsection{Step~\ref{step:algebra}:  Simplification of determinants}
\label{sec:harmonic-edge-step3}
We will now apply the classical steepest descent method in the form \eqref{eq:steepest-descent-reexpand} (taking into account that $x=x_\mathrm{h}(t)+O(\epsilon^{2/3})$) to expand the integrals in rows $3$ through $N+1$, and also \eqref{eq:W0-h-two-term-rewrite} and \eqref{eq:W1-h-final} to expand those in the first and second row respectively, of the matrices $\mathbf{A}(t,x;\epsilon)$ and $\mathbf{B}(t,x;\epsilon)$.
According to \eqref{eq:W0-h-two-term-rewrite}, the first row of each matrix therefore decomposes into a sum of two row vectors plus an error, while from \eqref{eq:W1-h-final} the second row of each matrix decomposes into a sum of three row vectors plus an error.  The determinants $N(t,x;\epsilon):=\det(\mathbf{A}(t,x;\epsilon))$ and $D(t,x;\epsilon):=\det(\mathbf{B}(t,x;\epsilon))$ then are similarly decomposed via multilinearity.  

We define a common nonzero scalar factor $E_\hrm(t,x;\epsilon)$ for both determinants by
\begin{equation}
    \begin{split}
        E_\hrm(t,x;\epsilon) &:= \frac{2\pi\epsilon^{1/3}}{\nu_\hrm(t)}\ee^{-\ii\pi/3}\ee^{-\ii[h(y_\hrm(t);t,x)+r_\hrm(t,x)]/\epsilon}\\
        &\Quad[2] \times\sqrt{\frac{2\pi\epsilon}{h''(y_2(t);t,x_\hrm(t))}}\ee^{-\ii\pi/4}\ee^{-\ii h(y_2(t,x);t,x)/\epsilon}\\
        &\Quad[2]\times\prod_{m=2}^N\ee^{\ii\varphi_m(t,x)}\sqrt{\frac{2\pi\epsilon}{|h''(z_m(t,x);t,x)|}}\ee^{-\ii h(z_m(t,x);t,x)/\epsilon}.
    \end{split}
    \label{eq:HarExpFactor}
\end{equation}
We recall that $z_m(t,x)$, $m=2,\dots,N$ are the simple complex critical points of $h(\diamond;t,x)$ in the upper half-plane and for each, $\varphi_m(t,x)$ is the tangent angle of the steepest descent path with which $W_m$ traverses the saddle point. Then, using row multilinearity and accounting for obvious cancellations, we obtain
\begin{equation}
    \frac{N(t,x;\epsilon)}{E_\hrm(t,x;\epsilon)}=\Ai(\ee^{-\ii\pi/3}X_\hrm)N_{\hrm,2}(t) 
    +
    \epsilon^{1/6}\ee^{\ii\phi(t,x)/\epsilon}\frac{B_\hrm(t)}{\nu_\hrm(t)}V(X_\hrm) N_{\hrm,\hrm}^{(1)}(t) + \bigO{\epsilon^{1/3}},
\label{eq:HarmNumExpansion1}
\end{equation}
where $V(X_\hrm)$ is defined in \eqref{eq:VofXh} below while
$N_{\hrm,2}(t)$ and $N_{\hrm,\hrm}^{(1)}(t)$ are the determinants
\begin{equation}
    N_{\hrm,2}(t):=\begin{vmatrix}       
    u_0(y_\hrm(t)) & (y_\hrm(t)-p_1)^{-1} & \cdots & (y_\hrm(t)-p_N)^{-1}\\
    u_0(y_2(t)) & (y_2(t)-p_1)^{-1} & \cdots & (y_2(t)-p_N)^{-1}\\
    u_0(z_{2,\hrm}(t)) & (z_{2,\hrm}(t)-p_1)^{-1} & \cdots & (z_{2,\hrm}(t)-p_N)^{-1}\\
    \vdots & \vdots & \ddots & \vdots\\
    u_0(z_{N,\hrm}(t)) & (z_{N,\hrm}(t)-p_1)^{-1} & \cdots & (z_{N,\hrm}(t)-p_N)^{-1}
    \end{vmatrix}
\end{equation}
and
\begin{equation}
    N_{\hrm,\hrm}^{(1)}(t):=\begin{vmatrix}
        u_0'(y_\hrm(t)) & -(y_\hrm(t)-p_1)^{-2} & \cdots & -(y_\hrm(t)-p_N)^{-2}\\
        u_0(y_\hrm(t)) & (y_\hrm(t)-p_1)^{-1} & \cdots  & (y_\hrm(t)-p_N)^{-1}\\
        u_0(z_{2,\hrm}(t)) & (z_{2,\hrm}(t)-p_1)^{-1} & \cdots & (z_{2,\hrm}(t)-p_N)^{-1}\\
        \vdots & \vdots & \ddots & \vdots\\
        u_0(z_{N,\hrm}(t)) & (z_{N,\hrm}(t)-p_1)^{-1} & \cdots & (z_{N,\hrm}(t)-p_N)^{-1}
    \end{vmatrix}.
\end{equation}
Exactly as in Section~\ref{sec:soliton-edge-step3}, the critical point equations
\begin{equation}
\begin{split}
    0 &= \frac{y_\hrm(t) - x_\hrm(t)}{2t} + u_0(y_\hrm(t)),  \\
    0 &= \frac{y_2(t) - x_\hrm(t)}{2t} + u_0(y_2(t)), \\
    0 &= \frac{z_{j,\hrm}(t) - x_\hrm(t)}{2t} + u_0(z_{j,\hrm}(t)),\quad 2 \leq j \leq N,
\end{split}
\label{eq:HarmCriticalConds}
\end{equation}
allow the determinant $N_{\hrm,2}(t)$ to be expressed in terms of the determinants $C^{(0)}(\mathbf{w},\mathbf{p})$ and $C^{(1)}(\mathbf{w},\mathbf{p})$ defined in \eqref{eq:CauchyDetFirstCoeff}--\eqref{eq:CauchyDetSecondCoeff} as follows:
\begin{equation}
    N_{\hrm,2}(t)=\frac{x_\hrm(t)}{2t}C^{(0)}(\mathbf{w}_{\hrm,2}(t),\mathbf{p})-\frac{1}{2t}C^{(1)}(\mathbf{w}_{\hrm,2}(t),\mathbf{p}),
\end{equation}
where
\begin{equation}
        \mathbf{w}_{\hrm,2}(t):=(y_\hrm(t),y_2(t),z_{2,\hrm}(t),\dots,z_{N,\hrm}(t))^\top
\label{eq:w-hrm-2}
\end{equation}
and $\mathbf{p}$ is given by \eqref{eq:p-vector}.  Using also the double critical point equation satisfied by $y_\hrm(t)$ at $(t,x_\hrm(t))$:
\begin{equation}
    0 = h''(y_\hrm(t);t, x_\hrm(t)) = \frac{1}{2 t} + u_0'(y_\hrm(t)),
\end{equation}
we then get that 
\begin{equation}
    N^{(1)}_{\hrm,\hrm}(t)=\frac{x_\hrm(t)}{2t}\frac{\partial C^{(0)}}{\partial w_0}(\mathbf{w}_{\hrm,\hrm}(t),\mathbf{p})-\frac{1}{2t}\frac{\partial C^{(1)}}{\partial w_0}(\mathbf{w}_{\hrm,\hrm}(t),\mathbf{p}),
\end{equation}
where
\begin{equation}
 \mathbf{w}_{\hrm,\hrm}(t):=(y_\hrm(t),y_\hrm(t),z_{2,\hrm}(t),\dots,z_{N,\hrm}(t))^\top.    
 \label{eq:w-hrm-hrm}
\end{equation}

The expansion of $D(t,x;\epsilon)/E_\hrm(t,x;\epsilon)$ is similar, but the determinants $N_{\hrm,2}(t)$ and $N_{\hrm,\hrm}^{(1)}(t)$ are replaced by others differing only in the first column, where the function $u_0(\diamond)$ is replaced by the constant function $1$.  Hence the latter determinants are expressible directly in terms of the Cauchy determinant $C^{(0)}(\mathbf{w},\mathbf{p})$ defined in \eqref{eq:CauchyDetFirstCoeff}:
\begin{multline}
    \frac{D(t,x;\epsilon)}{E_\hrm(t,x;\epsilon)}=\Ai(\ee^{-\ii\pi/3}X_\hrm)C^{(0)}(\mathbf{w}_{\hrm,2}(t),\mathbf{p}) 
    \\+
    \epsilon^{1/6}\ee^{\ii\phi(t,x)/\epsilon}\frac{B_\hrm(t)}{\nu_\hrm(t)}V(X_\hrm) \frac{\partial C^{(0)}}{\partial w_0}(\mathbf{w}_{\hrm,\hrm}(t),\mathbf{p}) + \bigO{\epsilon^{1/3}},
\label{eq:HarmDenExpansion1}
\end{multline}

In both expansions \eqref{eq:HarmNumExpansion1} and \eqref{eq:HarmDenExpansion1}, $V(X_\hrm)$ denotes the bilinear expression
\begin{equation}
    V(X_\hrm):=\ee^{-\ii\pi/4}\Ai'(\ee^{-\ii\pi/3}X_\hrm)\Ai(\ee^{\ii\pi/3}X_\hrm)-\ee^{5\ii\pi/12}\Ai(\ee^{-\ii\pi/3}X_\hrm)\Ai'(\ee^{\ii\pi/3}X_\hrm).
    \label{eq:VofXh}
\end{equation}
Note that $V(X_\hrm)$ is proportional to the Wronskian of solutions $\Ai(z\ee^{-2\pi\ii/3})$ and $\Ai(z\ee^{2\pi\ii/3})$ of Airy's equation $y''=zy$ evaluated for $z=-X_\hrm$ and this Wronskian is evaluated explicitly in \cite[Eqn.\@ 9.2.9]{DLMF}.  Hence in fact $V(X_\hrm)$ is independent of $X_\hrm$ and is simply
\begin{equation}
    V(X_\hrm)=\frac{\ee^{7\ii\pi/12}}{2\pi}.
\end{equation}

Using this, and noting that the leading term $\Ai(\ee^{-\ii\pi/3}X_\hrm)C^{(0)}(\mathbf{w}_{\hrm,2}(t),\mathbf{p})$ is nonvanishing for real $X_\hrm$, it is straightforward to expand the reciprocal of \eqref{eq:HarmDenExpansion1} and multiply by \eqref{eq:HarmNumExpansion1} to obtain
\begin{multline}
    \frac{N(t,x;\epsilon)}{D(t,x;\epsilon)}=\frac{x_\hrm(t)-C^{(1)}(\mathbf{w}_{\hrm,2}(t),\mathbf{p})/C^{(0)}(\mathbf{w}_{\hrm,2}(t),\mathbf{p})}{2t} \\+
    \frac{\epsilon^{1/6}\ee^{7\ii\pi/12}\ee^{\ii\phi(t,x)/\epsilon}B_\hrm(t)}{4\pi t\nu_\hrm(t)\Ai(\ee^{-\ii\pi/3}X_\hrm)}\left[\frac{C^{(1)}(\mathbf{w}_{\hrm,2}(t),\mathbf{p})}{C^{(0)}(\mathbf{w}_{\hrm,2}(t),\mathbf{p})^2}\frac{\partial C^{(0)}}{\partial w_0}(\mathbf{w}_{\hrm,\hrm}(t),\mathbf{p})\right.\\
    \left.-\frac{1}{C^{(0)}(\mathbf{w}_{\hrm,2}(t),\mathbf{p})}\frac{\partial C^{(1)}}{\partial w_0}(\mathbf{w}_{\hrm,\hrm}(t),\mathbf{p})\right] + \bigO{\epsilon^{1/3}}.
\end{multline}
According to Lemma~\ref{lem:Cauchy-determinants} and using \eqref{eq:w-hrm-2}, we have
\begin{equation}
    \frac{C^{(1)}(\mathbf{w}_{\hrm,2}(t),\mathbf{p})}{C^{(0)}(\mathbf{w}_{\hrm,2}(t),\mathbf{p})} = y_\hrm(t)+y_2(t) + \sum_{j=2}^Nz_{j,\hrm}(t)-\sum_{n=1}^Np_n,
\end{equation}
because $\Delta(\mathbf{w}_{\hrm,2}(t),\mathbf{p})\neq 0$.  Similarly, applying Lemma~\ref{lem:Cauchy-determinants} using also \eqref{eq:w-hrm-hrm} and noting that $\Delta(\mathbf{w}_{\hrm,\hrm}(t),\mathbf{p})=0$ one obtains
\begin{multline}
 \frac{C^{(1)}(\mathbf{w}_{\hrm,2}(t),\mathbf{p})}{C^{(0)}(\mathbf{w}_{\hrm,2}(t),\mathbf{p})^2}\frac{\partial C^{(0)}}{\partial w_0}(\mathbf{w}_{\hrm,\hrm}(t),\mathbf{p})-\frac{1}{C^{(0)}(\mathbf{w}_{\hrm,2}(t),\mathbf{p})}\frac{\partial C^{(1)}}{\partial w_0}(\mathbf{w}_{\hrm,\hrm}(t),\mathbf{p})\\=(y_2(t)-y_\hrm(t))\frac{1}{\Delta(\mathbf{w}_{\hrm,2}(t),\mathbf{p})}\frac{\partial \Delta}{\partial w_0}(\mathbf{w}_{\hrm,\hrm}(t),\mathbf{p}).   
\end{multline}
Then, a direct calculation using the product definition \eqref{eq:BODeltaDef} gives
\begin{equation}
\frac{1}{\Delta(\mathbf{w}_{\hrm,2}(t),\mathbf{p})}\frac{\partial \Delta}{\partial w_0}(\mathbf{w}_{\hrm,\hrm}(t),\mathbf{p}) = \frac{1}{y_\hrm(t)-y_2(t)}\cdot\prod_{j=2}^N\frac{z_{j,\hrm}(t)-y_\hrm(t)}{z_{j,\hrm}(t)-y_2(t)}\cdot\prod_{n=1}^N\frac{y_2(t)-p_n}{y_\hrm(t)-p_n}.    
\end{equation}
Therefore,
\begin{multline}
    \frac{N(t,x;\epsilon)}{D(t,x;\epsilon)}=\frac{1}{2t}\left[x_\hrm(t)-y_\hrm(t)-y_2(t)-\sum_{j=2}^Nz_{j,\hrm}(t)+\sum_{n=1}^Np_n\right]\\
    +\frac{\epsilon^{1/6}\ee^{-5\ii\pi/12}\ee^{\ii\phi(t,x)/\epsilon}B_\hrm(t)}{4\pi t\nu_\hrm(t)\Ai(\ee^{-\ii\pi/3}X_\hrm)}\prod_{j=2}^N\frac{z_{j,\hrm}(t)-y_\hrm(t)}{z_{j,\hrm}(t)-y_2(t)}\cdot\prod_{n=1}^N\frac{y_2(t)-p_n}{y_\hrm(t)-p_n} + \bigO{\epsilon^{1/3}}.
\label{eq:HarmNoverDFinal}
\end{multline}

Combining the definition of $B_\hrm(t)$ in \eqref{eq:Bh-phi} with \eqref{eq:hdoubleprime-harmonic-edge-y2} and using \eqref{eq:HarmNuDef} with \eqref{eq:htripleprime-harmonic-edge-yh} shows that, aside from the factor of $\epsilon^{1/6}/|\Ai(\ee^{-\ii\pi/3}X_\hrm)|=\epsilon^{1/6}/|\Ai(\ee^{\ii\pi/3}X_\hrm)|$ (for $X_\hrm\in\mathbb{R}$), the correction term in \eqref{eq:HarmNoverDFinal} has modulus
\begin{equation}
\begin{split}
    \frac{B_\hrm(t)}{2\pi t\nu_\hrm(t)}\left|\prod_{j=2}^N\frac{z_{j,\hrm}(t)-y_\hrm(t)}{z_{j,\hrm}(t)-y_2(t)}\cdot\prod_{n=1}^N\frac{y_2(t)-p_n}{y_\hrm(t)-p_n}\right| &= \frac{1}{2}\sqrt{\frac{2(u_0(y_2(t))-u_0(y_\hrm(t)))}{\pi k_\hrm(t)}}\\
    &=\frac{1}{2}\sqrt{\frac{2(u_2^\Brm(t)-u_\hrm^\Brm(t))}{\pi k_\hrm(t)}},
\end{split}
\end{equation}
where we also used the critical point conditions \eqref{eq:HarmCriticalConds} for $y_\hrm(t)$ and $y_2(t)$, and we recall $k_\hrm(t)=2t[\frac{1}{2}u_0''(y_\hrm(t))]^{1/3}=2t[\frac{1}{2}h'''(y_\hrm(t))]^{1/3}=2t\nu_\hrm(t)$.  Since $\Ai(\ee^{-\ii\pi/3}X_\hrm)=\Ai(\ee^{\ii\pi/3}X_\hrm)^*$ for $X_\hrm\in\mathbb{R}$, the phase of the correction term in \eqref{eq:HarmNoverDFinal} is then $\phi(t,x)/\epsilon +\varphi(t,X_\hrm)$, where
\begin{equation}
    \varphi(t,X_\hrm):=\arg\left(\Ai(\ee^{\ii\pi/3}X_\hrm)\right)-\frac{5\pi}{12}+\Lambda_\hrm(t)-\Lambda_2(t),
\end{equation}
with
\begin{equation}
\begin{split}
    \Lambda_\hrm(t)&:=\sum_{j=2}^N\arg\left(z_{j,\hrm}(t)-y_\hrm(t)\right) -\sum_{n=1}^N\arg\left(y_\hrm(t)-p_n\right)\\
    \Lambda_2(t)&:=\sum_{j=2}^N\arg\left(z_{j,\hrm}(t)-y_2(t)\right)-\sum_{n=1}^N\arg\left(y_2(t)-p_n\right).
\end{split}
\label{eq:HarmProofLambdah2}
\end{equation}
Also, we observe that for the leading term,
\begin{multline}
    2\mathrm{Re}\left(\frac{1}{2t}\left[x_\hrm(t)-y_\hrm(t)-y_2(t)-\sum_{j=2}^Nz_{j,\hrm}(t)+\sum_{n=1}^Np_n\right]\right)\\
    \begin{aligned}
    &=\frac{1}{2t}\left[2x_\hrm(t)-2y_\hrm(t)-2y_2(t)-\sum_{j=2}^N\left(z_{j,\hrm}(t)+z_{j,\hrm}(t)^*\right) +\sum_{n=1}^N\left(p_n+p_n^*\right)\right]\\
    &= \frac{x_\hrm(t)-y_2(t)}{2t}\\
    &= u_0(y_2(t))\\
    &=u_2^\Brm(t),
    \end{aligned}
\end{multline}
where we used Lemma~\ref{lem:PoleSPRelation} with real critical points $\mathcal{Y}=\{y_2(t),y_\hrm(t),y_\hrm(t)\}$ (repeated real root $y_\hrm(t)$) and $\mathcal{Z}=\{z_{2,\hrm}(t),\dots,z_{N,\hrm(t)}\}$, and then applied the critical point condition in \eqref{eq:HarmCriticalConds} for $y_2(t)$.  Therefore, according to Theorem~\ref{thm:IVP-solution}, \eqref{eq:HarmNoverDFinal} gives $u(t,x;\epsilon)$ in the form 
\begin{equation}
    u(t,x;\epsilon)=u_2^\Brm(t) + \epsilon^{1/6}\sqrt{\frac{2\left(u_2^\Brm(t)-u_\hrm^\Brm(t)\right)}{\pi k_\hrm(t)}}\frac{\cos(\epsilon^{-1}\phi(t,x) + \varphi(t,X_\hrm))}{|\Ai(\ee^{\ii\pi/3}X_\hrm)|} + \bigO{\epsilon^{1/3}}.
\label{eq:Harm-proof-u}
\end{equation}

To complete the proof of Theorem~\ref{thm:harmonic-edge}, it only remains to identify $\phi(t,x)$ with a suitable perturbation of $\Theta(t,x)$ defined in \eqref{eq:BOHarmonicPhase-0}, and to show that $\Lambda_\hrm(t)-\Lambda_2(t)$ defined in \eqref{eq:HarmProofLambdah2} coincides with $\Lambda(t)$ as given by the expressions \eqref{eq:BOHarmonicPhaseh2}--\eqref{eq:BOHarmonicgdef}.  

Since $r_\hrm(t,x_\hrm(t))=0$ and according to \eqref{eq:CannonicalFormCoeffsDerivs} the $C^\infty$ function $r_\hrm$ satisfies $\partial_x r_\hrm(t,x_\hrm(t))$, it follows that for bounded $X_\hrm$ we have $r_\hrm(t,x)=O(\epsilon^{4/3})$ because $x=x_\hrm(t)+O(\epsilon^{2/3})$.  Therefore from the definition of $\phi(t,x)$ in \eqref{eq:Bh-phi},
\begin{equation}
    \phi(t,x)=h(y_2(t,x);t,x)-h(y_\hrm(t);t,x) + \bigO{\epsilon^{4/3}}
\end{equation}
holds for bounded $X_\hrm$, uniformly for $t\in K$.  The quantity $h(y_2(t,x;t,x)-h(y_\hrm(t);t,x)$ is analytic in $x$ and so by Taylor expansion about $x=x_\hrm(t)$ using $x=x_\hrm(t)+O(\epsilon^{2/3})$ gives
\begin{multline}
    \phi(t,x)=h(y_2(t);t,x_\hrm(t))-h(y_\hrm(t);t,x_\hrm(t)) \\+ (x-x_\hrm(t))\frac{\partial}{\partial x}\left[h(y_2(t,x);t,x)-h(y_\hrm(t);t,x)\right]\bigg|_{x=x_\hrm(t)} + \bigO{\epsilon^{4/3}}.
\label{eq:phi-expand}
\end{multline}
The leading terms can be written as
\begin{multline}
    h(y_2(t);t,x_\hrm(t))-h(y_\hrm(t);t,x_\hrm(t))\\
    \begin{aligned}
        &=\int_{y_\hrm(t)}^{y_2(t)}h'(y;t,x_\hrm(t))\,\dd y\\
    &=\int_{y_\hrm(t)}^{y_2(t)}\frac{y-x_\hrm(t)}{2t}\,\dd y + \int_{y_\hrm(t)}^{y_2(t)}u_0(y)\,\dd y\\
    &= \frac{(y_2(t)-x_\hrm(t))^2}{4t}-\frac{(y_\hrm(t)-x_\hrm(t))^2}{4t} +\int_{y_\hrm(t)}^{y_2(t)}u_0(y)\,\dd y\\
    &= \left(u_0(y_2(t))^2-u_0(y_\hrm(t))^2\right)t+ \int_{y_\hrm(t)}^{y_2(t)}u_0(y)\,\dd y\\
    &= \left(u_2^\Brm(t)^2-u_\hrm^\Brm(t)^2\right)t + \int_{y_\hrm(t)}^{y_2(t)}u_0(y)\,\dd y,
    \end{aligned}
\label{eq:phi-leading-terms}
\end{multline}
where we used the critical point conditions for $y_2(t)$ and $y_\hrm(t)$ from \eqref{eq:HarmCriticalConds}. Here, in the beginning we integrate along the contour arc of $W_1$ between the two critical points $y_\hrm(t)$ and $y_2(t)$ since we know that it does not cross any branch cuts of $h(\diamond;t,x)$ defined by 
\eqref{eq:h-def}. Then, because there are no poles of $u_0$ between this arc of $W_1$ and the real line, we can deform the contour to the real line in the last step.  
Similarly, computing the derivative using $h'(y_2(t);t,x_\hrm(t))=0$, only differentiation of $h$ with respect to the parameter $x$ is required and one obtains
\begin{equation}
\begin{split}
    \frac{\partial}{\partial x}\left[h(y_2(t,x);t,x)-h(y_\hrm(t);t,x)\right]\bigg|_{x=x_\hrm(t)}&=\frac{x_\hrm(t)-y_2(t)}{2t} -\frac{x_\hrm(t)-y_\hrm(t)}{2t} \\ &= u_0(y_2(t))-u_0(y_\hrm(t)) \\ &= u_2^\Brm(t)-u_\hrm^\Brm(t),
\end{split}
\label{eq:phi-deriv-terms}
\end{equation}
where we again used the critical point conditions \eqref{eq:HarmCriticalConds}.  Using \eqref{eq:phi-leading-terms} and \eqref{eq:phi-deriv-terms} in \eqref{eq:phi-expand} then shows that
\begin{equation}
    \epsilon^{-1}\phi(t,x)=\epsilon^{-1}\Theta(t,x) + \bigO{\epsilon^{1/3}}
\end{equation}
holds uniformly for $X_\hrm$ bounded and $t\in K$, where $\Theta(t,x)$ is as defined in \eqref{eq:BOHarmonicPhase-0}.  Therefore, one may replace $\phi(t,x)$ with $\Theta(t,x)$ in \eqref{eq:Harm-proof-u} at the cost of a modification of the error term of the same order.

Finally, to show that $\Lambda_\hrm(t)-\Lambda_2(t)$ as defined in \eqref{eq:HarmProofLambdah2} agrees with $\Lambda(t)$ defined in \eqref{eq:BOHarmonicPhaseh2}--\eqref{eq:BOHarmonicgdef}, first observe that $g_\hrm(y,t)$ defined in \eqref{eq:BOHarmonicgdef} can be equivalently represented as
\begin{equation}
    g_\hrm(y,t)=\frac{\displaystyle \prod_{j=2}^N(y-z_{j,\hrm}(t))(y-z_{j,\hrm}(t)^*)}{\displaystyle\prod_{k=1}^N(y-p_k)(y-p_k^*)},
\end{equation}
and just as in Section~\ref{sec:soliton-edge-step3}, the function $y\mapsto L_\hrm(y):=\ln((y^2+1)g_\hrm(y,t))$ is in $L^2(\mathbb{R})$ and its complementary Cauchy-Szeg\H{o} projection is 
\begin{equation}
    L_\hrm(y)-\Pi L_\hrm(y)=\log(y-\ii) + \sum_{j=2}^N\log(y-z_{j,\hrm}(t)) -\sum_{k=1}^N\log(y-p_k)
\end{equation}
where the logarithms may be taken as principal branches.  Hence, comparing with \eqref{eq:HarmProofLambdah2}, we see that
\begin{equation}
    \begin{split}
        \Lambda_\hrm(t)&=\mathrm{Im}\left(L_\hrm(y_\hrm(t))-\Pi L_\hrm(y_\hrm(t))-\log(y_\hrm(t)-\ii) +(N-1)\ii \pi\right)\\
        \Lambda_2(t)&=\mathrm{Im}\left(L_\hrm(y_2(t))-\Pi L_\hrm(y_2(t))-\log(y_2(t)-\ii) +(N-1)\ii \pi\right).
    \end{split}
\end{equation}
Exactly as in Section~\ref{sec:soliton-edge-step3}, a singular integral calculation yields, for any $y\in\mathbb{R}$,
\begin{equation}
\begin{split}
    \mathrm{Im}\left(L_\hrm(y)-\Pi L_\hrm(y)\right) &= \frac{1}{2\pi}\int_0^{+\infty}\left(L_\hrm(y+s)-L_\hrm(y-s)\right)\frac{\dd s}{s}\\
    &=\frac{1}{2\pi}\int_0^{+\infty}\ln\left(\frac{(y+s)^2+1}{(y-s)^2+1}\right)\frac{\dd s}{s} +\frac{1}{2\pi}\int_0^{+\infty}\ln\left(\frac{g_\hrm(y+s,t)}{g_\hrm(y-s,t)}\right)\frac{\dd s}{s}.
\end{split}
\end{equation}
The first integral can be evaluated explicitly:
\begin{equation}
    \frac{1}{2\pi}\int_0^{+\infty}\ln\left(\frac{(y+s)^2+1}{(y-s)^2+1}\right)\frac{\dd s}{s}=\frac{\pi}{2}+\arg(y-\ii),\quad y\in\mathbb{R}.
\end{equation}
Consequently we obtain
\begin{equation}
    \begin{split}
        \Lambda_\hrm(t)&=(N-\tfrac{1}{2})\pi + \frac{1}{2\pi}\int_0^{+\infty}\ln\left(\frac{g_\hrm(y_\hrm(t)+s,t)}{g_\hrm(y_\hrm(t)-s,t)}\right)\frac{\dd s}{s}\\
        \Lambda_2(t)&=(N-\tfrac{1}{2})\pi + \frac{1}{2\pi}\int_0^{+\infty}\ln\left(\frac{g_\hrm(y_2(t)+s,t)}{g_\hrm(y_2(t)-s,t)}\right)\frac{\dd s}{s},
    \end{split}
\end{equation}
and it follows that $\Lambda(t)=\Lambda_\hrm(t)-\Lambda_2(t)$ as desired. This completes the proof of Theorem~\ref{thm:harmonic-edge}.

\section{Analysis near the gradient catastrophe point.  Proof of Theorem~\ref{thm:BO-catastrophe}}
\label{sec:GradientCatastropheProof}

Finally, we prove Theorem~\ref{thm:BO-catastrophe}, again following the three steps of Section~\ref{sec:structure-of-proofs}.
\subsection{Step~\ref{step:normal-form}:  Conformal map and normal form of the exponent}
\label{sec:catastrophe-step1}
When $(t,x)\approx (t_\cc,x_\cc)$, there are three critical points of $h(z;t,x)$ in a small neighborhood $D$ of $z=y_\cc$, and as shown in Figure~\ref{fig:catastrophe-contours} the contours $\widetilde{W}_0$ and $W_1$ for the integrals in the first two rows respectively of both matrices $\mathbf{A}(t,x;\epsilon)$ and $\mathbf{B}(t,x;\epsilon)$ pass through $D$.  To evaluate these integrals, we will now invoke Lemma~\ref{lem:normal-form} with $m=2$ and $n=3$, taking $z_0=y_\cc$, $\mathbf{x}=(t,x)$, $\mathbf{x}_0=(t_\cc,x_\cc)$ and $H(z;\mathbf{x})=h(z;t,x)$.  By the discussion in Section~\ref{sec:catastrophe-contours}, the parameter $\nu=\nu_\cc\neq 0$ in \eqref{eq:nu-power} satisfies
\begin{equation}
    \nu_\cc^4 = \frac{1}{6}h''''(y_\cc;t_\cc,x_\cc)=\frac{1}{6}u_0'''(y_\cc)>0.
\label{eq:nu-cc}
\end{equation}
We choose the positive real fourth root $\nu_\cc>0$, so that the conformal map $w=\mathcal{W}_\cc(z;t,x)$ whose existence is asserted in Lemma~\ref{lem:normal-form} is Schwarz-symmetric and monotone increasing, taking the real line near $y_\cc$ to an interval containing $w=0$.  The same Lemma then yields that under $w=\mathcal{W}_\cc(z;t,x)$, 
\begin{equation}
    h(z;t,x)=h(y_\cc;t,x)+\frac{1}{4}w^4 + q_\cc(t,x)w^2 + r_\cc(t,x)w + s_\cc(t,x),\quad w=\mathcal{W}_\cc(z;t,x),
\label{eq:normal-form-catastrophe}
\end{equation}
for real-valued functions $q_\cc$, $r_\cc$, and $s_\cc$ of $(t,x)\in \mathcal{X}$ that are of class $C^\infty$ and that satisfy $q_\cc(t_\cc,x_\cc)=r_\cc(t_\cc,x_\cc)=s_\cc(t_\cc,x_\cc)=0$.

Next, we compute the partial derivatives of these three functions at $(t_\cc,x_\cc)$ using \eqref{eq:partial-a0-crit}--\eqref{eq:partial-ak-crit}.
We use the facts that
\begin{equation}
    \partial_t h(z;t,x)=-\frac{(z-x)^2}{4t^2} \quad\text{and}\quad \partial_x h(z;t,x)= -\frac{z-x}{2t}
\end{equation}
and obtain from \eqref{eq:partial-a0-crit} that
\begin{equation}
    \partial_ts_\cc(t_\cc,x_\cc)=-\frac{1}{4t_\cc^2}\frac{1}{2\pi\ii}\oint_L\frac{(z-x_\cc)^2\mathcal{W}_\cc'(z;t_\cc,x_\cc)}{\mathcal{W}_\cc(z;t_\cc,x_\cc)}\,\dd z +\frac{(y_\cc-x_\cc)^2}{4t_\cc^2}
\end{equation}
and
\begin{equation}
    \partial_xs_\cc(t_\cc,x_\cc)=-\frac{1}{2t_\cc}\frac{1}{2\pi\ii}\oint_L\frac{(z-x_\cc)\mathcal{W}'_\cc(z;t_\cc,x_\cc)}{\mathcal{W}_\cc(z;t_\cc,x_\cc)}\,\dd z +\frac{y_\cc-x_\cc}{2t_\cc},
\end{equation}
while \eqref{eq:partial-ak-crit} gives
\begin{equation}
    \partial_tr_\cc(t_\cc,x_\cc)=-\frac{1}{4t_\cc^2}\frac{1}{2\pi\ii}\oint_L\frac{(z-x_\cc)^2\mathcal{W}_\cc'(z;t_\cc,x_\cc)}{\mathcal{W}_\cc(z;t_\cc,x_\cc)^2}\,\dd z,
\end{equation}
\begin{equation}
    \partial_tq_\cc(t_\cc,x_\cc)=-\frac{1}{4t_\cc^2}\frac{1}{2\pi\ii}\oint_L\frac{(z-x_\cc)^2\mathcal{W}_\cc'(z;t_\cc,x_\cc)}{\mathcal{W}_\cc(z;t_\cc,x_\cc)^3}\,\dd z,
\end{equation}
\begin{equation}
    \partial_xr_\cc(t_\cc,x_\cc)=-\frac{1}{2t_\cc}\frac{1}{2\pi\ii}\oint_L\frac{(z-x_\cc)\mathcal{W}_\cc'(z;t_\cc,x_\cc)}{\mathcal{W}_\cc(z;t_\cc,x_\cc)^2}\,\dd z,
\end{equation}
\begin{equation}
    \partial_xq_\cc(t_\cc,x_\cc)=-\frac{1}{2t_\cc}\frac{1}{2\pi\ii}\oint_L\frac{(z-x_\cc)\mathcal{W}_\cc'(z;t_\cc,x_\cc)}{\mathcal{W}_\cc(z;t_\cc,x_\cc)^3}\,\dd z.
\end{equation}
Then using the Taylor expansion $\mathcal{W}_\cc(z;t_\cc,x_\cc)=\nu_\cc(z-y_\cc) + c_\cc(z-y_\cc)^2 + O((z-y_\cc)^3)$, a residue computation gives
\begin{equation}
    \partial_ts_\cc(t_\cc,x_\cc)=\partial_x s_\cc(t_\cc,x_\cc)=0,
\end{equation}
\begin{equation}
    \partial_tr_\cc(t_\cc,x_\cc)=-\frac{y_\cc-x_\cc}{2t_\cc^2\nu_\cc},\quad\partial_x r_\cc(t_\cc,x_\cc)=-\frac{1}{2t_\cc\nu_\cc},
\end{equation}
and
\begin{equation}
    \partial_tq_\cc(t_\cc,x_\cc)=\frac{1}{4t_\cc^2\nu_\cc^2}\left(\frac{2c_\cc}{\nu_\cc}(y_\cc-x_\cc)-1\right),\quad \partial_xq_\cc(t_\cc,x_\cc)=\frac{c_\cc}{2t_\cc\nu_\cc^3}.
\end{equation}
Since $0=h'(y_\cc;t_\cc,x_\cc)=(y_\cc-x_\cc)/(2t_\cc)+u_0(y_\cc)$, writing $u_\cc:=u_0(y_\cc)$, the two nonzero $t$ derivatives can be written in the form
\begin{equation}
    \partial_tr_\cc(t_\cc,x_\cc)=\frac{u_\cc}{t_\cc\nu_\cc},\quad
    \partial_tq_\cc(t_\cc,x_\cc)=-\frac{c_\cc u_\cc}{t_\cc\nu_\cc^3}-\frac{1}{4t_\cc^2\nu_\cc^2}.
\end{equation}
We then obtain the corresponding Taylor expansions of $q_\cc$, $r_\cc$, and $s_\cc$ about $(t_\cc,x_\cc)$ in the form
\begin{equation}
\begin{split}
q_\cc(t,x)&=\frac{c_\cc}{2t_\cc\nu_\cc^3}\left[x-x_\cc-2u_\cc(t-t_\cc)\right] -\frac{1}{4t_\cc^2\nu_\cc^2}(t-t_\cc)+ \cdots\\
r_\cc(t,x)&=-\frac{1}{2t_\cc\nu_\cc}\left[x-x_\cc-2u_\cc(t-t_\cc)\right]+\cdots\\
s_\cc(t,x)&=\cdots
\end{split}
\label{eq:qrs-catastrophe-Taylor}
\end{equation}
where in each case the ellipsis indicates terms of quadratic and higher order in $x-x_\cc$ and $t-t_\cc$.

\subsection{Step~\ref{step:expand}:  Expansions of integrals}
\label{sec:catastrophe-step2}
The next step is to expand the integrals in the first and second rows of $\mathbf{A}(t,x;\epsilon)$ and $\mathbf{B}(t,x;\epsilon)$ over the contours $W_0$ and $W_1$ respectively as shown in Figure~\ref{fig:catastrophe-contours}.  These have the form
\begin{equation}
    I_j(t,x;\epsilon)=\int_{W_j}\ee^{-\ii h(z;t,x)/\epsilon}f(z)\,\dd z,\quad j=0,1,
\end{equation}
with $f(z)=u_0(z)$, $f(z)\equiv 1$, or $f(z)=(z-p_n)^{-1}$ for $n=1,\dots,N$.  Factoring out constants of unit modulus and truncating the contour outside a small neighborhood $D$ of $z=y_\cc$ at the cost of an exponentially small error gives
\begin{equation}
    \ee^{\ii[h(y_\cc;t,x)+s_\cc(t,x)]/\epsilon}I_j(t,x;\epsilon)=
    \int_{W_j\cap D}\ee^{-\ii[h(z;t,x)-h(y_\cc;t,x)-s_\cc(t,x)]/\epsilon}f(z)\,\dd z + O(\epsilon^{\infty}),\quad\epsilon\downarrow 0.
\end{equation}
Introducing the conformal mapping $w=\mathcal{W}_\cc(z;t,x)$ and using \eqref{eq:normal-form-catastrophe} gives, in the same limit $\epsilon\downarrow 0$,
\begin{multline}
\ee^{\ii[h(y_\cc;t,x)+s_\cc(t,x)]/\epsilon}I_j(t,x;\epsilon)\\
=\int_{\mathcal{W}_\cc(W_j\cap D;t,x)}\ee^{-\ii[w^4/4+q_\cc(t,x)w^2+r_\cc(t,x)w]/\epsilon}f(\mathcal{W}_\cc^{-1}(w;t,x))\frac{\dd\mathcal{W}_\cc^{-1}}{\dd w}(w;t,x)\,\dd w + O(\epsilon^\infty). 
\end{multline}
Since $\nu_\cc>0$, the arc $\mathcal{W}_\cc(W_0\cap D;t,x)$ enters a neighborhood of the origin $w=0$ from the direction $\arg(w)=\frac{7}{8}\pi$ and exits in the direction $\arg(w)=-\frac{1}{8}\pi$.  Similarly, the arc $\mathcal{W}_\cc(W_1\cap D;t,x)$ enters the same neighborhood of $w=0$ from the direction $\arg(w)=\frac{7}{8}\pi$ and exits in the direction $\arg(w)=\frac{3}{8}\pi$.  See Figure~\ref{fig:catastrophe-contours}.  By Taylor expansion about $w=0$ we define coefficients $\Tc_\cc^{(m)}[f](t,x)$ by
\begin{equation}
    f(\mathcal{W}_\cc^{-1}(w;t,x))\frac{\dd\mathcal{W}_\cc^{-1}}{\dd w}(w;t,x)=\sum_{m=0}^M\Tc_\cc^{(m)}[f](t,x)w^m + O(w^{M+1}).  
\label{eq:kcc-define}
\end{equation}
It follows that 
\begin{multline}
    \ee^{\ii[h(y_\cc;t,x)+s_\cc(t,x)]/\epsilon}I_j(t,x;\epsilon)\\
    =\sum_{m=0}^M\Tc_\cc^{(m)}[f](t,x)\int_{\mathcal{W}_\cc(W_j\cap D;t,x)}\ee^{-\ii[w^4/4 + q_\cc(t,x)w^2+r_\cc(t,x)w]/\epsilon}w^m\,\dd w \\
    +O\left(\int_{\mathcal{W}_\cc(W_j\cap D;t,x)}\left|\ee^{-\ii[w^4/4+q_\cc(t,x)w^2+r_\cc(t,x)w]/\epsilon}\right||w|^{M+1}\,|\dd w|\right) + O(\epsilon^\infty).
\end{multline}
Rescaling by $w=\ee^{-\ii\pi/4}(4\epsilon)^\frac{1}{4}\zeta$ and setting 
\begin{equation}
    Y:=-2\epsilon^{-1/2}q_\cc(t,x)\quad\text{and}\quad Z:=-\sqrt{2}\epsilon^{-3/4}r_\cc(t,x)
\end{equation}
gives
\begin{multline}
\ee^{\ii[h(y_\cc;t,x)+s_\cc(t,x)]/\epsilon}I_j(t,x;\epsilon)\\
    =\sum_{m=0}^M \ee^{-\ii\pi (m+1)/4}(4\epsilon)^{(m+1)/4}\Tc_\cc^{(m)}[f](t,x)    
    \int_{\ee^{\ii\pi/4}(4\epsilon)^{-1/4}\mathcal{W}_\cc(W_j\cap D;t,x)}\ee^{\ii (\zeta^4 - \ii Y\zeta^2+\ee^{-\ii\pi/4}Z\zeta)}\zeta^m\,\dd\zeta\\
    +O\left(\epsilon^{(M+2)/4}\int_{\ee^{\ii\pi/4}(4\epsilon)^{-1/4}\mathcal{W}_\cc(W_j\cap D;t,x)}\left|\ee^{\ii (\zeta^4 - \ii Y\zeta^2+\ee^{-\ii\pi/4}Z\zeta)}\right||\zeta|^{M+1}\,|\dd\zeta|\right) + O(\epsilon^\infty).
\end{multline}
It follows from the definition of the coordinates $X$ and $T$ in \eqref{eq:catastrophe-coordinates} and the Taylor expansions \eqref{eq:qrs-catastrophe-Taylor} that $(T,X)$ bounded in $\mathbb{R}^2$ implies that also 
\begin{equation}
    Y = T + O(\epsilon^{1/4}) \quad\text{and}\quad Z=X+O(\epsilon^{1/4}),
\label{eq:YZ-perturb}
\end{equation}
so the coordinates $Y$ and $Z$ are also bounded.  Then, the path of integration in each integral can be extended to $\zeta=\infty$ in suitable directions at the cost of terms absorbable into the $O(\epsilon^\infty)$ term already present.  Taking $M=2$ we therefore obtain
\begin{multline}
  \ee^{\ii[h(y_\cc;t,x)+s_\cc(t,x)]/\epsilon}I_0(t,x;\epsilon)  =
  \ee^{-\ii\pi/4}(4\epsilon)^{1/4}\Tc_\cc^{(0)}[f](t,x)P(-\ii Y,\ee^{-\ii\pi/4}Z)\\ -2 \epsilon^{1/2}\Tc_\cc^{(1)}[f](t,x)\frac{\partial P}{\partial\beta}(-\ii Y,\ee^{-\ii\pi/4}Z)\\
  +\ee^{\ii\pi/4}(4\epsilon)^{3/4}\Tc_\cc^{(2)}[f](t,x)\frac{\partial^2 P}{\partial\beta^2}(-\ii Y,\ee^{-\ii\pi/4}Z)+ O(\epsilon)
\label{eq:I0-catastrophe}
\end{multline}
and
\begin{multline}
  \ee^{\ii[h(y_\cc;t,x)+s_\cc(t,x)]/\epsilon}I_1(t,x;\epsilon)  =
  \ee^{-\ii\pi/4}(4\epsilon)^{1/4}\Tc_\cc^{(0)}[f](t,x)Q(-\ii Y,\ee^{-\ii\pi/4}Z) \\-2 \epsilon^{1/2}\Tc_\cc^{(1)}[f](t,x)\frac{\partial Q}{\partial\beta}(-\ii Y,\ee^{-\ii\pi/4}Z)\\
  +\ee^{\ii\pi/4}(4\epsilon)^{3/4}\Tc_\cc^{(2)}[f](t,x)\frac{\partial^2Q}{\partial \beta^2}(-\ii Y,\ee^{-\ii\pi/4}Z)+ O(\epsilon),
\label{eq:I1-catastrophe}
\end{multline}
wherein 
\begin{equation}
    P(\alpha,\beta):=\int_{\ee^{-7\pi\ii/8}\infty}^{\ee^{\ii\pi/8}\infty}\ee^{\ii (\zeta^4+\alpha\zeta^2+\beta\zeta)}\,\dd\zeta\quad\text{and}\quad
    Q(\alpha,\beta):=\int_{\ee^{-7\pi\ii/8}\infty}^{\ee^{5\pi\ii/8}\infty}\ee^{\ii (\zeta^4 +\alpha\zeta^2+\beta\zeta)}\,\dd\zeta.
\label{eq:PQ-define}
\end{equation}
In particular, by a contour rotation, we see that $P(\alpha,\beta)$ coincides with the Pearcey diffraction integral defined in \cite[36.2.14]{DLMF}.

\subsection{Step~\ref{step:algebra}:  Simplification of the determinants}
\label{sec:catastrophe-step3}
Now we use the classical steepest descent method in the form \eqref{eq:steepest-descent} to expand the integrals in rows $3$ through $N+1$, together with \eqref{eq:I0-catastrophe} and \eqref{eq:I1-catastrophe} to expand the integrals in rows $1$ and $2$ respectively, of the matrices $\mathbf{A}(t,x;\epsilon)$ and $\mathbf{B}(t,x;\epsilon)$.  In particular, \eqref{eq:I0-catastrophe} and \eqref{eq:I1-catastrophe} allow the representation of each of the first two rows of these matrices as a sum of three row vectors and an error.  We then use multilinearity of the determinant to expand $N(t,x;\epsilon):=\det(\mathbf{A}(t,x;\epsilon))$ and $D(t,x;\epsilon):=\det(\mathbf{B}(t,x;\epsilon))$ according to the indicated decomposition of the first two rows.

Define a common scalar factor $E_\cc(t,x;\epsilon)$ by
\begin{equation}
E_\cc(t,x;\epsilon):=\ee^{-2\ii[h(y_\cc;t,x)+s_\cc(t,x)]/\epsilon}\prod_{m=2}^N\ee^{\ii\varphi_m(t,x)}\sqrt{\frac{2\pi\epsilon}{|h''(z_m(t,x);t,x)|}}\ee^{-\ii h(z_m(t,x);t,x)/\epsilon},
\end{equation}
where we recall from Section~\ref{sec:catastrophe-contours} that $z_2(t,x),\dots,z_N(t,x)$ denote the (simple) complex critical points of $z\mapsto h(z;t,x)$ in the upper half-plane, and $\varphi_m(t,x)$ denotes the tangent angle of steepest descent passage of the contour $W_m$ as it traverses the saddle point $z=z_m(t,x)$.  Also, define elementary determinants by
\begin{equation}
    N^{(i,j)}_\cc(t,x):=\begin{vmatrix} \Tc_\cc^{(i)}[u_0](t,x) & \Tc_\cc^{(i)}[(\diamond-p_1)^{-1}](t,x) & \cdots &\Tc_\cc^{(i)}[(\diamond-p_N)^{-1}](t,x)\\
    \Tc_\cc^{(j)}[u_0](t,x) & \Tc_\cc^{(j)}[(\diamond-p_1)^{-1}](t,x)& \cdots &\Tc_\cc^{(j)}[(\diamond-p_N)^{-1}](t,x)\\
    u_0(z_2(t,x)) & (z_2(t,x)-p_1)^{-1} & \cdots &(z_2(t,x)-p_N)^{-1}\\
    \vdots & \vdots & \ddots & \vdots\\
    u_0(z_N(t,x)) & (z_N(t,x)-p_1)^{-1} & \cdots &(z_N(t,x)-p_N)^{-1}\end{vmatrix},
\end{equation}
and
\begin{equation}
    D^{(i,j)}_\cc(t,x):=\begin{vmatrix} \Tc_\cc^{(i)}[1](t,x) & \Tc_\cc^{(i)}[(\diamond-p_1)^{-1}](t,x) & \cdots &\Tc_\cc^{(i)}[(\diamond-p_N)^{-1}](t,x)\\
    \Tc_\cc^{(j)}[1](t,x) & \Tc_\cc^{(j)}[(\diamond-p_1)^{-1}](t,x)& \cdots &\Tc_\cc^{(j)}[(\diamond-p_N)^{-1}](t,x)\\
    1 & (z_2(t,x)-p_1)^{-1} & \cdots &(z_2(t,x)-p_N)^{-1}\\
    \vdots & \vdots & \ddots & \vdots\\
    1 & (z_N(t,x)-p_1)^{-1} & \cdots &(z_N(t,x)-p_N)^{-1}\end{vmatrix},
\end{equation}
for $i,j=0,1,2$.  Clearly these are antisymmetric in the indices $i,j$.  Using the antisymmetry, it follows that
\begin{multline}
    \frac{N(t,x;\epsilon)}{E_\cc(t,x;\epsilon)}=\ee^{-\ii\pi/4}(4\epsilon)^{3/4}N_\cc^{(1,0)}(t,x)W_\beta[P,Q](-\ii Y,\ee^{-\ii\pi/4}Z)\\
    -4\epsilon N_\cc^{(2,0)}(t,x)\frac{\partial W_\beta[P,Q]}{\partial\beta}(-\ii Y,\ee^{-\ii\pi/4}Z) + O(\epsilon^{5/4}),
\end{multline}
and
\begin{multline}
    \frac{D(t,x;\epsilon)}{E_\cc(t,x;\epsilon)}=\ee^{-\ii\pi/4}(4\epsilon)^{3/4}D_\cc^{(1,0)}(t,x)W_\beta[P,Q](-\ii Y,\ee^{-\ii\pi/4}Z)\\
    -4\epsilon D_\cc^{(2,0)}(t,x)\frac{\partial W_\beta[P,Q]}{\partial\beta}(-\ii Y,\ee^{-\ii\pi/4}Z) + O(\epsilon^{5/4}),    
\end{multline}
where $W_\beta[P,Q](\alpha,\beta)$ denotes the Wronskian of $P,Q$ with respect to $\beta$:
\begin{equation}
    W_\beta[P,Q](\alpha,\beta):=P(\alpha,\beta)\frac{\partial Q}{\partial\beta}(\alpha,\beta)-\frac{\partial P}{\partial\beta}(\alpha,\beta)Q(\alpha,\beta).
\label{eq:WronskianPQ}
\end{equation}
Since $N_\cc^{(i,j)}(t,x)=N_\cc^{(i,j)}(t_\cc,x_\cc)+O(\epsilon^{1/2})$ and $D_\cc^{(i,j)}(t,x)=D_\cc^{(i,j)}(t_\cc,x_\cc)+O(\epsilon^{1/2})$ both hold for all indices $(i,j)$ when $T$ and $X$ defined by \eqref{eq:catastrophe-coordinates} are bounded, for different error terms of the same order we can write
\begin{multline}
    \frac{N(t,x;\epsilon)}{E_\cc(t,x;\epsilon)}=\ee^{-\ii\pi/4}(4\epsilon)^{3/4}N_\cc^{(1,0)}(t_\cc,x_\cc)W_\beta[P,Q](-\ii Y,\ee^{-\ii\pi/4}Z)\\
    -4\epsilon N_\cc^{(2,0)}(t_\cc,x_\cc)\frac{\partial W_\beta[P,Q]}{\partial\beta}(-\ii Y,\ee^{-\ii\pi/4}Z) + O(\epsilon^{5/4}),
\label{eq:Numerator-catastrophe}
\end{multline}
and
\begin{multline}
    \frac{D(t,x;\epsilon)}{E_\cc(t,x;\epsilon)}=\ee^{-\ii\pi/4}(4\epsilon)^{3/4}D_\cc^{(1,0)}(t_\cc,x_\cc)W_\beta[P,Q](-\ii Y,\ee^{-\ii\pi/4}Z)\\
    -4\epsilon D_\cc^{(2,0)}(t_\cc,x_\cc)\frac{\partial W_\beta[P,Q]}{\partial\beta}(-\ii Y,\ee^{-\ii\pi/4}Z) + O(\epsilon^{5/4}).    
\label{eq:Denominator-catastrophe}
\end{multline}

\subsubsection{The coefficients $N_\cc^{(i,0)}(t_\cc,x_\cc)$ and $D_\cc^{(i,0)}(t_\cc,x_\cc)$}
Now we evaluate the determinants $N_\cc^{(i,0)}(t_\cc,x_\cc)$ and $D_\cc^{(i,0)}(t_\cc,x_\cc)$, for $i=1,2$.  
Note that if $\mathcal{W}_\cc(z;t_\cc,x_\cc)=\nu_\cc(z-y_\cc)+c_\cc(z-y_\cc)^2+d_\cc(z-y_\cc)^3+O((z-y_\cc)^4)$, then the inverse of the conformal mapping  has the corresponding expansion $\mathcal{W}^{-1}_\cc(w;t_\cc,x_\cc)=y_\cc+w/\nu_\cc -c_\cc w^2/\nu_\cc^3 +(2c_\cc^2/\nu_\cc^5-d_\cc/\nu_\cc^4)w^3 +O(w^4)$.  It then follows that for any function $z\mapsto f(z)$ analytic at $z=y_\cc$, 
\begin{multline}
f(\mathcal{W}_\cc^{-1}(w;t_\cc,x_\cc))\frac{\dd\mathcal{W}^{-1}}{\dd w}(w;t_\cc,x_\cc) =
\frac{f(y_\cc)}{\nu_\cc} + \frac{\nu_\cc f'(y_\cc)-2c_\cc f(y_\cc)}{\nu_\cc^3}w \\+
\frac{\nu_\cc^2 f''(y_\cc)-6c_\cc\nu_\cc f'(y_\cc) + [12c_\cc^2-6d_\cc\nu_\cc]f(y_\cc)}{2\nu_\cc^5}w^2+O(w^3),\quad w\to 0.
\end{multline}
Hence, referring to \eqref{eq:kcc-define},
\begin{equation}
\begin{split}
    \Tc^{(0)}_\cc[f](t_\cc,x_\cc)&=\frac{f(y_\cc)}{\nu_\cc},\\
    \Tc^{(1)}_\cc[f](t_\cc,x_\cc)&=\frac{\nu_\cc f'(y_\cc)-2c_\cc f(y_\cc)}{\nu_\cc^3},\\
    \Tc^{(2)}_\cc[f](t_\cc,x_\cc)&=\frac{\nu_\cc^2 f''(y_\cc)-6c_\cc\nu_\cc f'(y_\cc) + [12c_\cc^2-6d_\cc\nu_\cc]f(y_\cc)}{2\nu_\cc^5}.
\end{split}
\end{equation}
By antisymmetry of a determinant in the first two rows, the contribution of the term proportional to $f(y_\cc)$ in $\Tc^{(1)}[f](t_\cc,x_\cc)$ will cancel from $N^{(1,0)}_\cc(t_\cc,x_\cc)$ and $D^{(1,0)}_\cc(t_\cc,x_\cc)$ and therefore one sees that
\begin{equation}
    N_\cc^{(1,0)}(t_\cc,x_\cc)=\frac{1}{\nu_\cc^3}\begin{vmatrix} u_0'(y_\cc) & -(y_\cc-p_1)^{-2} & \cdots & -(y_\cc-p_N)^{-2}\\
    u_0(y_\cc) & (y_\cc-p_1)^{-1} & \cdots & (y_\cc-p_N)^{-1}\\
    u_0(z_2(t_\cc,x_\cc)) & (z_2(t_\cc,x_\cc)-p_1)^{-1} &\cdots & (z_2(t_\cc,x_\cc)-p_N)^{-1}\\
    \vdots  & \vdots & \ddots & \vdots\\
    u_0(z_N(t_\cc,x_\cc)) & (z_N(t_\cc,x_\cc)-p_1)^{-1} & \cdots & (z_N(t_\cc,x_\cc)-p_N)^{-1}\end{vmatrix},
\label{eq:Nc01}
\end{equation}
and, from the definition of $C^{(0)}(\mathbf{w}, \mathbf{p})$ \eqref{eq:CauchyDetFirstCoeff},
\begin{equation}
    D_\cc^{(1,0)}(t_\cc,x_\cc) = \frac{1}{\nu_\cc^3} \frac{\partial C^{(0)}}{\partial w_0}(\mathbf{w}_\cc, \mathbf{p}),
\end{equation}
where
\begin{equation}
 \mathbf{w}_\cc := (y_\cc, y_\cc, z_2(t_\cc, x_\cc), \dots, z_N(t_\cc, x_\cc))^\top ,
 \label{eq:w-c}
\end{equation}
and $\mathbf{p}$ is defined in \eqref{eq:p-vector}. 
Similarly, in computing $N^{(2,0)}_\cc(t_\cc,x_\cc)$ and $D^{(2,0)}_\cc(t_\cc,x_\cc)$, the term in $\Tc^{(2)}[f](t_\cc,x_\cc)$ proportional to $f(y_\cc)$ will cancel, while the term proportional to $f'(y_\cc)$ will give a contribution proportional to $N^{(1,0)}_\cc(t_\cc,x_\cc)$ and $D^{(1,0)}_\cc(t_\cc,x_\cc)$ respectively:
\begin{multline}
    N_\cc^{(2,0)}(t_\cc,x_\cc)=-\frac{3c_\cc}{\nu_\cc^2}N_\cc^{(1,0)}(t_\cc,x_\cc)\\
    +\frac{1}{2\nu_\cc^4}\begin{vmatrix}u_0''(y_\cc) & 2(y_\cc-p_1)^{-3} & \cdots & 2(y_\cc-p_N)^{-3}\\
    u_0(y_\cc) & (y_\cc-p_1)^{-1} & \cdots & (y_\cc-p_N)^{-1}\\
        u_0(z_2(t_\cc,x_\cc)) & (z_2(t_\cc,x_\cc)-p_1)^{-1} &\cdots & (z_2(t_\cc,x_\cc)-p_N)^{-1}\\
    \vdots  & \vdots & \ddots & \vdots\\
    u_0(z_N(t_\cc,x_\cc)) & (z_N(t_\cc,x_\cc)-p_1)^{-1} & \cdots & (z_N(t_\cc,x_\cc)-p_N)^{-1}\end{vmatrix},
\label{eq:Nc02}
\end{multline}
and
\begin{equation}
    D_\cc^{(2,0)}(t_\cc,x_\cc)=-\frac{3c_\cc}{\nu_\cc^2}D_\cc^{(1,0)}(t_\cc,x_\cc) + \frac{1}{2\nu_\cc^4} \frac{\partial^2 C^{(0)}}{\partial w_0^2}(\mathbf{w}_\cc, \mathbf{p}) .
\end{equation}
To evaluate the determinants on the right-hand sides of \eqref{eq:Nc01} and \eqref{eq:Nc02}, we use the fact that at the critical point $(t_\cc,x_\cc)$, we have the three conditions $h'(y_\cc;t_\cc,x_\cc)=0$, $h''(y_\cc;t_\cc,x_\cc)=0$, and $h'''(y_\cc;t_\cc,x_\cc)=0$, which can be rewritten explicitly as
\begin{equation}
    u_0(y_\cc)=\frac{x_\cc-y_\cc}{2t_\cc},\quad u_0'(y_\cc)=-\frac{1}{2t_\cc},\quad u_0''(y_\cc)=0.
\end{equation}
Similarly, because each simple critical point $z_m(t_\cc,x_\cc)$, $m=2,\dots,N$, satisfies the identity $h'(z_m(t_\cc,x_\cc);t_\cc,x_\cc)=0$, 
\begin{equation}
u_0(z_m(t_\cc,x_\cc))=\frac{x_\cc-z_m(t_\cc,x_\cc)}{2t_\cc},\quad m=2,\dots,N.
\end{equation}
Therefore, we observe that the determinants $N_\cc^{(1,0)}(t_\cc,x_\cc)$ and $N_\cc^{(2,0)}(t_\cc,x_\cc)$ are related to derivatives of the determinants $C^{(0)}(\mathbf{w},\mathbf{p})$ and $C^{(1)}(\mathbf{w},\mathbf{p})$ defined in Lemma~\ref{lem:Cauchy-determinants}: 
\begin{equation}
    \begin{split}
    N_\cc^{(1,0)}(t_\cc,x_\cc) &= \frac{1}{2t_\cc\nu_\cc^3} \left( x_\cc \frac{\partial C^{(0)}}{\partial w_0}(\mathbf{w}_\cc, \mathbf{p}) + \frac{\partial C^{(1)}}{\partial w_0}(\mathbf{w}_\cc, \mathbf{p}) \right) \\
    &= \frac{1}{2t_\cc\nu_\cc^3}\left(x_\cc-2y_\cc - \sum_{m=2}^N z_m(t_\cc,x_\cc) + \sum_{n=1}^N p_n \right) \frac{\partial\Delta}{\partial w_0}(\mathbf{w}_\cc, \mathbf{p}) ,
    \end{split}
    \label{eq:N10c-final}
\end{equation}
where we used the fact that $\Delta(\mathbf{w}_\cc, \mathbf{p}) = 0$, and
\begin{equation}
    \begin{split}
        N_\cc^{(2,0)}(t_\cc,x_\cc) &= -\frac{3c_\cc}{\nu_\cc^2}N_\cc^{(1,0)}(t_\cc,x_\cc) +\frac{1}{4t_\cc\nu_\cc^4} \left( x_\cc \frac{\partial^2 C^{(0)}}{\partial w_0^2}(\mathbf{w}_\cc, \mathbf{p}) + \frac{\partial^2 C^{(1)}}{\partial^2 w_0}(\mathbf{w}_\cc, \mathbf{p}) \right) \\
        &= -\frac{3c_\cc}{\nu_\cc^2}N_\cc^{(1,0)}(t_\cc,x_\cc) - \frac{1}{2t_\cc\nu_\cc^4} \frac{\partial\Delta}{\partial w_0}(\mathbf{w}_\cc,\mathbf{p}) \\
        &\Quad[2] +\frac{1}{4t_\cc\nu_\cc^4}\left( x_\cc - 2y_\cc - \sum_{m=2}^N z_m(t_\cc,x_\cc) + \sum_{n=1}^N p_n \right) \frac{\partial^2\Delta}{\partial w_0^2}(\mathbf{w}_\cc, \mathbf{p}) ,
    \end{split}
    \label{eq:N20c-final}
\end{equation}
after evaluating the derivatives of the determinants $C^{(0)}(\mathbf{w},\mathbf{p})$ and $C^{(1)}(\mathbf{w},\mathbf{p})$ in terms of $\Delta(\mathbf{w}, \mathbf{p})$ according to Lemma~\ref{lem:Cauchy-determinants}. Next, we observe that the determinants $D_\cc^{(1,0)}(t_\cc,x_\cc)$ and $D_\cc^{(2,0)}(t_\cc,x_\cc)$ are related to derivatives of only the determinant $C^{(0)}(\mathbf{w},\mathbf{p})$:
\begin{equation}
    \begin{split}
        D^{(1,0)}_\cc(t_\cc,x_\cc) &= \frac{1}{\nu_\cc^3} \frac{\partial C^{(0)}}{\partial w_0}(\mathbf{w}_\cc, \mathbf{p})\\
        &=\frac{1}{\nu_\cc^3} \frac{\partial \Delta}{\partial w_0}(\mathbf{w}_\cc, \mathbf{p}),
    \end{split}
    \label{eq:D10c}
\end{equation}
and
\begin{equation}
    \begin{split}
        D^{(2,0)}_\cc(t_\cc,x_\cc) &= -\frac{3c_\cc}{\nu_\cc^2} D_\cc^{(1,0)}(t_\cc,x_\cc) + \frac{1}{2\nu_\cc^4} \frac{\partial^2 C^{(0)}}{\partial w_0^2}(\mathbf{w}_\cc, \mathbf{p}) \\
        &= -\frac{3c_\cc}{\nu_\cc^2} D_\cc^{(1,0)}(t_\cc,x_\cc) + \frac{1}{2\nu_\cc^4} \frac{\partial^2\Delta}{\partial w_0^2}(\mathbf{w}_\cc, \mathbf{p}) .
    \end{split}
    \label{eq:D20c}
\end{equation}

Combining the definition \eqref{eq:BODeltaDef} with \eqref{eq:D10c}, one sees that $D_\cc^{(1,0)}(t_\cc,x_\cc)\neq 0$.
It then follows from \eqref{eq:D10c} and \eqref{eq:N10c-final} that 
\begin{equation}
    \frac{N_\cc^{(1,0)}(t_\cc,x_\cc)}{D_\cc^{(1,0)}(t_\cc,x_\cc)} = \frac{1}{2t_\cc}\left(x_\cc-2y_\cc-\sum_{m=2}^Nz_m(t_\cc,x_\cc)+\sum_{n=1}^Np_n\right),
\label{eq:catastrophe-leading-term}
\end{equation}
and also using \eqref{eq:D20c} and \eqref{eq:N20c-final} that
\begin{equation}
\frac{N_\cc^{(2,0)}(t_\cc,x_\cc)D_\cc^{(1,0)}(t_\cc,x_\cc)-N_\cc^{(1,0)}(t_\cc,x_\cc)D_\cc^{(2,0)}(t_\cc,x_\cc)}{D_\cc^{(1,0)}(t_\cc,x_\cc)^2}=-\frac{1}{2t_\cc\nu_\cc}.
\label{eq:catastrophe-next-term}
\end{equation}

\subsubsection{Contour integral representation of the Wronskian $W_\beta[P,Q](\alpha,\beta)$}
The Wronskian $W_\beta[P,Q](\alpha,\beta)$ defined by \eqref{eq:WronskianPQ} would appear to be a double contour integral based on the single-integral representations of $P$ and $Q$ given in \eqref{eq:PQ-define}.  
The goal of this section is to prove the following.
\begin{proposition}
The Wronskian $W_\beta[P,Q](\alpha,\beta)$ has the contour integral representation
\begin{equation}
    W_\beta[P,Q](\alpha,\beta)=2^{-1/6}\pi^{1/2}\ee^{-\ii\alpha^2/2}\int_{\infty\ee^{-7\ii\pi/8}}^{\infty\ee^{-3\ii\pi/8}}\ee^{\ii (\zeta^4-\ii\alpha\zeta^2+\ee^{-3\ii\pi/4}\beta\zeta)}\,\dd\zeta.
\end{equation}
\label{prop:Wronskian}
\end{proposition}
To prove this, first note that both functions $P(\alpha,\beta)$ and $Q(\alpha,\beta)$ defined in \eqref{eq:PQ-define} satisfy the same third-order ordinary differential equation with respect to $\beta$:
\begin{equation}
4\frac{\partial^3P}{\partial\beta^3} -2\alpha\frac{\partial P}{\partial\beta}-\ii\beta P=4\frac{\partial^3Q}{\partial\beta^3}-2\alpha\frac{\partial Q}{\partial\beta}-\ii\beta Q=0.
\label{eq:PQ-ODEs}
\end{equation}
Directly from \eqref{eq:PQ-ODEs}, it then follows that the Wronskian $W_\beta[P,Q](\alpha,\beta)$ itself satisfies a different third-order linear equation:
\begin{equation}
    4\frac{\partial^3 W_\beta}{\partial\beta^3} -2\alpha\frac{\partial W_\beta}{\partial\beta} +\ii\beta W_\beta=0.
\label{eq:W-ODE}
\end{equation}
By Laplace transform methods, this equation has a general solution that is a superposition of a basis of three single contour integral solutions, with coefficients depending on $\alpha$:
\begin{equation}
    W_\beta[P,Q](\alpha,\beta)=\sum_{j=-1}^1 A_j(\alpha)\int_{\Gamma_j}\ee^{\ii(\zeta^4-\ii\alpha\zeta^2+\ee^{-3\ii\pi/4}\beta\zeta)}\,\dd\zeta,
\label{eq:W-general-solution}
\end{equation}
where the $\Gamma_j$ are three infinite integration contours all originating at $\zeta=\infty\ee^{-7\ii\pi/8}$, with $\Gamma_j$ terminating at $\zeta=\infty\ee^{\ii\pi/8}\ee^{\ii j\pi/2}$, $j=-1,0,1$.  Indeed, it is straightforward to verify that the individual terms in \eqref{eq:W-general-solution} satisfy \eqref{eq:W-ODE}.  That \eqref{eq:W-general-solution} is in fact the general solution of \eqref{eq:W-ODE} can be shown by computing a $3\times 3$ Wronskian. 
\begin{lemma}
Let $\mathcal{D}(\alpha,\beta)$ be defined by
\begin{multline}
    \mathcal{D}(\alpha,\beta):=\\\begin{vmatrix}
        \displaystyle \int_{\Gamma_{-1}}\ee^{\ii (\zeta^4-\ii\alpha\zeta^2+\ee^{-3\ii\pi/4}\beta\zeta)}\,\dd\zeta &
        \displaystyle \int_{\Gamma_{0}}\ee^{\ii (\zeta^4-\ii\alpha\zeta^2+\ee^{-3\ii\pi/4}\beta\zeta)}\,\dd\zeta &
        \displaystyle \int_{\Gamma_{1}}\ee^{\ii (\zeta^4-\ii\alpha\zeta^2+\ee^{-3\ii\pi/4}\beta\zeta)}\,\dd\zeta\\
        \displaystyle \partial_\beta\int_{\Gamma_{-1}}\ee^{\ii (\zeta^4-\ii\alpha\zeta^2+\ee^{-3\ii\pi/4}\beta\zeta)}\,\dd\zeta &
        \displaystyle \partial_\beta\int_{\Gamma_{0}}\ee^{\ii (\zeta^4-\ii\alpha\zeta^2+\ee^{-3\ii\pi/4}\beta\zeta)}\,\dd\zeta &
        \displaystyle \partial_\beta\int_{\Gamma_{1}}\ee^{\ii (\zeta^4-\ii\alpha\zeta^2+\ee^{-3\ii\pi/4}\beta\zeta)}\,\dd\zeta\\
        \displaystyle \partial_\beta^2\int_{\Gamma_{-1}}\ee^{\ii (\zeta^4-\ii\alpha\zeta^2+\ee^{-3\ii\pi/4}\beta\zeta)}\,\dd\zeta&
        \displaystyle \partial_\beta^2\int_{\Gamma_{0}}\ee^{\ii (\zeta^4-\ii\alpha\zeta^2+\ee^{-3\ii\pi/4}\beta\zeta)}\,\dd\zeta &
        \displaystyle \partial_\beta^2\int_{\Gamma_{1}}\ee^{\ii (\zeta^4-\ii\alpha\zeta^2+\ee^{-3\ii\pi/4}\beta\zeta)}\,\dd\zeta
    \end{vmatrix}.
\end{multline}
Then, for each $(\alpha,\beta)\in\mathbb{C}^2$, $\mathcal{D}(\alpha,\beta)\neq 0$.
\label{lem:linearly-independent}
\end{lemma}
\begin{proof}
Fix $\alpha\in\mathbb{C}$.  Because there is no second derivative term in \eqref{eq:W-ODE}, it follows from Abel's Theorem that the Wronskian $\beta\mapsto \mathcal{D}(\alpha,\beta)$
is independent of $\beta$.  Hence 
\begin{equation}
    \mathcal{D}(\alpha,\beta)=\lim_{Z\to+\infty}\mathcal{D}(\alpha,\ee^{-\ii\pi/4}Z).
\end{equation}
We will compute this limit explicitly and show that it is nonzero.
Taking the derivatives and evaluating for $\beta=\ee^{-\ii\pi/4}Z$ gives
\begin{multline}
 \ee^{\ii\pi/4}\mathcal{D}(\alpha,\ee^{-\ii\pi/4}Z)=\\\begin{vmatrix}
     \displaystyle\int_{\Gamma_{-1}}\ee^{\ii(\zeta^4-\ii\alpha\zeta^2-Z\zeta)}\,\dd\zeta & 
     \displaystyle\int_{\Gamma_{0}}\ee^{\ii(\zeta^4-\ii\alpha\zeta^2-Z\zeta)}\,\dd\zeta &
     \displaystyle\int_{\Gamma_{1}}\ee^{\ii(\zeta^4-\ii\alpha\zeta^2-Z\zeta)}\,\dd\zeta \\
     \displaystyle\int_{\Gamma_{-1}}\zeta\ee^{\ii(\zeta^4-\ii\alpha\zeta^2-Z\zeta)}\,\dd\zeta & 
     \displaystyle\int_{\Gamma_{0}}\zeta\ee^{\ii(\zeta^4-\ii\alpha\zeta^2-Z\zeta)}\,\dd\zeta &
     \displaystyle\int_{\Gamma_{1}}\zeta\ee^{\ii(\zeta^4-\ii\alpha\zeta^2-Z\zeta)}\,\dd\zeta \\
     \displaystyle\int_{\Gamma_{-1}}\zeta^2\ee^{\ii(\zeta^4-\ii\alpha\zeta^2-Z\zeta)}\,\dd\zeta & 
     \displaystyle\int_{\Gamma_{0}}\zeta^2\ee^{\ii(\zeta^4-\ii\alpha\zeta^2-Z\zeta)}\,\dd\zeta &
     \displaystyle\int_{\Gamma_{1}}\zeta^2\ee^{\ii(\zeta^4-\ii\alpha\zeta^2-Z\zeta)}\,\dd\zeta 
 \end{vmatrix}.      
\end{multline}
Next, rescaling the integration variable by $\zeta=Z^{1/3}\xi$ gives
\begin{multline}
 \ee^{\ii\pi/4}Z^{-2}\mathcal{D}(\alpha,\ee^{-\ii\pi/4}Z)=\\\begin{vmatrix}
     \displaystyle\int_{\Gamma_{-1}}\ee^{\ii Z^{4/3}\varphi(\xi;-\ii Z^{-2/3}\alpha)}\,\dd\xi & 
     \displaystyle\int_{\Gamma_{0}}\ee^{\ii Z^{4/3}\varphi(\xi;-\ii Z^{-2/3}\alpha)}\,\dd\xi &
     \displaystyle\int_{\Gamma_{1}}\ee^{\ii Z^{4/3}\varphi(\xi;-\ii Z^{-2/3}\alpha)}\,\dd\xi \\
     \displaystyle\int_{\Gamma_{-1}}\xi\ee^{\ii Z^{4/3}\varphi(\xi;-\ii Z^{-2/3}\alpha)}\,\dd\xi & 
     \displaystyle\int_{\Gamma_{0}}\xi\ee^{\ii Z^{4/3}\varphi(\xi;-\ii Z^{-2/3}\alpha)}\,\dd\xi &
     \displaystyle\int_{\Gamma_{1}}\xi\ee^{\ii Z^{4/3}\varphi(\xi;-\ii Z^{-2/3}\alpha)}\,\dd\xi \\
     \displaystyle\int_{\Gamma_{-1}}\xi^2\ee^{\ii Z^{4/3}\varphi(\xi;-\ii Z^{-2/3}\alpha)}\,\dd\xi & 
     \displaystyle\int_{\Gamma_{0}}\xi^2\ee^{\ii Z^{4/3}\varphi(\xi;-\ii Z^{-2/3}\alpha)}\,\dd\xi &
     \displaystyle\int_{\Gamma_{1}}\xi^2\ee^{\ii Z^{4/3}\varphi(\xi;-\ii Z^{-2/3}\alpha)}\,\dd\xi 
 \end{vmatrix},     
\end{multline}
where a phase is defined by
\begin{equation}
    \varphi(\xi;\delta):=\xi^4-\xi +\delta\xi^2
\end{equation}
and we used the property that the rescaled integration contours are Cauchy-equivalent to the original ones.  When $Z\gg 1$, the phase is a perturbation of $\varphi(\xi;0)=\xi^4-\xi$ which has three simple critical points at $\xi=\xi_0:=4^{-1/3}$ and $\xi=\xi_{\pm 1}:=\xi_0\ee^{\pm 2\ii\pi/3}$.  Moreover, since the level set $\mathrm{Re}(\ii\varphi(\xi;\delta))=0$ is a perturbation of the limiting level set for $\delta=0$, one can show that the classical steepest descent method applies in which one takes the contour $\Gamma_j$ to pass over the critical point of $\varphi(\xi;\delta)$ closest to $\xi=\xi_j$, $j=-1,0,1$.  See the left-hand panel of Figure~\ref{fig:landscapes}.
\begin{figure}[h]
\begin{center}
\includegraphics[width=0.45\linewidth]{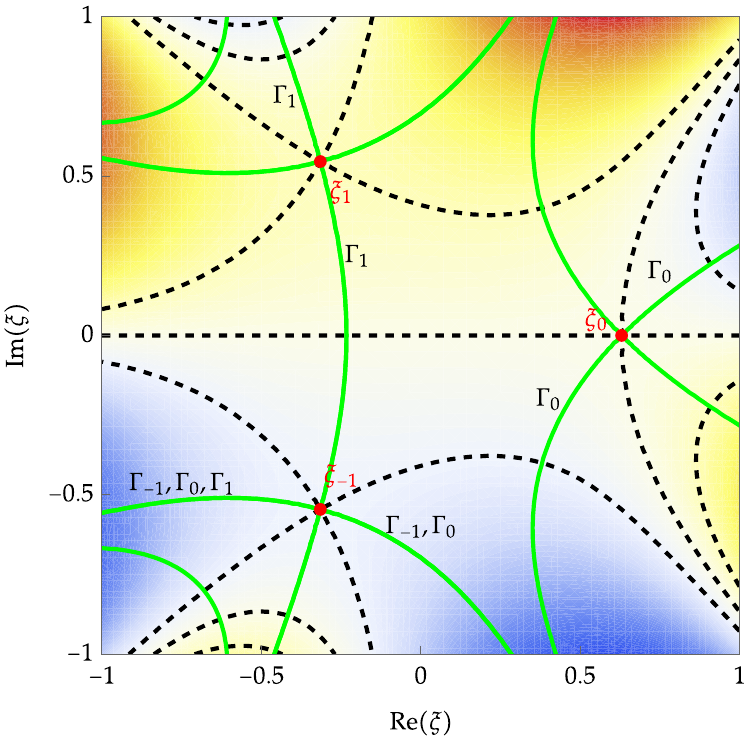}\hfill%
\includegraphics[width=0.45\linewidth]{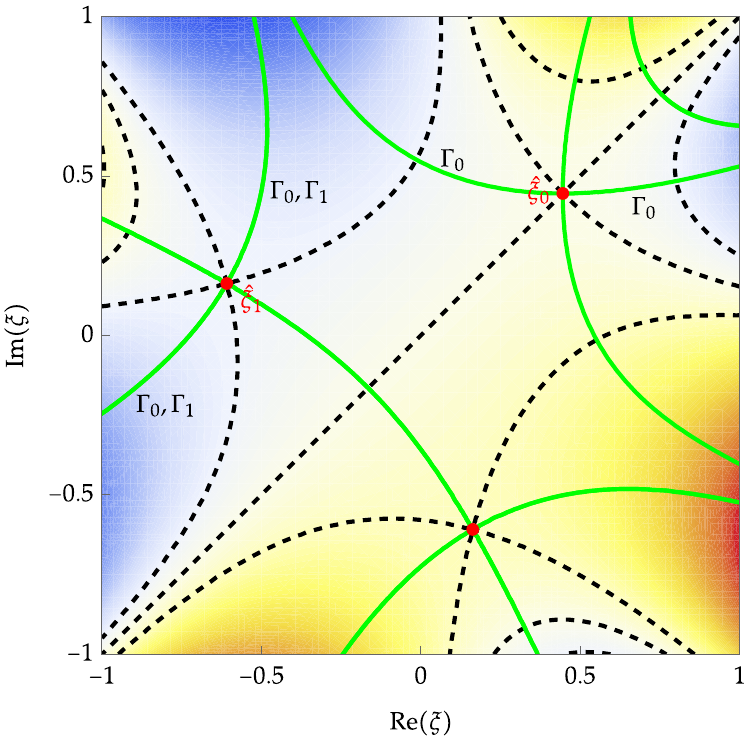}\hfill
\end{center}
\caption{Left:  the landscape of $\mathrm{Re}(\ii\varphi(\xi;0))$ (hot/cool colors for positive/negative values) with the critical points $\xi_j$, $j=-1,0,1$ shown along with the corresponding level curves of $\mathrm{Re}(\ii\varphi(\xi;0))$ and steepest descent curves (level curves of $\mathrm{Im}(\ii\varphi(\xi;0))$).  The contour $\Gamma_{-1}$ is deformed to be close to the limiting steepest descent curve passing over $\xi_{-1}$ from valley to valley; $\Gamma_0$ is deformed into $\Gamma_{-1}$ followed by a curve close to the limiting steepest descent path over $\xi_0$ from valley to valley (and the contribution from $\xi_0$ is exponentially dominant); $\Gamma_1$ is deformed to follow the steepest descent curve uphill to $\xi_{-1}$ and then turn left and continue over the saddle point $\xi_1$, which gives an exponentially dominant contribution.  Right:  the same for $\phi(\xi)$ instead of $\varphi(\xi;0)$.  Here the integration contour $\Gamma_1$ for $Q(0,\ee^{-\ii\pi/4}Z)$ and its $\beta$-derivative is deformed into the steepest descent curve passing over the saddle point $\hat{\xi}_1$ from valley to valley, and the contour $\Gamma_0$ for $P(0,\ee^{-\ii\pi/4}Z)$ and its $\beta$-derivative is deformed into $\Gamma_1$ followed by the steepest descent curve passing over the saddle point $\hat{\xi}_0$, which gives an exponentially dominant contribution.}
\label{fig:landscapes}
\end{figure}
Since $|\varphi''(\xi_j;0)|=12\cdot 4^{-1/3}$, $j=-1,0,1$, and since the corresponding second derivatives at the perturbed critical points differ by $O(Z^{-2/3})$, 
we obtain
\begin{multline}
    \ee^{\ii\pi/4}Z^{-2}\mathcal{D}(\alpha,\ee^{-\ii\pi/4}Z)=\\
    \frac{1}{2 Z^2}\left(\frac{2\pi}{3}\right)^{3/2}\ee^{3\ii\pi/4}\ee^{\ii Z^{4/3}[\varphi_{-1}+\varphi_0+\varphi_{1}]}\left(\begin{vmatrix}1 & 1 & 1\\\xi_{-1} & \xi_0 & \xi_1\\\xi_{-1}^2 & \xi_0^2 & \xi_1^2\end{vmatrix}+O(Z^{-2/3})\right),\quad Z\to+\infty.
\end{multline}
Here, the quantity $\varphi_j$ refers to the value of the phase $\varphi(\xi;\delta)$ at the $\delta$-dependent critical point converging to $\xi_j$ as $Z\to+\infty$.  The sum can be written as a contour integral:
\begin{equation}
    \varphi_{-1}+\varphi_0+\varphi_1=\frac{1}{2\pi\ii}\oint\frac{\varphi(\xi;\delta)\varphi''(\xi;\delta)}{\varphi'(\xi;\delta)}\,\dd\xi
\end{equation}
where the integration contour is a positively-oriented circle of sufficiently large radius to contain all of the critical points.  Evaluating the integral by taking the residue at $\xi=\infty$ gives
\begin{equation}
    \varphi_{-1}+\varphi_0+\varphi_1=-\frac{\delta^2}{2}.
\end{equation}
Then, since $\delta=-\ii Z^{-2/3}\alpha$, we see that $\ii Z^{4/3}[\varphi_{-1}+\varphi_0+\varphi_1]=\frac{1}{2}\ii\alpha^2$.
Finally, since the Vandermonde determinant is nonzero because the three limiting critical points are distinct, we have proved that
\begin{equation}
    \mathcal{D}(\alpha,\beta)=\lim_{Z\to+\infty}\mathcal{D}(\alpha,\ee^{-\ii\pi/4}Z) = \frac{1}{2}\ii\left(\frac{2\pi}{3}\right)^{3/2}\ee^{\ii\alpha^2/2}\begin{vmatrix} 1 & 1 & 1\\\xi_{-1} & \xi_0 & \xi_1\\
    \xi_{-1}^2 & \xi_{0}^2 & \xi_1^2\end{vmatrix}\neq 0,
\end{equation}
so the proof is finished.
\end{proof}

The functions $P(\alpha,\beta)$ and $Q(\alpha,\beta)$ defined in \eqref{eq:PQ-define} also both satisfy the free-particle Schr\"odinger equation:
\begin{equation}
    \ii\frac{\partial P}{\partial\alpha}-\frac{\partial^2P}{\partial\beta^2}=\ii\frac{\partial Q}{\partial\alpha}-\frac{\partial^2Q}{\partial\beta^2}=0.
\end{equation}
Using these and also \eqref{eq:PQ-ODEs} to differentiate $W_\beta[P,Q](\alpha,\beta)$ gives a time-dependent Schr\"odin\-ger equation:
\begin{equation}
    \ii\frac{\partial W_\beta}{\partial\alpha} +\frac{\partial^2W_\beta}{\partial\beta^2} -\alpha W_\beta = 0.
\end{equation}
Substituting into this equation the representation \eqref{eq:W-general-solution} and using the fact that
\begin{equation}
    \left(\ii\partial_\alpha+\partial_\beta^2\right)\int_{\Gamma_j}\ee^{\ii (\zeta^4-\ii\alpha\zeta^2+\ee^{-3\ii\pi/4}\beta\zeta)}\,\dd\zeta=0,
\end{equation}
we obtain the identity
\begin{equation}
    \sum_{j=-1}^1\left[\ii A_j'(\alpha)-\alpha A_j(\alpha)\right]\int_{\Gamma_j}\ee^{\ii (\zeta^4-\ii\alpha\zeta^2+\ee^{-3\ii\pi/4}\beta\zeta)}\,\dd\zeta = 0.
\end{equation}
By Lemma~\ref{lem:linearly-independent}, the three contour integrals are linearly independent functions of $\beta$, so it follows that each coefficient $A_j(\alpha)$ satisfies
\begin{equation}
    A_j(\alpha)= A_j(0)\ee^{-\ii\alpha^2/2}.
\end{equation}

The next task is to determine the values of the three constants $A_{-1}(0)$, $A_0(0)$, and $A_1(0)$.  For this purpose, we may set $\alpha=0$ in \eqref{eq:W-general-solution}, and then using the same arguments as in the proof of Lemma~\ref{lem:linearly-independent} one gets that
\begin{equation}
    W_\beta[P,Q](0,\ee^{-\ii\pi/4}Z)=\sum_{j=-1}^1 A_j(0)\frac{2^{-1/6}3^{-1/2}\pi^{1/2}}{Z^{1/3}}\ee^{\ii\theta_j}\ee^{\ii Z^{4/3}\varphi(\xi_j;0)}(1+O(Z^{-4/3})),\quad Z\to +\infty,
\end{equation}
where $\theta_{-1}=-\pi/12$, $\theta_0=\pi/4$, and $\theta_1=7\pi/12$ are the steepest-descent tangent angles with which $\Gamma_{-1}$, $\Gamma_0$, and $\Gamma_1$ traverse the critical points $\xi_{-1}$, $\xi_0$, and $\xi_1$ respectively.  
Note that 
\begin{equation}
    \ii\varphi(\xi_{\pm 1};0)=\pm 2^{-11/3}3^{3/2} + 2^{-11/3}3\ii,\quad\ii\varphi(\xi_0;0)=-2^{-8/3}3\ii,
\end{equation}
so that $\ee^{\ii Z^{4/3}\varphi(\xi_{-1};0)}$, $\ee^{\ii Z^{4/3}\varphi(\xi_0;0)}$, and $\ee^{\ii Z^{4/3}\varphi(\xi_1;0)}$ are respectively decaying, oscillatory, and growing exponentials as $Z\to+\infty$.
Carrying out a similar asymptotic study of $P(0,\ee^{-\ii\pi/4}Z)$, $Q(0,\ee^{-\ii\pi/4}Z)$, $\partial_\beta P(0,\ee^{-\ii\pi/4}Z)$, and $\partial_\beta Q(0,\ee^{-\ii\pi/4}Z)$, we begin by rescaling the integration variable in \eqref{eq:PQ-define} by $\zeta=Z^{1/3}\xi$ to get
\begin{equation}
    \begin{split}
        P(0,\ee^{-\ii\pi/4} Z)&=Z^{1/3}\int_{\Gamma_0}\ee^{\ii Z^{4/3}\phi(\xi)}\,\dd\xi\\
        \partial_\beta P(0,\ee^{-\ii\pi/4} Z)&=\ii Z^{2/3}\int_{\Gamma_0}\xi\ee^{\ii Z^{4/3}\phi(\xi)}\,\dd\xi\\
        Q(0,\ee^{-\ii\pi/4}Z)&=Z^{1/3}\int_{\Gamma_1}\ee^{\ii Z^{4/3}\phi(\xi)}\,\dd\xi\\
        \partial_\beta Q(0,\ee^{-\ii\pi/4}Z)&=\ii Z^{2/3}\int_{\Gamma_1}\xi\ee^{\ii Z^{4/3}\phi(\xi)}\,\dd\xi,
    \end{split}
\end{equation}
where $\phi(\xi):=\xi^4+\ee^{-\ii\pi/4}\xi$.  Referring to the right-hand panel of Figure~\ref{fig:landscapes}, we obtain that the dominant contributions to $P$ and $\partial_\beta P$ comes from the critical point $\xi=\hat{\xi}_0:=\ee^{\ii\pi/4}4^{-1/3}$ traversed with tangent angle $\theta=0$, while the dominant contributions to $Q$ and $\partial_\beta Q$ come instead from the critical point $\xi=\hat{\xi}_1:=\ee^{11\ii\pi/12}4^{-1/3}$ traversed with tangent angle $\theta=\frac{1}{3}\pi$, and therefore as $Z\to+\infty$,
\begin{equation}
    \begin{split}
        P(0,\ee^{-\ii\pi/4}Z)&=Z^{1/3}\sqrt{\frac{2\pi}{12\cdot 4^{-2/3} Z^{4/3}}}\ee^{\ii Z^{4/3}\phi(\hat{\xi}_0)}(1+O(Z^{-4/3})),\\
        \partial_\beta P(0,\ee^{-\ii\pi/4}Z)&=\ii Z^{2/3}\sqrt{\frac{2\pi}{12\cdot 4^{-2/3}Z^{4/3}}}\hat{\xi}_0\ee^{\ii Z^{4/3}\phi(\hat{\xi}_0)}(1+O(Z^{-4/3})),\\
        Q(0,\ee^{-\ii\pi/4}Z)&=Z^{1/3}\sqrt{\frac{2\pi}{12\cdot 4^{-2/3}Z^{4/3}}}\ee^{\ii\pi/3}\ee^{\ii Z^{4/3}\phi(\hat{\xi}_1)}(1+O(Z^{-4/3})),\\
        \partial_\beta Q(0,\ee^{-\ii\pi/4}Z)&=\ii Z^{2/3}\sqrt{\frac{2\pi}{12\cdot 4^{-2/3}Z^{4/3}}}\ee^{\ii\pi/3}\hat{\xi}_1\ee^{\ii Z^{4/3}\phi(\hat{\xi}_1)}(1+O(Z^{-4/3})).
    \end{split}
\end{equation}
Using these in the definition \eqref{eq:WronskianPQ} along with the observation that
\begin{equation}
    \ii (\phi(\hat{\xi}_0)+\phi(\hat{\xi}_1))= -2^{-11/3}3^{3/2}+2^{-11/3}3\ii,
\end{equation}
gives that as $Z\to+\infty$,
\begin{equation}
    W_\beta[P,Q](0,\ee^{-\ii\pi/4}Z)=\frac{2^{-1/3}3^{-1/2}\pi}{Z^{1/3}}\ee^{-\ii\pi/12}\ee^{Z^{4/3}(-2^{-11/3}3^{3/2}+2^{-11/3}3\ii)}(1+O(Z^{-4/3})).
\label{eq:WPQ-expand}
\end{equation}
This is a solution of \eqref{eq:W-ODE} that is decaying exponentially as $Z\to+\infty$, so we are forced to take $A_0(0)=A_1(0)=0$ in \eqref{eq:W-general-solution}, and then to match the asymptotic between \eqref{eq:W-general-solution} and \eqref{eq:WPQ-expand} we must choose
\begin{equation}
    A_{-1}(0)=2^{-1/6}\pi^{1/2}.
\end{equation}
This completes the proof of Proposition~\ref{prop:Wronskian}.

Using Proposition~\ref{prop:Wronskian}, we obtain 
\begin{equation}
\begin{split}
    W_\beta[P,Q](-\ii Y,\ee^{-\ii\pi/4}Z)&=2^{-1/6}\pi^{1/2}\ee^{\ii Y^2/2}\int_{\infty\ee^{-7\ii\pi/8}}^{\infty\ee^{-3\ii\pi/8}}\ee^{\ii(\zeta^4 - Y\zeta^2-Z\zeta)}\,\dd\zeta \\ &= 2^{-1/6}\pi^{1/2}\ee^{\ii Y^2/2}\tau(Y,Z),
\end{split}
\end{equation}
where on the second line we recall the definition \eqref{eq:tau-define}.

\subsubsection{Completion of the proof of Theorem~\ref{thm:BO-catastrophe}}
Recall the assumption that $\tau(T,X)\neq 0$.  It then follows from \eqref{eq:YZ-perturb} that also $\tau(Y,Z)\neq 0$ for $\epsilon>0$ sufficiently small.  Dividing through by $W_\beta[P,Q](-\ii Y,\ee^{-\ii\pi/4}Z)\neq 0$ and $D_\cc^{(1,0)}(t_\cc,x_\cc)\neq 0$ along with the numerical factors $\ee^{-\ii\pi/4}(4\epsilon)^{3/4}$ in \eqref{eq:Numerator-catastrophe} and \eqref{eq:Denominator-catastrophe}, we obtain respectively
\begin{multline}
    \frac{N(t,x;\epsilon)}{E_\cc(t,x;\epsilon)\ee^{-\ii\pi/4}(4\epsilon)^{3/4}D_\cc^{(1,0)}(t_\cc,x_\cc)W_\beta[P,Q](-\ii Y,\ee^{-\ii\pi/4} Z)} = \\
    \frac{N_\cc^{(1,0)}(t_\cc,x_\cc)}{D_\cc^{(1,0)}(t_\cc,x_\cc)} -(4\epsilon)^{1/4}\ii\frac{N_\cc^{(2,0)}(t_\cc,x_\cc)}{D_\cc^{(1,0)}(t_\cc,x_\cc)}\frac{\partial}{\partial Z}\log(\tau(Y,Z)) + O(\epsilon^{1/2})
\end{multline}
and
\begin{multline}
    \frac{D(t,x;\epsilon)}{E_\cc(t,x;\epsilon)\ee^{-\ii\pi/4}(4\epsilon)^{3/4}D_\cc^{(1,0)}(t_\cc,x_\cc)W_\beta[P,Q](-\ii Y,\ee^{-\ii\pi/4} Z)} = \\
    1-(4\epsilon)^{1/4}\ii\frac{D_\cc^{(2,0)}(t_\cc,x_\cc)}{D_\cc^{(1,0)}(t_\cc,x_\cc)}\frac{\partial}{\partial Z}\log(\tau(Y,Z)) + O(\epsilon^{1/2}).
\label{eq:Denominator-catastrophe-renormalize}
\end{multline}
In both cases, we may use \eqref{eq:YZ-perturb} again to replace $(Y,Z)$ with $(T,X)$ on the right-hand side at the cost of an error absorbable into the existing $O(\epsilon^{1/2})$ terms.  Expanding the reciprocal of \eqref{eq:Denominator-catastrophe-renormalize} then shows that
\begin{multline}
    \frac{N(t,x;\epsilon)}{D(t,x;\epsilon)}=\frac{N_\cc^{(1,0)}(t_\cc,x_\cc)}{D_\cc^{(1,0)}(t_\cc,x_\cc)} \\+ (4\epsilon)^{1/4}\ii\frac{N_\cc^{(1,0)}(t_\cc,x_\cc)D_\cc^{(2,0)}(t_\cc,x_\cc)-N_\cc^{(2,0)}(t_\cc,x_\cc)D_\cc^{(1,0)}(t_\cc,x_\cc)}{D_\cc^{(1,0)}(t_\cc,x_\cc)^2}\frac{\partial}{\partial X}\log(\tau(T,X)) + O(\epsilon^{1/2}).
\end{multline}
Using \eqref{eq:catastrophe-leading-term} and \eqref{eq:catastrophe-next-term} then gives
\begin{equation}
    \frac{N(t,x;\epsilon)}{D(t,x;\epsilon)}=\frac{1}{2t_\cc}\left(x_\cc-2y_\cc-\sum_{m=2}^Nz_m(t_\cc,x_\cc)+\sum_{j=1}^Np_j\right) +\frac{(4\epsilon)^{1/4}}{2t_\cc\nu_\cc}\ii \frac{\partial}{\partial X}\log(\tau(T,X))+O(\epsilon^{1/2}).
\end{equation}
Then using Theorem~\ref{thm:IVP-solution},
\begin{equation}
\begin{split}
    u(t,x;\epsilon)&=2\mathrm{Re}\left(\frac{N(t,x;\epsilon)}{D(t,x;\epsilon)}\right)\\
    &=\frac{1}{2t_\cc}\left(2x_\cc-4y_\cc-\sum_{m=2}^N(z_m(t_\cc,x_\cc)+z_m(t_\cc,x_\cc)^*) +\sum_{j=1}^N(p_j+p_j^*)\right) \\
    &\qquad\qquad{}-\frac{(4\epsilon)^{1/4}}{t_\cc\nu_\cc}\mathrm{Im}\left(\frac{\partial}{\partial X}\log(\tau(T,X))\right) + O(\epsilon^{1/2}).
\end{split}
\end{equation}
Since according to the discussion in Section~\ref{sec:catastrophe-contours} the real critical points of $h(\diamond;t_\cc,x_\cc)$ are exhausted by the triple critical point $y_\cc$ while the critical points $z_2(t_\cc,x_\cc),\dots,z_N(t_\cc,x_\cc)$ in the upper half-plane are all simple, by applying Lemma~\ref{lem:PoleSPRelation} with $\mathcal{Y}=\{y_\cc,y_\cc,y_\cc\}$ and $\mathcal{Z}=\{z_2(t_\cc,x_\cc),\dots,z_N(t_\cc,x_\cc)\}$ and using $x_\cc-y_\cc=2t_\cc u_0(y_\cc)=2t_\cc u_\cc$ we obtain
\begin{equation}
    u(t,x;\epsilon)=u_\cc -\frac{(4\epsilon)^{1/4}}{t_\cc\nu_\cc}\mathrm{Im}\left(\frac{\partial}{\partial X}\log(\tau(T,X))\right)+O(\epsilon^{1/2}),\quad\epsilon\to 0.
\end{equation}
Recalling the definition of $\nu_\cc>0$ in \eqref{eq:nu-cc} and comparing with \eqref{eq:U-profile} completes the proof.

\bibliographystyle{siamplain}

\end{document}